\documentclass[11pt, a4paper]{article}
\usepackage{times}
\usepackage[dvips, hyperindex]{hyperref}
\usepackage[british,german]{babel}
\usepackage{enumerate, longtable}
\usepackage{amsmath, amscd, amsfonts, amssymb, latexsym, theorem, comment, graphicx}
\usepackage[all]{xypic}
\usepackage[T1]{fontenc}
\usepackage[latin1]{inputenc}

{\theorembodyfont{\itshape}   \newtheorem{thm}{Theorem}[section]}
{\theorembodyfont{\itshape}   }
{\theorembodyfont{\itshape}   \newtheorem{defi}[thm]{Definition}}
{\theorembodyfont{\itshape}   }
{\theorembodyfont{\itshape}   }
{\theorembodyfont{\itshape}   \newtheorem{prop}[thm]{Proposition}}
{\theorembodyfont{\itshape}   \newtheorem{cor}[thm]{Corollary}}
{\theorembodyfont{\itshape}   \newtheorem{question}[thm]{Question}}

\newcommand{\CC}{\mathbb{C}}
\newcommand{\FF}{\mathbb{F}}
\newcommand{\NN}{\mathbb{N}}
\newcommand{\QQ}{\mathbb{Q}}

\newcommand{\TT}{\mathbb{T}}
\newcommand{\ZZ}{\mathbb{Z}}

\newcommand{\fp}{\mathfrak{p}}
\newcommand{\fm}{\mathfrak{m}}

\newcommand{\cO}{\mathcal{O}}

\newcommand{\Qbar}{\overline{\QQ}}
\newcommand{\Zbar}{\overline{\ZZ}}
\newcommand{\Fbar}{\overline{\FF}}
\newcommand{\Tbar}{\overline{\TT}}

\newcommand{\That}{\widehat{\TT}}

\DeclareMathOperator{\Frac}{Frac}
\DeclareMathOperator{\Gal}{Gal}
\DeclareMathOperator{\Spec}{Spec}
\DeclareMathOperator{\MaxSpec}{MaxSpec}
\DeclareMathOperator{\MinSpec}{MinSpec}

\newcommand{\End}{{\rm End}}
\newcommand{\Aut}{{\rm Aut}}
\newcommand{\Hom}{{\rm Hom}}
\newcommand{\GL}{\mathrm{GL}}
\newcommand{\SL}{\mathrm{SL}}

\newcommand{\const}{\mathrm{const}}
\newcommand{\new}{\mathrm{new}}

\newcommand{\pf}{{\bf Proof. }}
\newcommand{\qed}{\hspace* {.5cm} \hfill $\Box$}

\newcommand{\mplot}[1]{\parbox[t]{7.8cm}{\stepcounter{figure} \hspace*{1cm} Figure \thefigure \\[-.9cm]
\includegraphics[width=5.5cm, angle=270]{#1}}}

\newcommand{\mplotr}[2]{\parbox[t]{7.8cm}{\refstepcounter{figure}\label{#2}\hspace*{1cm}
Figure \thefigure  \\[-.9cm]
\includegraphics[width=5.5cm, angle=270]{#1}}}

\begin{document}

\selectlanguage{british}

\title{A Computational Study of the\\ Asymptotic Behaviour of Coefficient
Fields of Modular Forms}
\author{Marcel Mohyla and Gabor Wiese\footnote{Authors' address:
Universität Duisburg-Essen, Institut für
Experimentelle Mathematik, Ellernstraße 29, 45326 Essen, Germany}}
\maketitle

\begin{abstract}
The article motivates, presents and describes large computer calculations
concerning the asymptotic behaviour of arithmetic properties of coefficient
fields of modular forms. The observations suggest certain patterns,
which deserve further study.

\end{abstract}

\section{Introduction}

A recent breakthrough in Arithmetic Geometry is the proof of the
Sato-Tate conjecture by Barnet-Lamb, Clozel, Geraghty, Harris, Shepherd-Barron and Taylor
(\cite{BLGHT}, \cite{CHT}, \cite{HSHT}, \cite{T}).
It states that the normalised Hecke eigenvalues $\frac{a_p(f)}{2p^{(k-1)/2}}$ on a holomorphic
newform~$f$ of weight $k \ge 2$ (and trivial Dirichlet character\footnote{We only make
this assumption for the sake of simplicity of the exposition.})
are equidistributed with respect to a certain measure
(the so-called Sato-Tate measure), when~$p$ runs through the set of prime numbers.
Let us use the name {\em horizontal} Sato-Tate for this situation\footnote{This name
was suggested by G.\ Chênevert.}.

The reversed situation, to be referred to as {\em vertical} horizontal Sato-Tate, was
successfully treated by Serre in~\cite{repas}. He fixes a prime~$p$ and allows
any sequence of positive integers $(N_n,k_n)$ with even $k_n$ and $p \nmid N_n$
such that $N_n+k_n$ tends to infinity and proves that $\frac{a_p(f)}{2p^{(k-1)/2}}$ is equidistributed with
respect to a certain measure depending on~$p$ (which is related to the Sato-Tate measure),
when $f$ runs through all the newforms in any of the spaces of cusp forms
of level $\Gamma_0(N_n)$ and weight~$k_n$.
A corollary is that for fixed positive even weight~$k$, the set
\begin{equation}\label{eqserre}
 \{[\QQ_f:\QQ] \;|\; f \textnormal{ newform of level } \Gamma_0(N_n) \textnormal{ and weight } k\}
\end{equation}
is unbounded for any sequence $N_n$ tending to infinity. Here, $\QQ_f$ denotes the number
field obtained from~$\QQ$ by adjoining all Hecke eigenvalues on~$f$.

In this article we perform a first computational study towards a (weak) arithmetic analogue of
vertical Sato-Tate, where the name {\em arithmetic} refers to taking
a finite place of~$\QQ$, as opposed to the infinite place used in usual
Sato-Tate (the assertion of usual Sato-Tate concerns the Hecke eigenvalues as real numbers).
An arithmetic analog of horizontal Sato-Tate is Chebotarev's density
theorem. Consider, for example, a normalised cuspidal Hecke eigenform~$f$
with attached Galois representation
$\rho_f: \Gal(\Qbar/\QQ)  \to \GL_2(\ZZ_\ell)$.\footnote{Again, it is only for simplicity
of the exposition that we are taking $\ZZ_\ell$ instead of $\Qbar_\ell$.}
Fix some $x \in \ZZ_\ell$ and let $n \in \NN$.
Let $G$ be the image of the composite representation
$\Gal(\Qbar/\QQ) \xrightarrow{\rho_f} \GL_2(\ZZ_\ell) \twoheadrightarrow
\GL_2(\ZZ/\ell^n\ZZ)$ and let $d(x)$ be the number of elements in~$G$ with trace
equal to~$x$ modulo~$\ell^n$.
Then the density of the set $\{ p \,|\, |a_p - x|_\ell \le \ell^{-n}\}$ is equal to
$\frac{d(x)}{|G|}$ by Chebotarev's density theorem; hence, the situation is completely clear.
Whereas at the infinite place horizontal Sato-Tate seems to be more difficult
than vertical Sato-Tate, the situation appears to be reversed for arithmetic analogs.
We are not going to propose such an analogue. But, we are going to study related questions
by means of computer calculations. For instance, as a motivation let us consider the set
\begin{equation}\label{eqzwei}
\{[\FF_{\ell,f}:\FF_\ell] \,|\, f \in S_{k}(N_n;\Fbar_\ell) \textnormal{ normalised Hecke eigenform} \}
\end{equation}
in analogy to Equation~\ref{eqserre}.
Here, $S_{k}(N_n;\Fbar_\ell)$ denotes the $\Fbar_\ell$-vector space of
cuspidal modular forms over~$\Fbar_\ell$ (see Section~\ref{sec:background} for definitions)
and $\FF_{\ell,f}$ is defined by adjoining to~$\FF_\ell$ all Hecke eigenvalues on~$f$.
It is easy to construct sequences $(N_n,k_n)$ for which the set in question is
infinite (see e.g.\ \cite{DiWi} and \cite{W}),
but it does not seem simple to obtain all natural numbers as degrees.
Most importantly, it seems to be unknown whether this set is infinite
when $(N_n)$ is the sequence of prime numbers, $k=2$ and $\ell > 2$.

Concerning properties of modular forms in positive characteristic, there is other,
much more subtle information than just the degree of $\FF_{\ell,f}$ to be studied,
e.g.\ congruences between modular forms.
In order to take the full information into account, in this article we examine
the properties of the $\FF_\ell$-Hecke algebras $\Tbar_k(N_n)$ on $S_{k}(N_n;\Fbar_\ell)$
asymptotically for fixed weight~$k$ (mostly~$2$) and running level~$N_n$ (mostly the set
of prime numbers) by means of experimentation.
More precisely, we investigate three quantities:

\begin{enumerate}[(a)]
\item\label{part:a} The deviation of $\Tbar_k(N_n)$ from
being semisimple. In Section~\ref{sec:ss}, we include a proposition relating nonsemisimplicity
to congruences, ramification and certain indices. Our experiments suggest that for odd
primes~$\ell$, the Hecke algebra $\Tbar_k(N_n)$ tends to be close to semisimple,
whereas the situation seems to be completely different for~$p=2$ (see Section~\ref{sec:a}).

\item\label{part:b} The average residue degree of $\Tbar_k(N_n)$. That is the
arithmetic mean of all the $\FF_{\ell,f}$ for f in $S_{k}(N_n;\Fbar_\ell)$.
Our computations (see Section~\ref{sec:b})
strongly suggest that this quantity is unbounded. More precisely, we seem to observe
a certain asymptotic behaviour, which we formulate as a question.

\item\label{part:c} The maximum residue degree of $\Tbar_k(N_n)$. That is the maximum of all
$\FF_{\ell,f}$ for f in $S_{k}(N_n;\Fbar_\ell)$. Our experiments suggest that
this quantity is 'asymptotically' proportional to the dimension of $S_{k}(N_n;\Fbar_\ell)$.

\end{enumerate}

In Section~\ref{sec:observations} we describe our computations and derive certain
questions from our observations. However, we do not attempt to propose
any heuristic explanations in this article. This will have to be the subject of
subsequent studies, building on refined and extended computations.

We see (\ref{part:b}) and (\ref{part:c}) as strong evidence for the infinity of the set
in Equation~\ref{eqzwei}, when $N_n$ runs through the primes.
Generally speaking, there appears to be some regularity in the
otherwise quite erratic behaviour of the examined quantities, lending some support to
the hope of finding a formulation of an arithmetic analogue of vertical Sato-Tate.

\subsection*{Acknowledgements}

G.~W.\ acknowledges partial support by the European Research Training Network
{\em Galois Theory and Explicit Methods} MRTN-CT-2006-035495
and by the Sonderforschungsbereich Transregio 45
{\em Periods, moduli spaces and arithmetic of algebraic varieties}
of the Deutsche Forschungsgemeinschaft.

\section{Background and notation}\label{sec:background}

We start by introducing the necessary notation and explaining the background.
For facts on modular forms, we refer to~\cite{DiamondIm}.
Let us fix an interger~$k$ and a congruence subgroup $\Gamma \subseteq \SL_2(\ZZ)$.
Denote by $S_k(\Gamma)$ the complex vector space of holomorphic cusp forms of weight~$k$
for~$\Gamma$. Define $\TT := \TT_k(\Gamma)$ to be the $\ZZ$-Hecke algebra of weight~$k$ for~$\Gamma$,
i.e.\ the subring of $\End_\CC(S_k(\Gamma))$ spanned
by all Hecke operators~$T_n$ for $n \in \NN$.
If $\Gamma = \Gamma_0(N)$ we simply write $\TT_k(N)$. We use similar notation
in other contexts, too.
It is an important theorem that $\TT$ is free of finite rank as a $\ZZ$-module,
hence has Krull dimension one as a ring, and that the map
\begin{equation}\label{eqq}
\Hom_\ZZ(\TT,\CC) \to S_k(\Gamma), \;\;\; \phi \mapsto \sum_{n=1}^\infty \phi(T_n) q^n
\end{equation}
with $q=q(z)=e^{2 \pi i z}$
defines an isomorphism of $\CC$-vector spaces, which is compatible with the natural Hecke
action. For any ring~$A$, define $S_k(\Gamma; A) := \Hom_\ZZ(\TT,A)$
equipped with the natural Hecke action (i.e.\ $\TT$-action), so that we have
$S_k(\Gamma) \cong S_k(\Gamma;\CC)$. We think of elements in $S_k(\Gamma;A)$
in terms of formal $q$-expansions, i.e.\ as formal power series in $A[[q]]$,
by an analog of Eq.~\ref{eqq}.
Note that normalised Hecke eigenforms, i.e.\ those
$f = \sum_{n=1}^\infty a_n q^n \in S_k(\Gamma;A)$ that satisfy $a_1=1$ and
$T_n f = a_n f$, precisely correspond to ring homomorphisms $\phi:\TT \to A$
with $\phi(T_n)=a_n$.
When $A$ is an integral domain, a normalised eigenfunction
gives rise to a prime ideal $\fp$ of~$\TT$, namely the kernel of~$\phi$.
We may think of $\TT/\fp$ as the smallest subring of~$A$ generated by the $a_n$
for $n \in \NN$: the {\em coefficient ring of~$f$ in~$A$}.
Note that $\Aut(A)$ acts on $S_k(\Gamma; A)$ by composing $\phi:\TT \to A$ with
$\sigma \in \Aut(A)$. Obviously, this action does not change the ideal~$\fp$
corresponding to an eigenform.

We fix a prime number~$p$. We put
$\That := \That_k(\Gamma) := \TT_k(\Gamma) \otimes_\ZZ {\ZZ_p}$,
$\That_\eta := \TT_k(\Gamma) \otimes_\ZZ {\QQ_p} = \That \otimes_{\ZZ_p} \QQ_p$
and
$\Tbar := \Tbar_k(\Gamma) := \TT_k(\Gamma) \otimes_\ZZ {\FF_p} = \That \otimes_{\ZZ_p} \FF_p$.
Note the isomorphisms
$S_k(\Gamma;\Zbar_p) \cong \Hom_{\ZZ_p}(\That,\Zbar_p)$,
$S_k(\Gamma;\Qbar_p) \cong \Hom_{\ZZ_p}(\That,\Qbar_p)$
and
$S_k(\Gamma;\Fbar_p) \cong \Hom_{\ZZ_p}(\That,\Fbar_p)\cong \Hom_{\FF_p}(\Tbar,\Fbar_p)$.

The $\Gal(\Qbar_p/\QQ_p)$-conjugacy classes of normalised eigenforms
in $S_k(\Gamma;\Qbar_p)$ (by which we mean the classes for the $\Gal(\Qbar_p/\QQ_p)$-action
described above)
are in bijection with the prime (and automatically
maximal) ideals of~$\That_\eta$ and also in bijection with the minimal prime ideals of~$\That$,
whose set is denoted by $\MinSpec(\That)$.
The second bijection is explicitly given by taking preimages for the injection
$\That \hookrightarrow \That_\eta$.
Note that $\That$ has Krull dimension one, meaning that any prime ideal
is either minimal (i.e.\ not containing any smaller prime ideal) or maximal.
Moreover, the $\Gal(\Fbar_p/\FF_p)$-conjugacy classes
of normalised eigenforms in $S_k(\Gamma;\Fbar_p)$ are in bijection
with $\Spec(\Tbar) = \MaxSpec(\Tbar)$.
Furthermore, $\Spec(\Tbar)$ is in natural bijection with $\MaxSpec(\That)$ under
taking preimages for the natural projection $\That \twoheadrightarrow \Tbar$.
By a result in commutative algebra, we have direct product decompositions
$$ \That = \prod_{\fm \in \MaxSpec(\That)} \That_\fm, \;\;\;
\Tbar = \prod_{\fm \in \MaxSpec(\That)} \Tbar_\fm \;\;\textrm{ and }\;\;
\That_\eta = \prod_{\fp \in \Spec(\That_\eta)} \That_{\eta,\fp},$$
where the factors are the localisations at the prime ideals
indicated by the subscripts.

\begin{defi}
We say that two $\fp_1, \fp_2 \in \MinSpec(\That)$ are {\em congruent} if they
lie in the same maximal ideal $\fm \in \MaxSpec(\That)$.
For $\fp \in \MinSpec(\That)$, we call $\That/\fp$ the {\em local coefficient ring} and
$L_\fp := \Frac(\That/\fp)$ the {\em local coefficient field}.
We say that $\fp \in \MinSpec(\That)$ is {\em ramified} if
$L_\fp$ is a ramified extension of~$\QQ_p$.
We denote by $i_\fp$ the index of $\That/\fp$ in the ring of integers of~$L_\fp$.
The residue field $\That/\fm \cong \Tbar/\fm$ will be denoted by $\FF_\fm$
and will be called the {\em residual coefficient field}.
\end{defi}

We now establish the connection with the usual understanding of the terms
in the definition.

Let $\Zbar \subset \Qbar \subset \CC$ be the algebraic integers and the algebraic
numbers, respectively. As $\TT$ is of finite $\ZZ$-rank, the set of normalised
eigenforms in $S_k(\Gamma)$ is the same as the set of normalised eigenforms
in $S_k(\Gamma; \Zbar)$.
Fix homorphisms
$$ \xymatrix@=1cm{
\Zbar \ar@{^(->}[r]^{\iota}  \ar@{->>}[dr]_\pi &
\Zbar_p \ar@{->>}[d]^\pi \\
& \Fbar_p} 
\hspace*{.5cm}\textnormal{ giving rise to }\hspace*{.5cm}
\xymatrix@=1cm{
S_k(\Gamma;\Zbar) \ar@{^(->}[r]^{\iota}  \ar@{->>}[dr]_\pi &
S_k(\Gamma;\Zbar_p) \ar@{->>}[d]^\pi \\
& S_k(\Gamma;\Fbar_p)}.$$
From this perspective, a holomorphic normalised Hecke eigenform
$f= \sum_{n=1}^\infty a_n q^n \in S_k(\Gamma)$
gives rise to an eigenform in $S_k(\Gamma; \Fbar_p)$,
called the {\em reduction of~$f$ modulo~$p$},
whose formal $q$-expansion is $\pi(f) := \sum_{n=1}^\infty \pi(a_n) q^n \in \Fbar_p[[q]]$.
The reduction corresponds to $\fm \in \MaxSpec(\That)$ and to $\fm \in \Spec(\Tbar)$
(we use the same symbol due to the natural bijection between the two sets).
The form~$f$ also gives rise to an eigenform in $S_k(\Gamma; \Qbar_p)$, which
corresponds to $\fp_f \in \MinSpec(\That)$ and to $\fp_f \in \Spec(\That_\eta)$
(the same symbol is used again due to the natural bijection explained above).

Let $g=\sum_{n=1}^\infty b_n q^n$ be another holomorphic normalised Hecke eigenform.
If $\pi(f) = \pi(g)$, then clearly $\fp_f \subset \fm$ and $\fp_g \subset \fm$, i.e.\
$\fp_f$ and $\fp_g$ are congruent.
Conversely, let $\fp_f,\fp_g \in \MinSpec(\That)$ such that $\fp_f \subset \fm$
and $\fp_g \subset \fm$ for some $\fm \in \MaxSpec(\That)$, so that $\fp_f$
and $\fp_g$ are congruent. The ideals $\fp_f$ and $\fp_g$ correspond to
$\Gal(\Qbar_p/\QQ_p)$-conjugacy classes in $S_k(\Gamma;\Qbar_p)$ and we can
choose $f,g \in S_k(\Gamma; \Zbar_p)$ corresponding to $\fp_f$ and
$\fp_g$ with $\pi(f) = \pi(g)$.
We illustrate the situation by the diagram
$$ \xymatrix@=1cm{
&& \That/\fp_f \ar@{->>}[drr] \ar@{^(->}[d]  \\
\That \ar@{->>}[urr] \ar@{->>}[drr] \ar@/^/[rr]^{f} \ar@/_/[rr]_{g}&&
 \Zbar_p\ar@{->>}[r]^\pi& \Fbar_p & \Tbar/\fm = \FF_\fm. \ar@{_(->}[l]\\
&& \That/\fp_g \ar@{->>}[urr] \ar@{_(->}[u]  }$$
Note that $f,g$ can already be found in $S_k(\Gamma;\Zbar) \subset S_k(\Gamma)$.
This justifies our usage of the term {\em congruence}.

Moreover, the local coefficient ring $\That/\fp$ can be identified
with $\ZZ_{p,f} := \ZZ_p[\iota(a_n) | n \in \NN]$ and its fraction field~$L_\fp$ with
$\QQ_{p,f} := \QQ_p(\iota(a_n) | n \in \NN)$, whence $i_\fp$ is the index
of $\ZZ_{p,f}$ in its normalisation.
Furthermore, the residual coefficient field, i.e.\ $\FF_\fm = \Tbar/\fm$,
can be interpreted as
$\FF_p[\pi(a_n) | n \in \NN]$.
The relation to the arithmetic of the coefficient field
$\QQ_f := \QQ(a_n | n \in \NN)$ and the coefficient ring
$\ZZ_f := \ZZ[a_n | n \in \NN]$ is apparent.

In order to conclude this background section,
we point out that in the case $k=2$,
the coefficient ring $\ZZ_f$ is the endomorphism ring of the abelian
variety~$A_f$ attached to~$f$. From that point of view, the following analysis can also be
interpreted as a study of the arithmetic of the
endomorphism algebras of $\GL_2$-abelian varieties.

\section{Semisimplicity of Hecke algebras}\label{sec:ss}

We recall that a finite dimensional commutative $K$-algebra, where $K$ is a field,
is {\em semisimple} if and only if it is isomorphic to a direct product of fields
(which are necessarily finite field extensions of~$K$).

In this section we first study the semisimplicity of the Hecke algebra~$\That_\eta$.
In the case when it is semisimple,
we relate the non-semisimplicity of the mod~$p$ Hecke algebra~$\Tbar$
to three phonomena: congruences between $\Gal(\Qbar_p/\QQ_p)$-conjugacy classes
of newforms, ramification at~$p$ of the coefficient fields of newforms and the
$p$-index of the local coefficient ring in the ring of integers of the local
coefficient field.

Let $f = \sum_{n=1}^\infty a_n(f) q^n \in S_k(\Gamma_1(M))^\new$
be a normalised Hecke eigenform and let $m$ be
any positive integer. We define the $\CC$-vector space $V_f(m)$ to be the
span of $\{f(q^d) \;|\; d \mid m\}$, where $d$ runs through all
positive divisors of~$m$ (including $1$ and~$m$). Newform theory states that
$$ S_k(\Gamma_1(N)) \cong
\bigoplus_{m \mid N} \bigoplus_{f\in S_k(\Gamma_1(N/m))^\new} V_f(m).$$
This is an isomorphism of Hecke modules. The Hecke operators $T_n$ for
$(n,m) = 1$ restricted to $V_f(m)$ are scalar matrices with $a_n(f)$ as
diagonal entries.
We now describe the Hecke operator $T_\ell$ on $V_f(m)$
for a prime~$\ell$. Suppose that $\ell^r \mid\mid m$. Let $\epsilon$ be the
Dirichlet character of~$f$.
Consider the $(r+1)\times(r+1)$-matrix
\begin{equation*}\label{bigmatrix}
A := A_f(m,\ell) := \left(
\begin{array}{cccccc}
   a_\ell(f)    & 1      & 0 & 0 & \ldots & 0      \\
   -\delta \epsilon(\ell) \ell^{k-1} & 0      & 1 & 0 & \ldots & 0      \\
   0         & 0      & 0 & 1 & \ldots & 0      \\
   \vdots    &        &   &   &        & \vdots \\
   0         & \ldots & 0 & 0 & 0      & 1      \\
   0         & \ldots & 0 & 0 & 0      & 0      \\
\end{array}
\right),
\end{equation*}
where $\delta=0$ if $\ell\mid N$ and $\delta=1$ otherwise.
The Hecke operator $T_\ell$ is then given on $V_f(m)$ (for a certain natural basis)
by a diagonal block matrix having only blocks equal to~$A$ on the diagonal,
where each block on the diagonal corresponds to a divisor of $N/\ell^r$.
Let $\TT$ be the Hecke algebra of $S_k(\Gamma_1(N))$ (as in Section~\ref{sec:background}).
The algebra $\TT_\QQ := \TT \otimes_\ZZ \QQ$ is semisimple if and only if
$\TT_\CC := \TT \otimes_\ZZ \CC$ is semisimple (if and only if $\That_\eta$ is semisimple).
By the above discussion,
$\TT_\CC$ is semisimple if and only if all the matrices~$A_f(m,\ell)$ that appear
are diagonalisable.

\begin{prop}
Assume the notation above, $M=N/m$ and $k\ge 2$. Moreover, if $k\ge 3$ assume Tate's conjecture
(see \cite{ColemanEdixhoven}, Section~1).
\begin{enumerate}[(a)]
\item Assume $\ell \nmid  M$. Then $A_f(m,\ell)$ is diagonalisable if and only if $r \le 2$.
\item Assume that $\ell \mid M$ and that either $\ell \mid\mid M$ or that $\epsilon$
cannot be defined mod~$M/\ell$.
Then $A_f(m,\ell)$ is diagonalisable if and only if $r \le 1$.
\item Assume that $\ell^2 \mid M$ and that $\epsilon$ can be defined modulo~$M/\ell$.
Then $A_f(m,\ell)$ is diagonalisable if and only if $r = 0$.
\end{enumerate}
\end{prop}

\pf
(a) Assume $r \ge 1$ (otherwise the result is trivial) and call $B$ the top left
$2 \times 2$-block of~$A = A_f(m,\ell)$. The characteristic polynomial of~$B$
is $g(X) = X^2 - a_\ell(f) X + \epsilon(\ell)\ell^{k-1}$.
We have $g(0) \neq 0$ and $g(X)$ has discriminant
$a_\ell(f)^2 - 4 \epsilon(\ell)\ell^{k-1}$, which is non-zero,
since $|a_\ell(f)| = 2 |\ell|^{(k-1)/2}$ would contradict \cite{ColemanEdixhoven},
Theorem~4.1. Consequently, $A$ is diagonalisable if and only if apart from~$B$ there is at most
one more row and column.

In cases (b) and~(c), note that $A$ is in Jordan form. The result is now immediate,
since $a_\ell(f)$ is non-zero for~(b) and zero for~(c) (see \cite{DeligneSerre}, 1.8).
\qed
\medskip

We have the immediate corollary (which is Theorem~4.2 in~\cite{ColemanEdixhoven}).

\begin{cor}\label{cor:nondegzero}
Let $N$ be cubefree and $k\ge 2$. If $k \ge 3$ assume Tate's conjecture
(see \cite{ColemanEdixhoven}, Section~1). Then the Hecke algebras
$\TT_k(\Gamma_1(N)) \otimes \QQ$ and $\That_k(\Gamma_1(N))_\eta$
as well as $\TT_k(N) \otimes \QQ$ and $\That_k(N)_\eta$ are semisimple.
\end{cor}

The principal result of this section is the following proposition on the
structure of the residual Hecke algebra. We assume
the notation laid out in Section~\ref{sec:background}, in particular, we work
with a general congruence subgroup~$\Gamma$.

\begin{prop}
Assume that $\That_\eta$ is semisimple (see, e.g., Corollary~\ref{cor:nondegzero}), i.e.\
$\That_\eta \cong \prod_{\fp \in \MinSpec(\That)} L_\fp$.
Then the residual Hecke algebra $\Tbar$ is semisimple if and only if
all of the following three conditions are satisfied:
\begin{enumerate}[(i)]
\item No two $\fp_1, \fp_2 \in \MinSpec(\That)$ are congruent.
\item None of the $\fp \in \MinSpec(\That)$ is ramified.
\item For all $\fp \in \MinSpec(\That)$, the index $i_\fp = 1$.
\end{enumerate}
\end{prop}

\pf
We first prove that (i), (ii) and (iii) imply the semisimplicity of~$\Tbar$.

The fact that there is no congruence means that in every $\fm \in \MaxSpec(\That)$
there is a unique $\fp \in \MinSpec(\That_\fm)$.
As $\That_\fm \otimes_{\ZZ_p} \QQ_p \cong L_\fp$, it follows that
$\That_\fm$ is a subring of~$L_\fp$.
Due to~(iii), $\That_\fm$ is the ring of integers of~$L_\fp$. Since by~(ii)
$L_\fp$ is unramified, we get that $\Tbar_\fm$ is the
residue field of the integers of~$L_\fp$.
This shows that $\Tbar$ is a product of finite fields, i.e.\ semisimple.

Now we prove the converse direction and assume that $\Tbar$ is semisimple.
Let $\fm \in \MaxSpec(\That)$.
Let $\fp_1, \fp_2, \dots, \fp_m \in \MinSpec(\That)$ be the distinct minimal primes
contained in~$\fm$.
Then $\That_\fm \otimes_{\ZZ_p} \QQ_p \cong L_{\fp_1} \times \dots \times L_{\fp_m}$.
Due to the non-degeneration
$\Tbar_\fm \cong \That_\fm \otimes_{\ZZ_p} \FF_p \cong \FF_{p^n}$ for some~$n$.
Since $\dim_{\QQ_p} \That_\fm \otimes_{\ZZ_p} \QQ_p = n$, we have
$[L_{\fp_i} : \QQ_p] \le n$ for $i=1,\dots,m$.

Let $\cO_i$ be the ring of integers of~$L_{\fp_i}$ for $i=1,\dots,m$.
It contains $\That_\fm/\fp_i$ with index~$i_{\fp_i}$. Tensoring the exact sequence
of $\ZZ_p$-modules
$$ 0 \to \That_\fm/\fp_i \to \cO_i \to \cO_i/(\That_\fm/\fp_i) \to 0$$
with~$\FF_p$ over~$\ZZ_p$ we obtain the exact sequence of $\FF_p$-vector spaces:
$$ \FF_{p^n} \to \cO_i \otimes_{\ZZ_p} \FF_p \to (\cO_i/(\That_\fm/\fp_i))\otimes_{\ZZ_p} \FF_p \to 0.$$
Since the map on the left is a ring homomorphism, it is injective.
Now $\dim_{\FF_p} \cO_i \otimes_{\ZZ_p} \FF_p \le n$
implies that $\cO_i$ is unramified and that $i_{\fp_i} = 1$.
Thus $[L_{\fp_i} : \QQ_p] = n$ for $i=1,\dots,m$ and, hence, $m=1$, concluding the proof.
\qed

\section{Observations and Questions}\label{sec:observations}

In this section, we explain and expose our computer experiments and we ask
some questions suggested by our studies.
All computer calculations were performed using {\sc Magma} (see \cite{magma}).

\subsection{Semisimplicity of the residual Hecke algebra}\label{sec:a}

A local finite-dimensional commutative $\FF_p$-algebra~$A$ is semisimple
if and only if it is simple, which is equivalent to~$A$ being field.
We take the dimension of the maximal ideal~$\fm$ of~$A$
as a measure for the deviation of~$A$ from being semisimple.
In particular, $A$ is a field if and only if $\fm$ has dimension~$0$.

For given prime~$p$, level~$N$ and weight~$k$ we study the {\em sum of the residue degrees} of
all prime ideals:
$$ a_{N,k}^{(p)} = \sum_{\fm \in \Spec(\Tbar_k(N))} [\FF_\fm:\FF_p].$$
Clearly, $a_{N,k}^{(p)}$ is less than or equal to the $\Fbar_p$-dimension of $S_k(N;\Fbar_p)$.
Hence, $\Tbar_k(N)$ is semisimple if and only if
$a_{N,k}^{(p)}$ is equal to this dimension.

We intend to study the asymptotic behaviour of the function $a_{N,k}^{(p)}$
for a fixed prime~$p$ and fixed weight~$k$ as a function of the level~$N$.
For simplicity, we let $N$ run through the prime numbers only in order to
avoid contributions from lower levels via the degeneracy maps. We should
point out that there can be contributions from lower weights: an eigenform
in $S_k(N;\Fbar_p)$ also lives in $S_{k+n(p-1)}(N;\Fbar_p)$ for all~$n \ge 0$
by multiplication by the Hasse invariant. Note that for $p>2$ and $k=2$,
as well as for $p>3$ and $k=4$ there is no such contribution.

Our computational findings are best illustrated by plotting graphs.
In each of the following plots, the prime~$p$ and the weight~$k$ are fixed and
on the $x$-axis we plot $d(N) := \dim_{\Fbar_p} S_k(N;\Fbar_p)$ and on the $y$-axis
the function $a_{N,k}^{(p)}$ as a function of~$N$, i.e.\ each~$N$ gives
rise to a dot at the appropriate place.
The green line in the graphs was determined as the linear function $\alpha \cdot d(N)$
which best fits the data (according to gnuplot and the least squares method).

In the weights that we considered we observed a behaviour for $p=2$
which seems to be completely different from the behaviour at all other primes.
We made plots for all odd primes less than~$100$ and weight~$2$ and present a selection here.
The graphs that we leave out look very similar.

\noindent\begin{longtable}{cc}
\mplotr{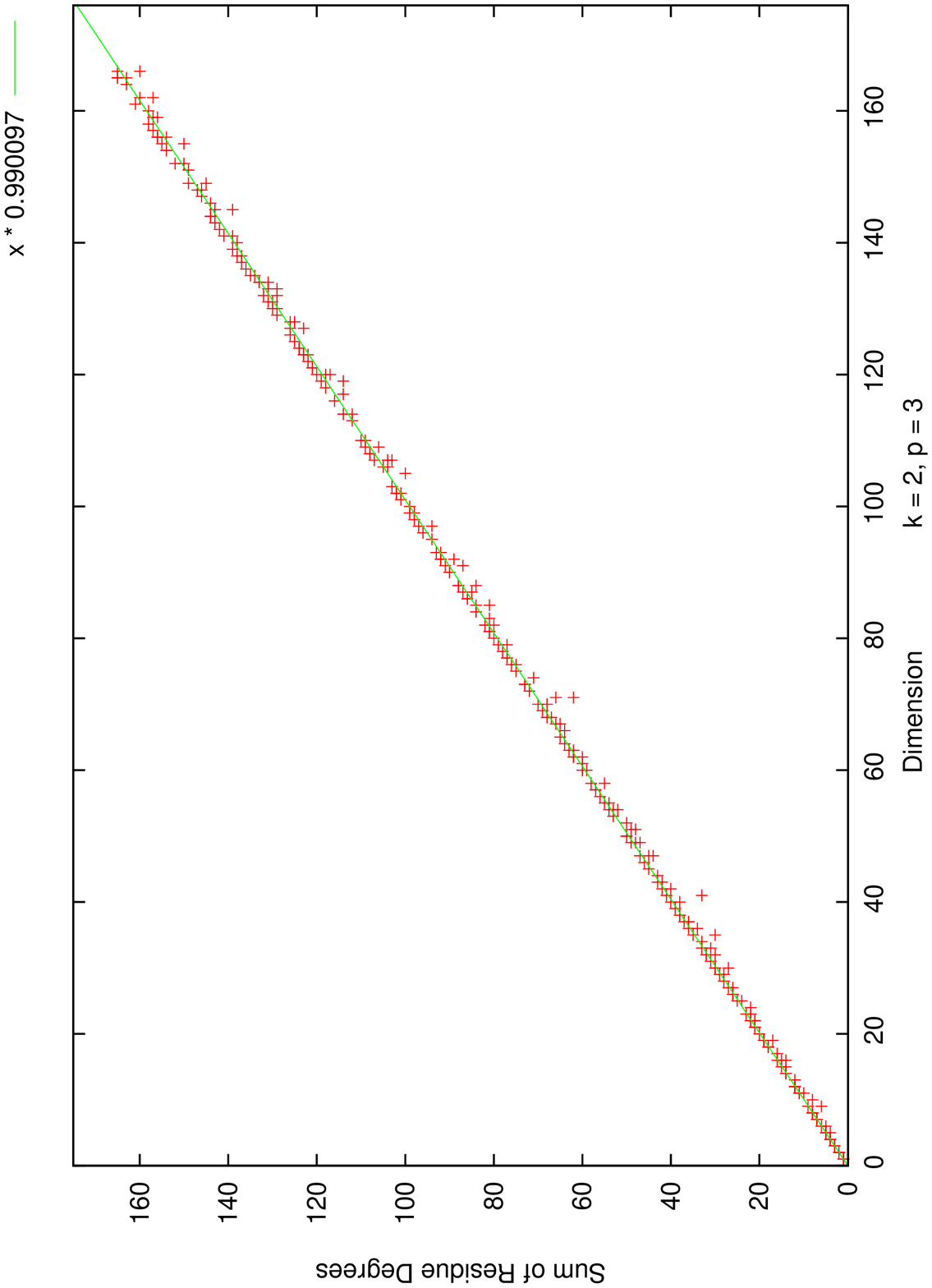}{a:2:3} &
\mplotr{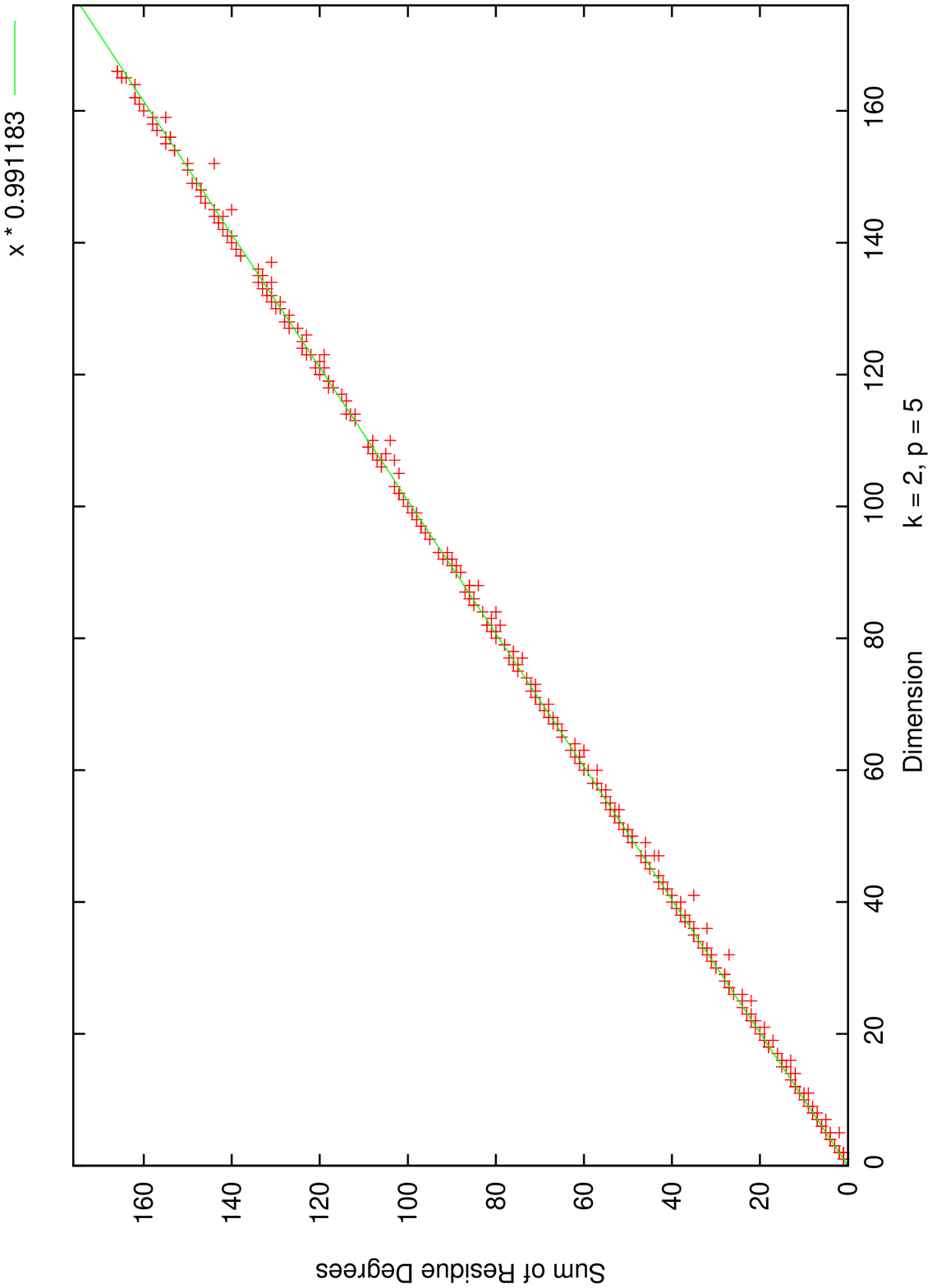}{a:2:5}\\
\mplotr{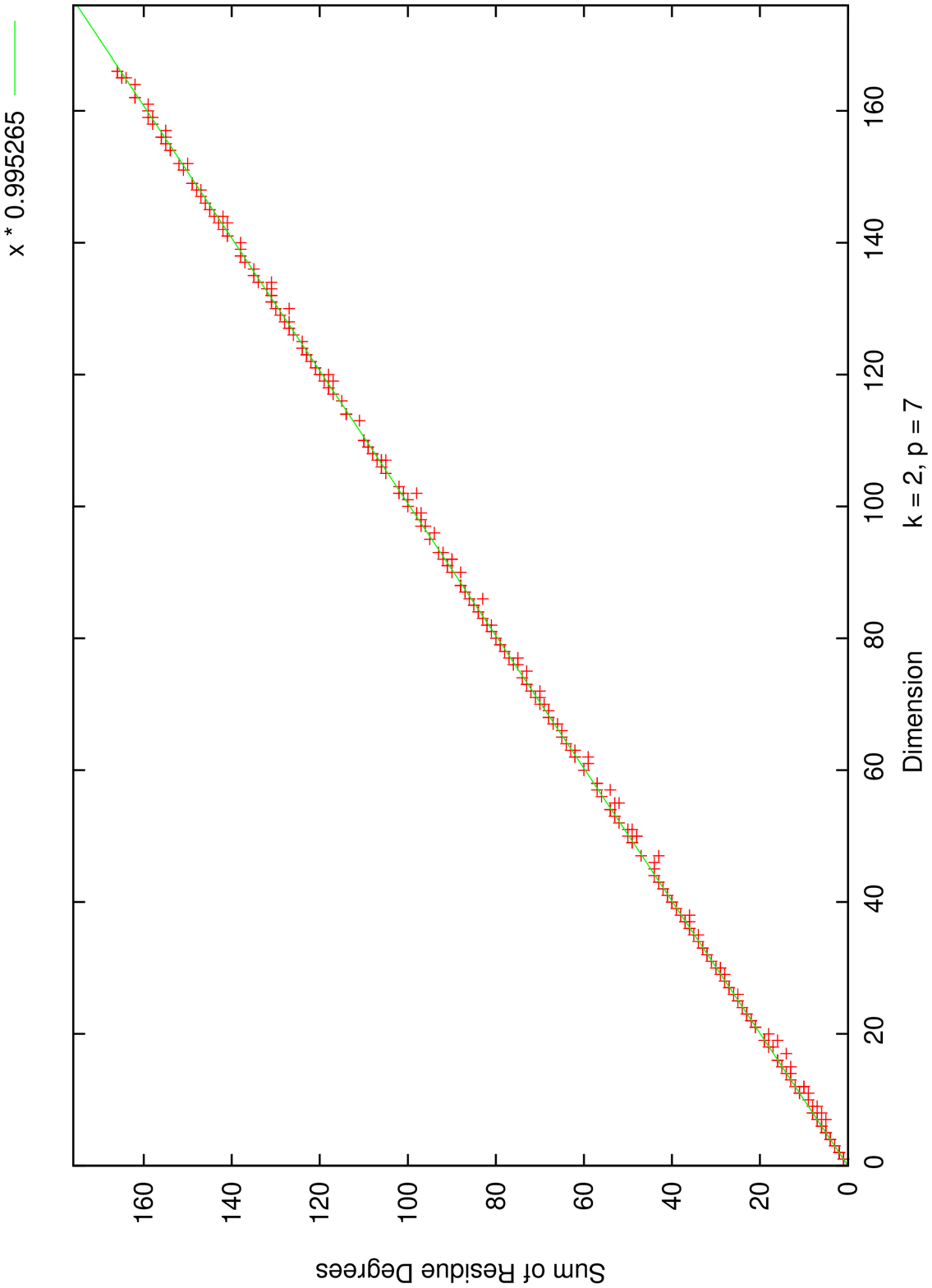}{a:2:7} & 
\mplotr{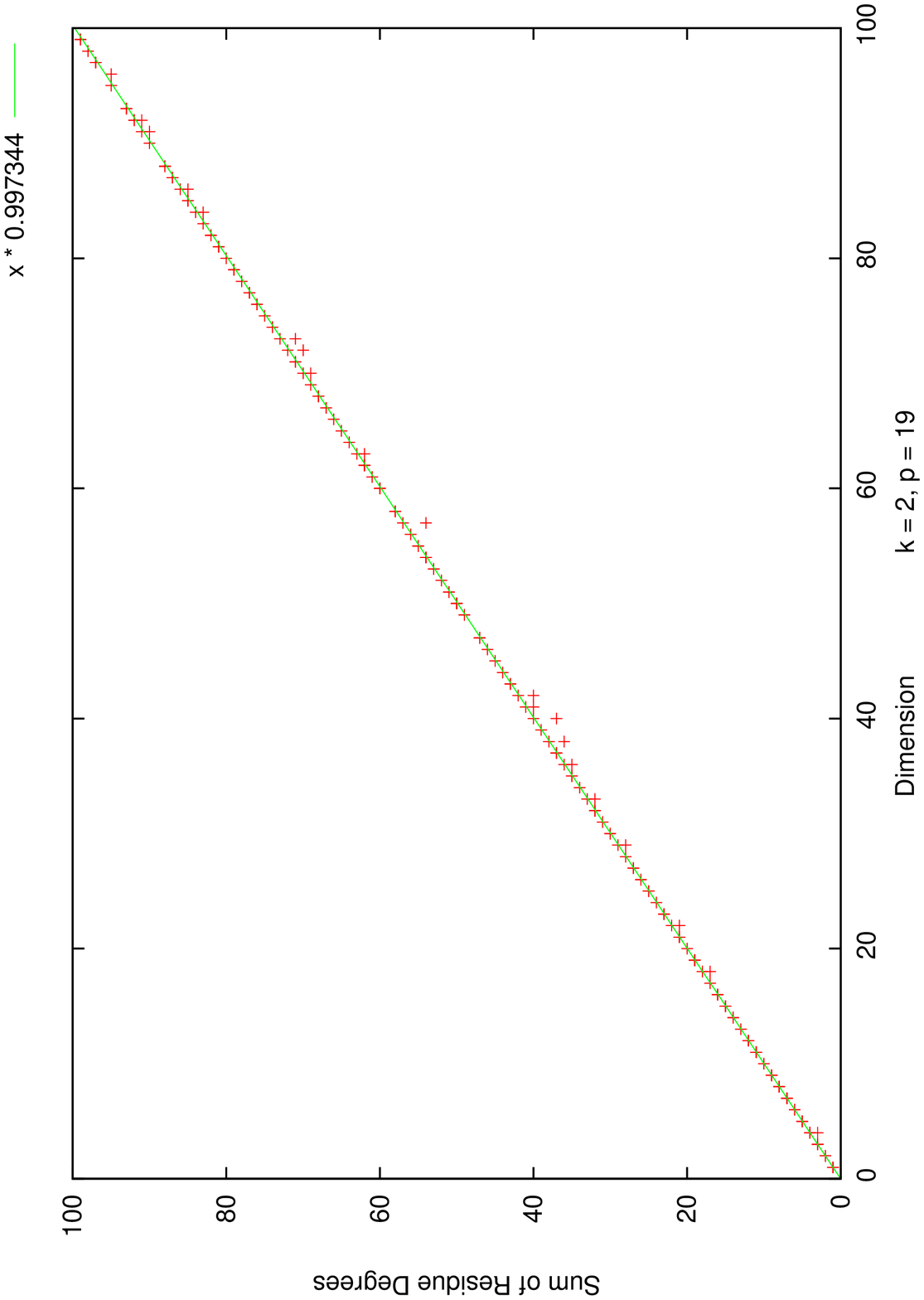}{a:2:19}\\
\mplotr{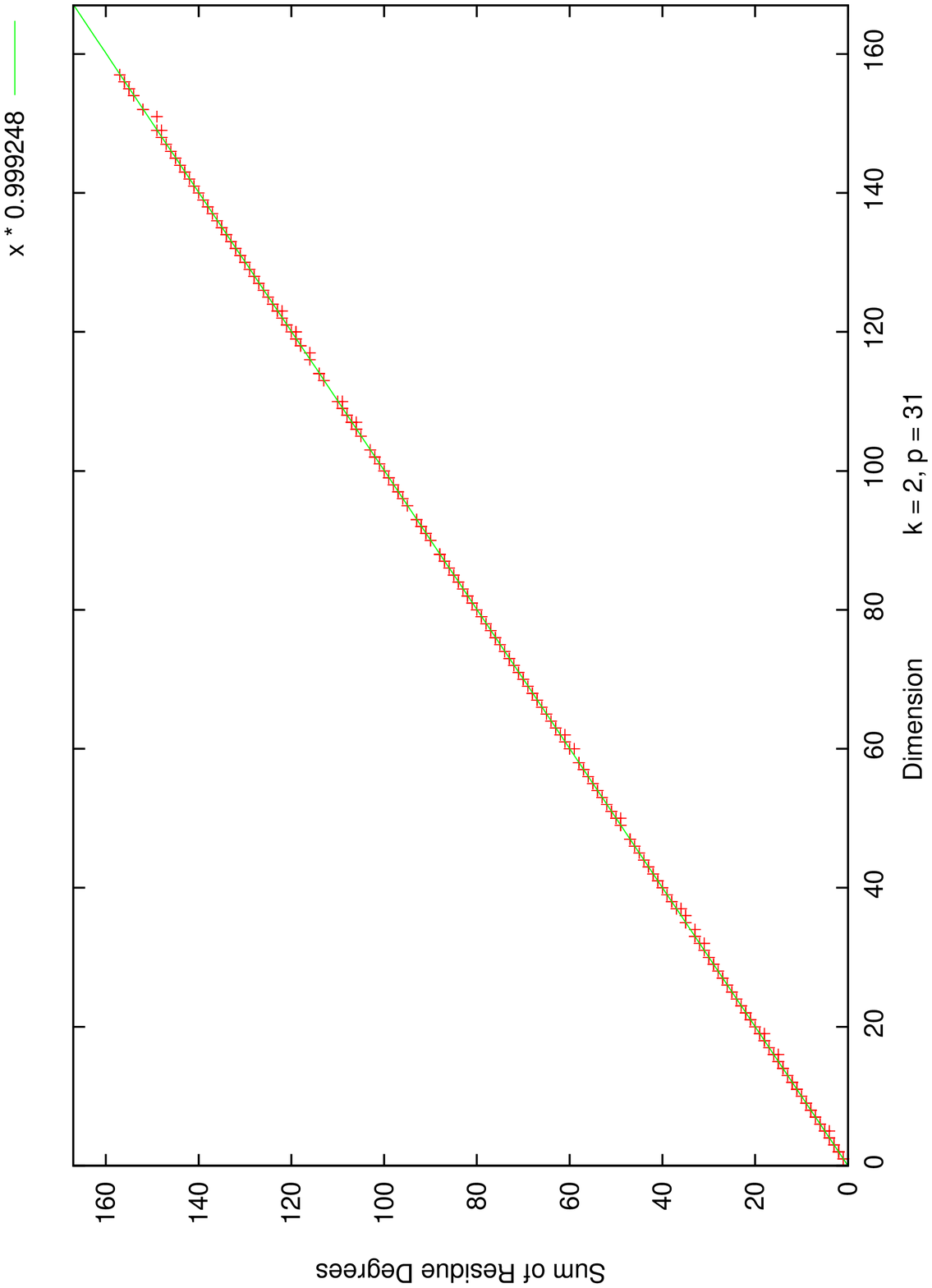}{a:2:31}&
\mplotr{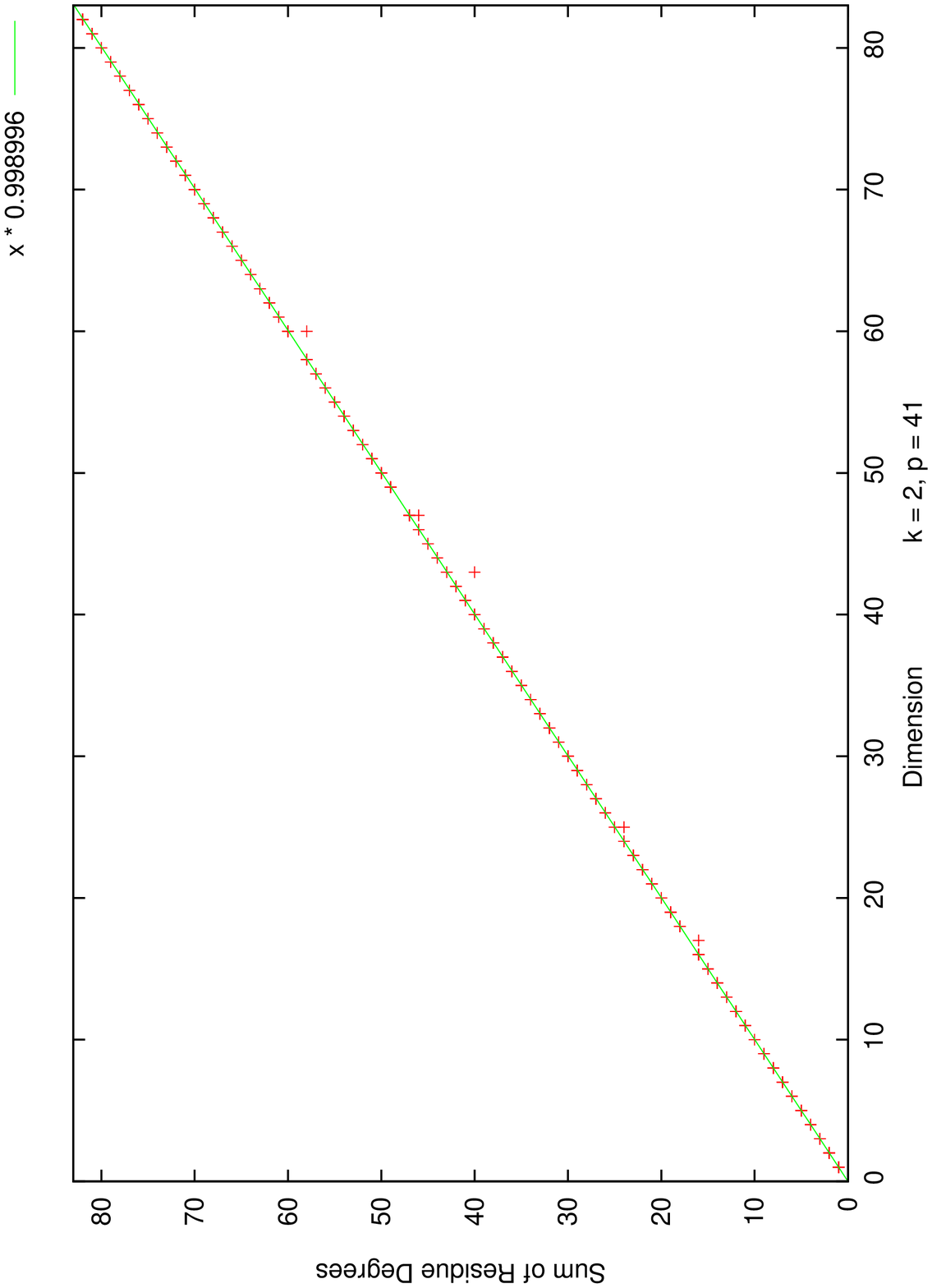}{a:2:41}\\
\mplotr{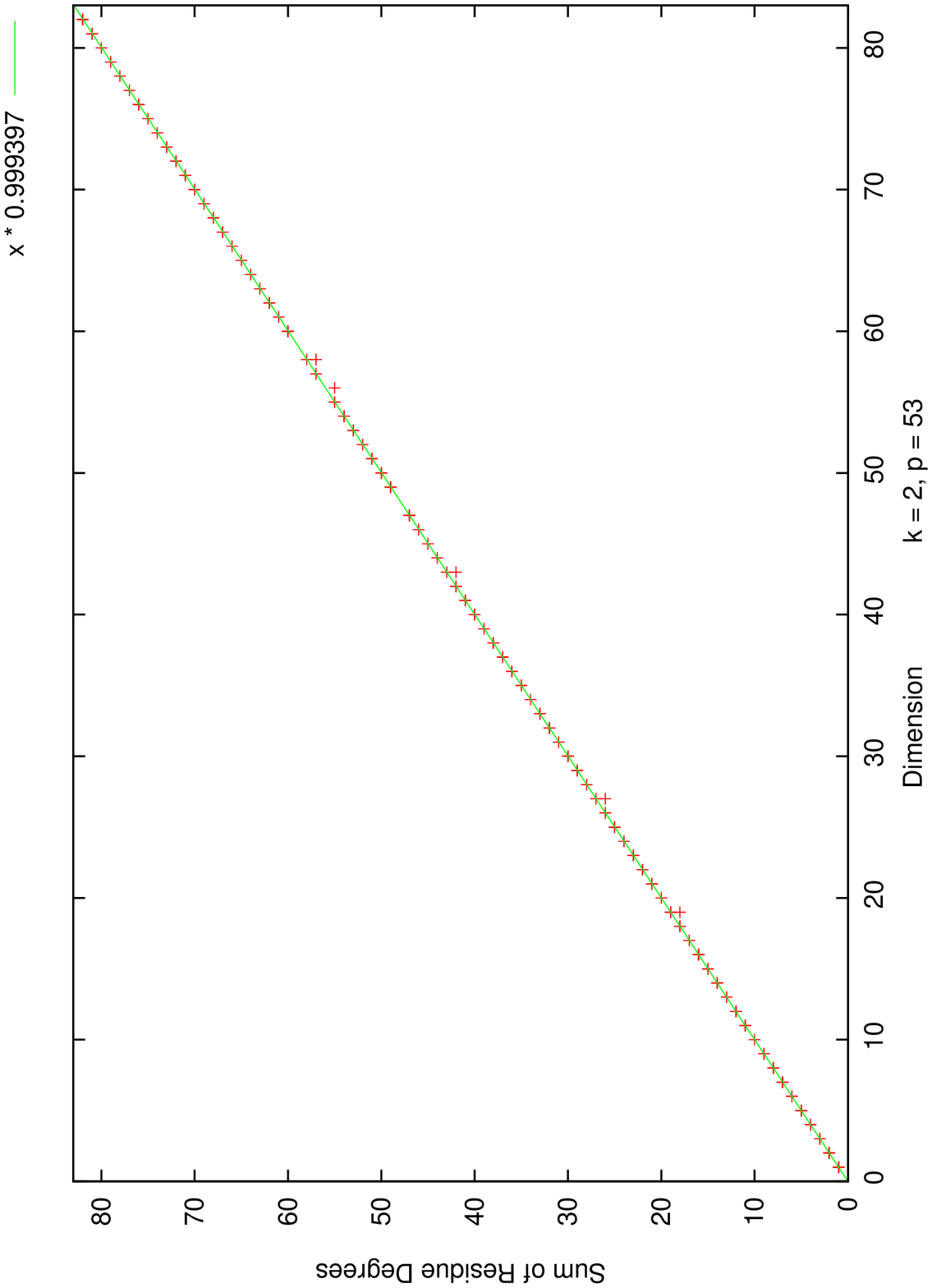}{a:2:53}& 
\mplotr{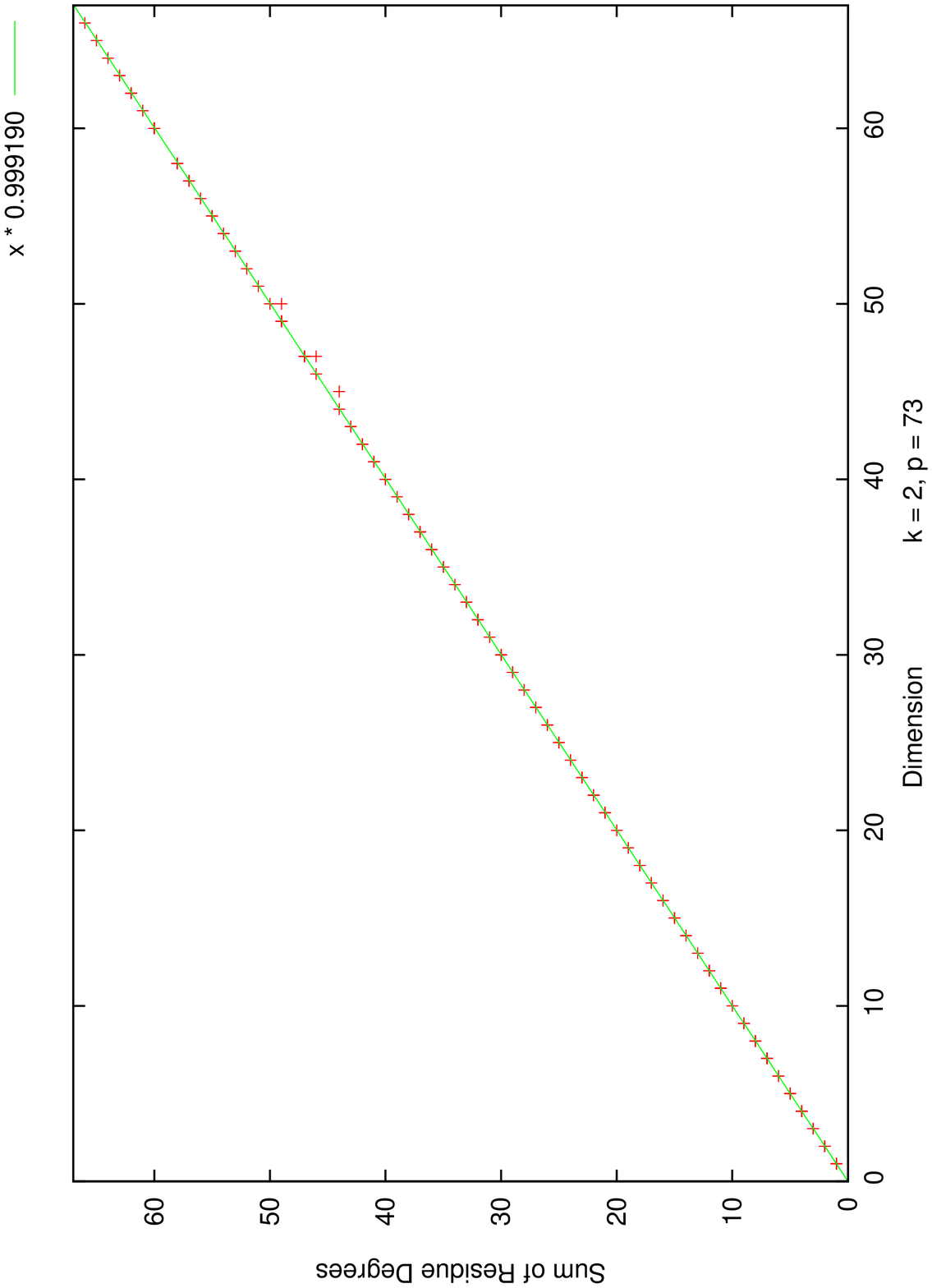}{a:2:73}\\
\mplotr{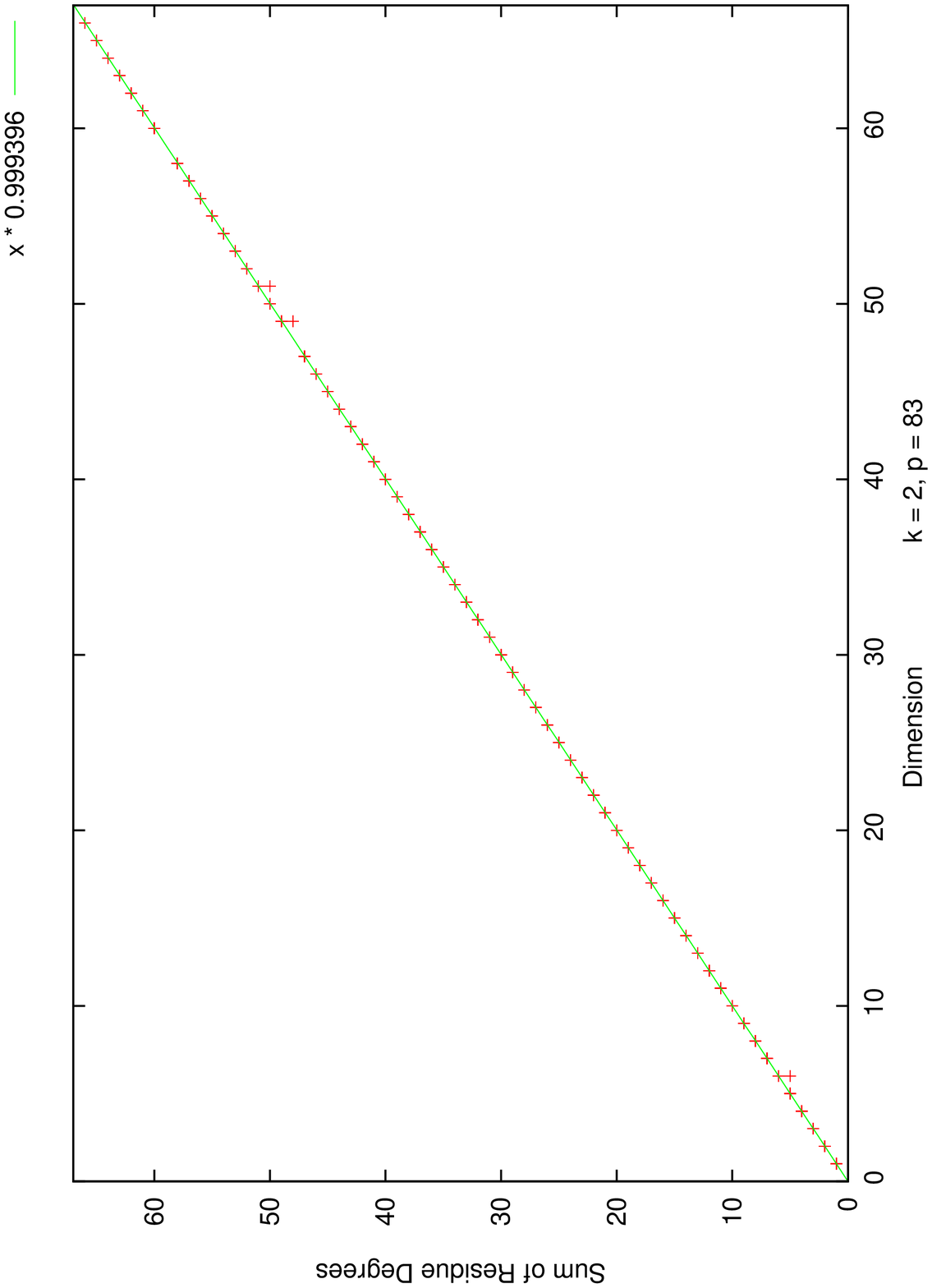}{a:2:83}&
\mplotr{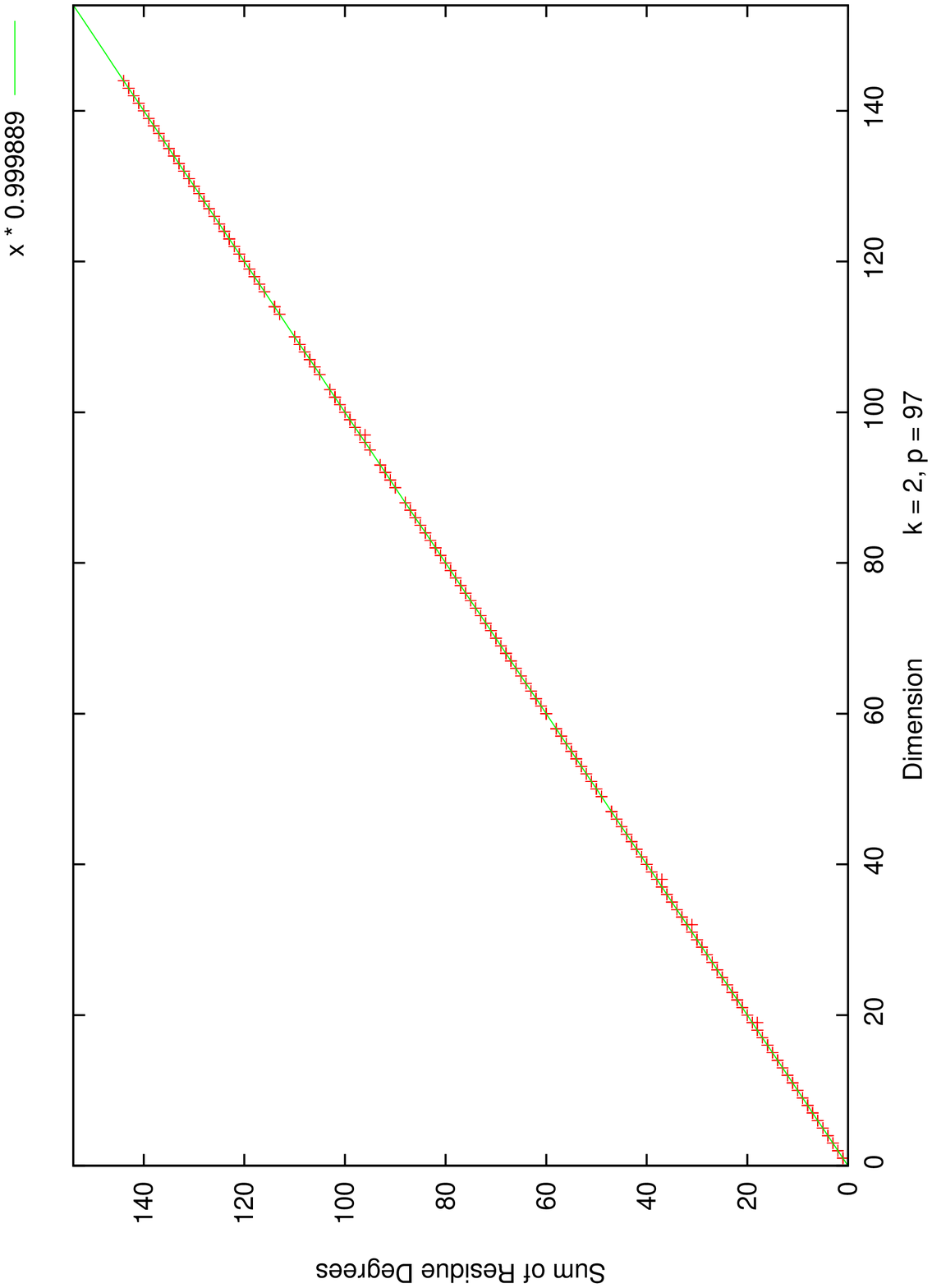}{a:2:97} \\
\end{longtable}

In Figures \ref{a:2:3}--\ref{a:2:97} the levels range over all primes up to a certain
bound (which is not the same for all~$p$).
We observe that non-semisimplicity seems to be a rather rare phenomenon which becomes
rarer for growing~$p$, as one might have guessed.
In the next figures, we analyse the cases $p=3,5$ still for weight~$2$ more closely by letting
the levels range through all primes between $3000$ and $10009$ subdivided into four
consecutive intervals.

\noindent\begin{longtable}{cc}
\mplot{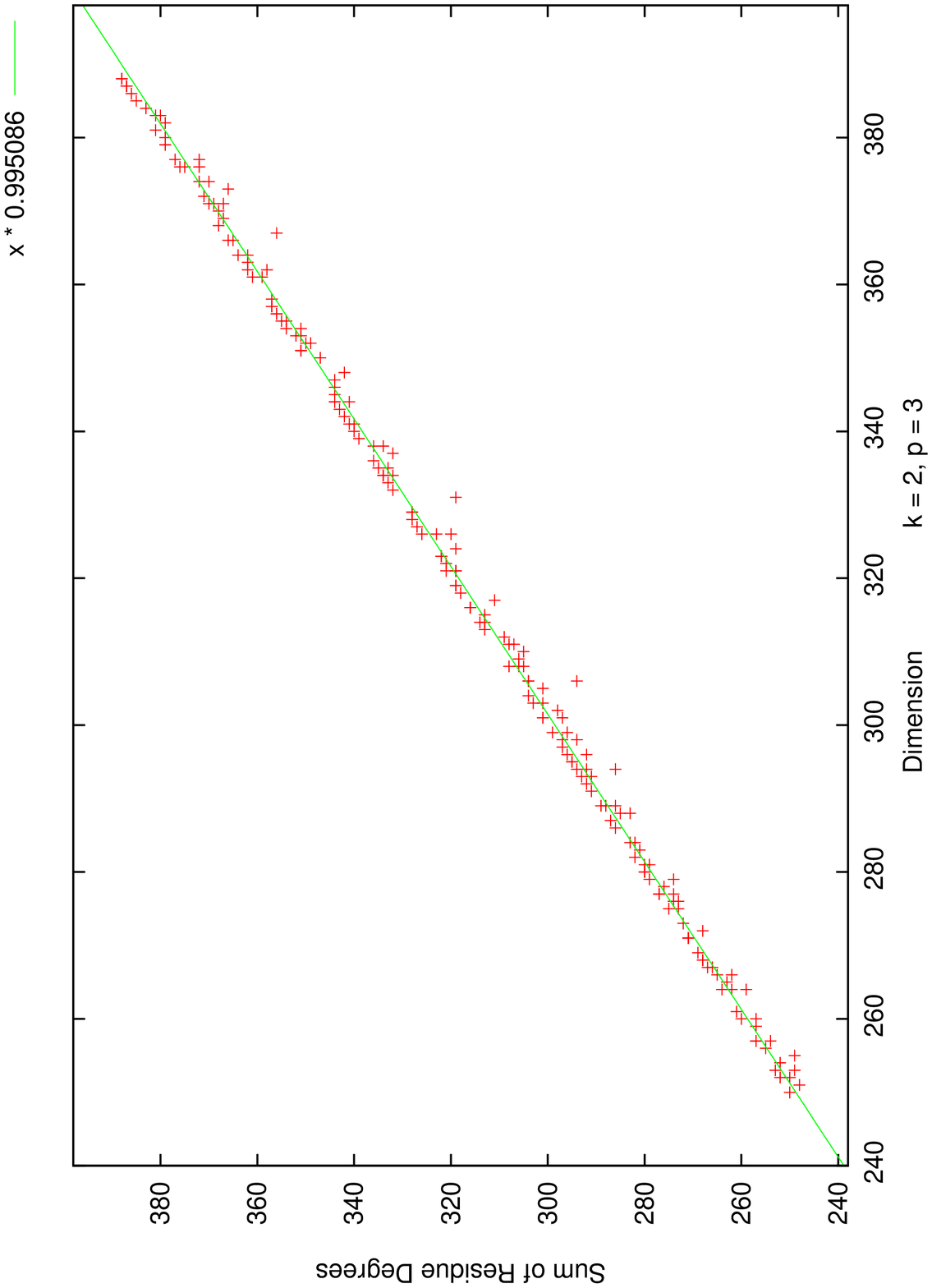} &
\mplot{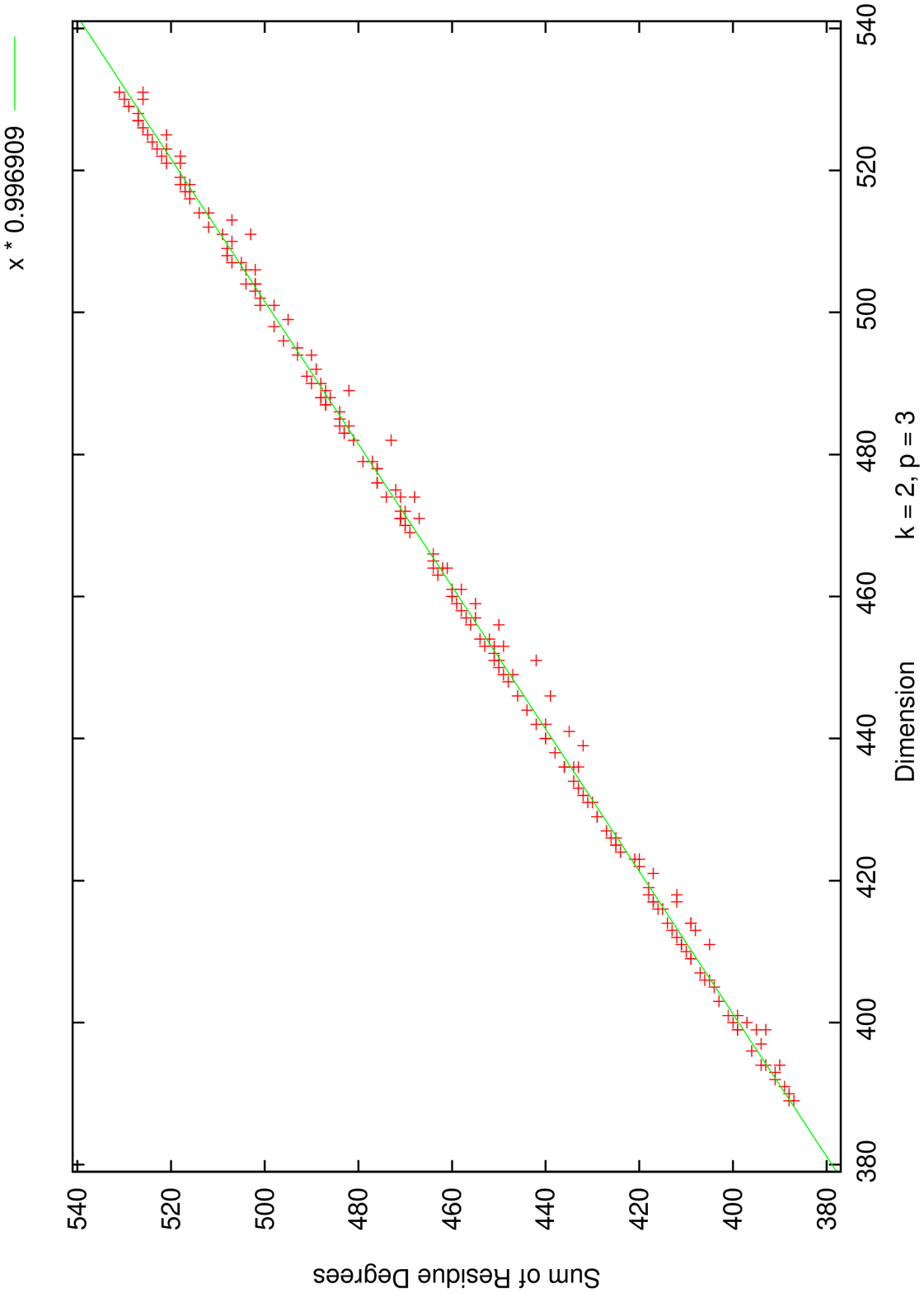} \\
\mplot{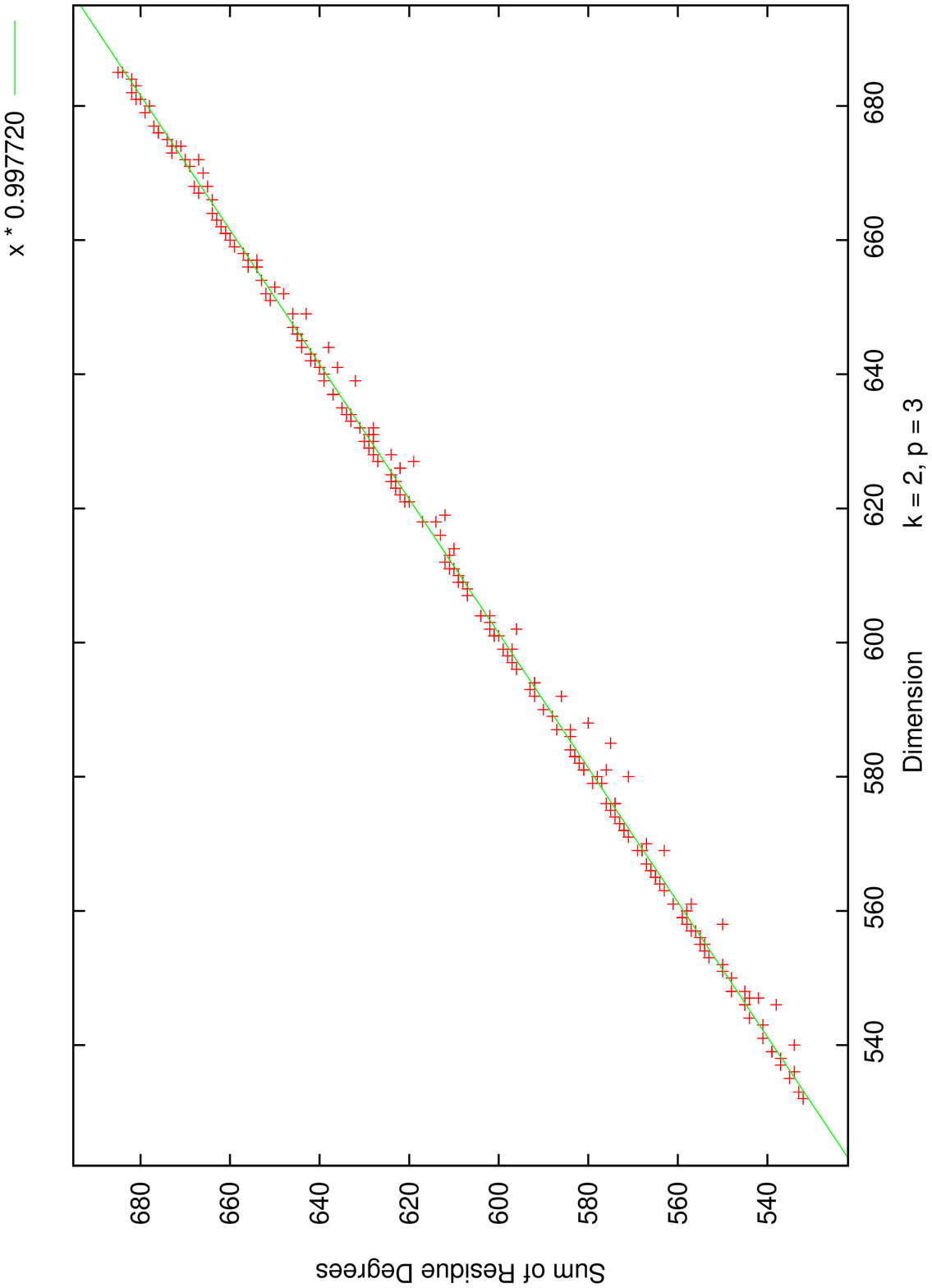} &
\mplot{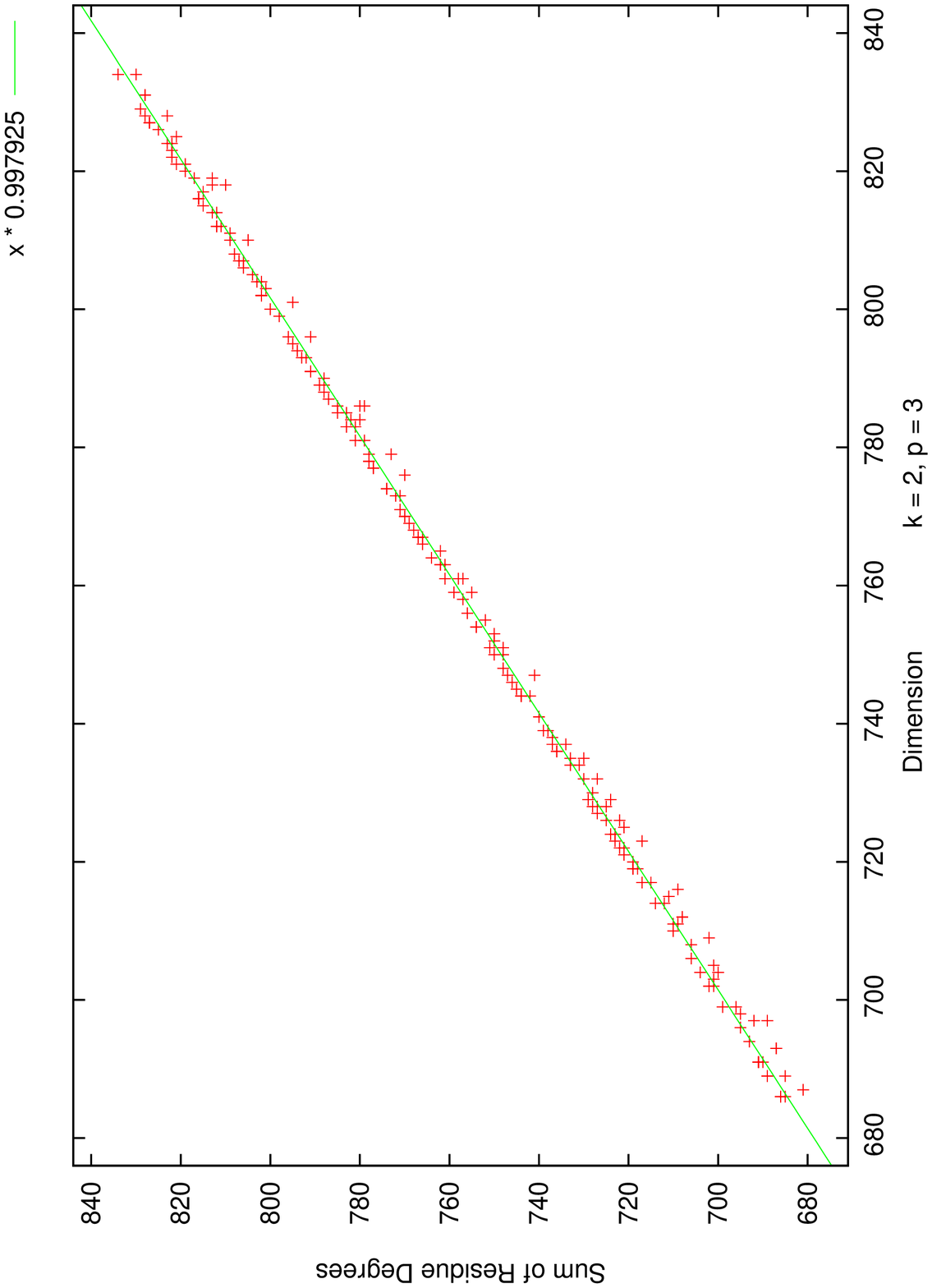} \\
\mplot{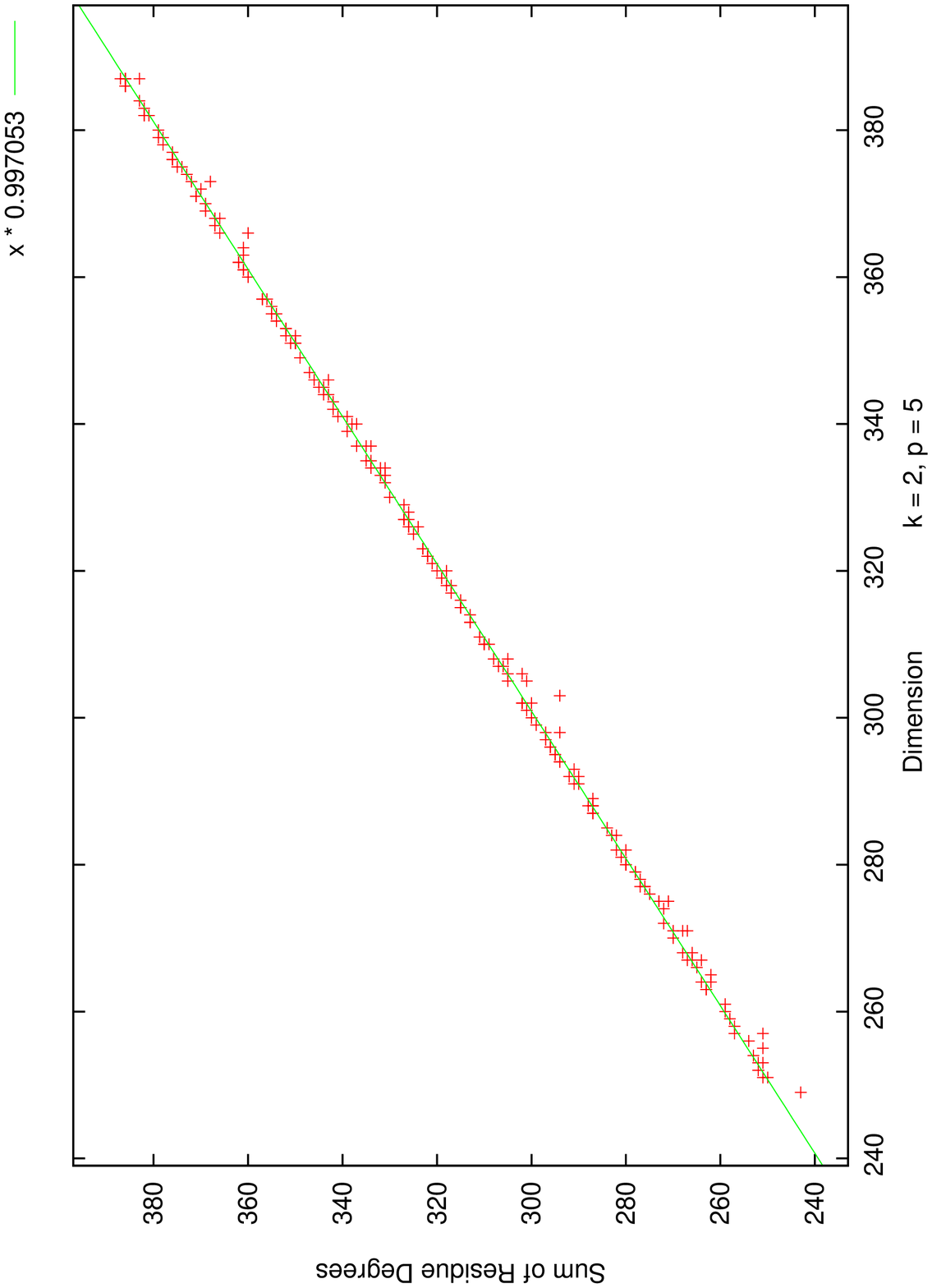} &
\mplot{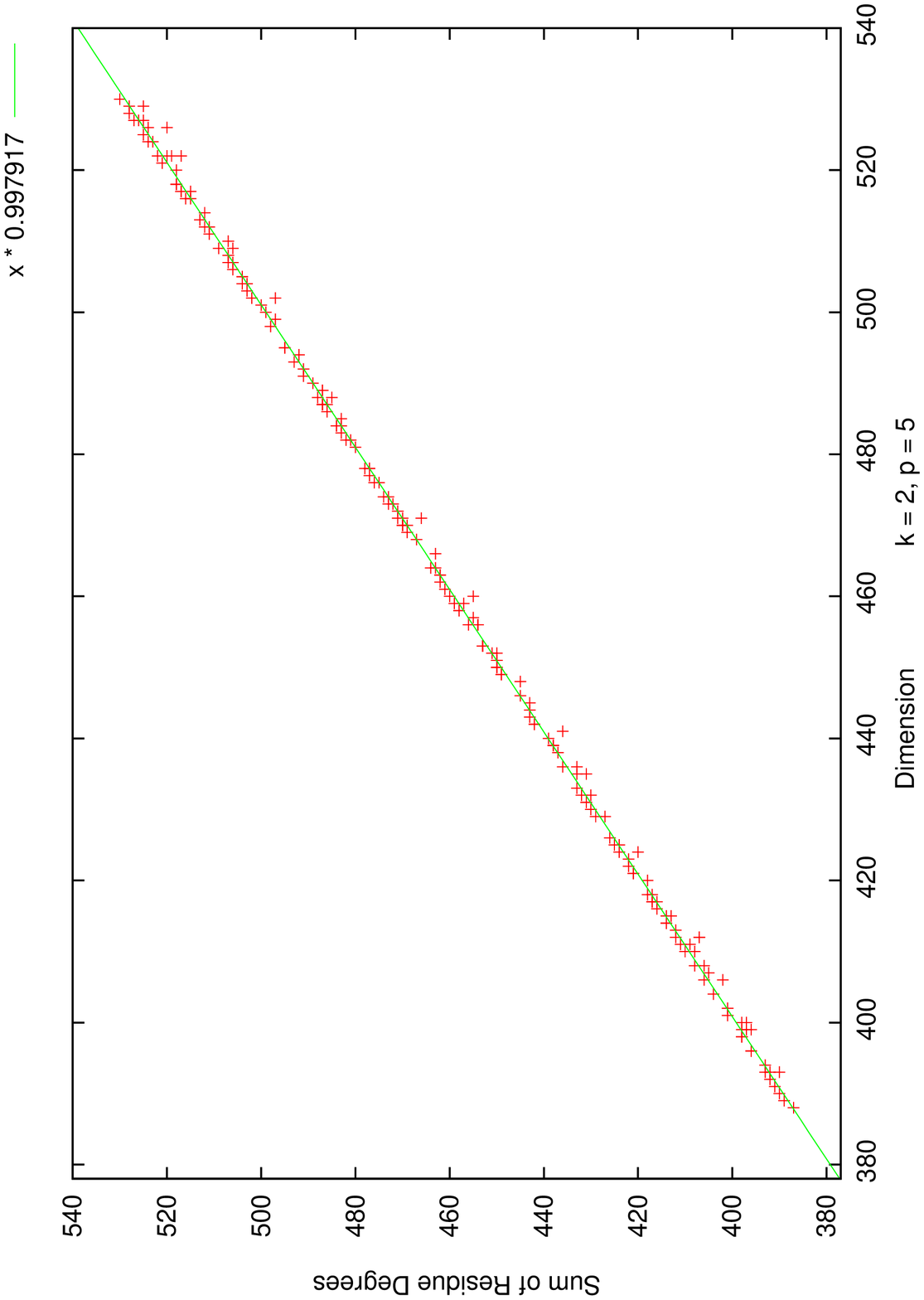} \\
\mplot{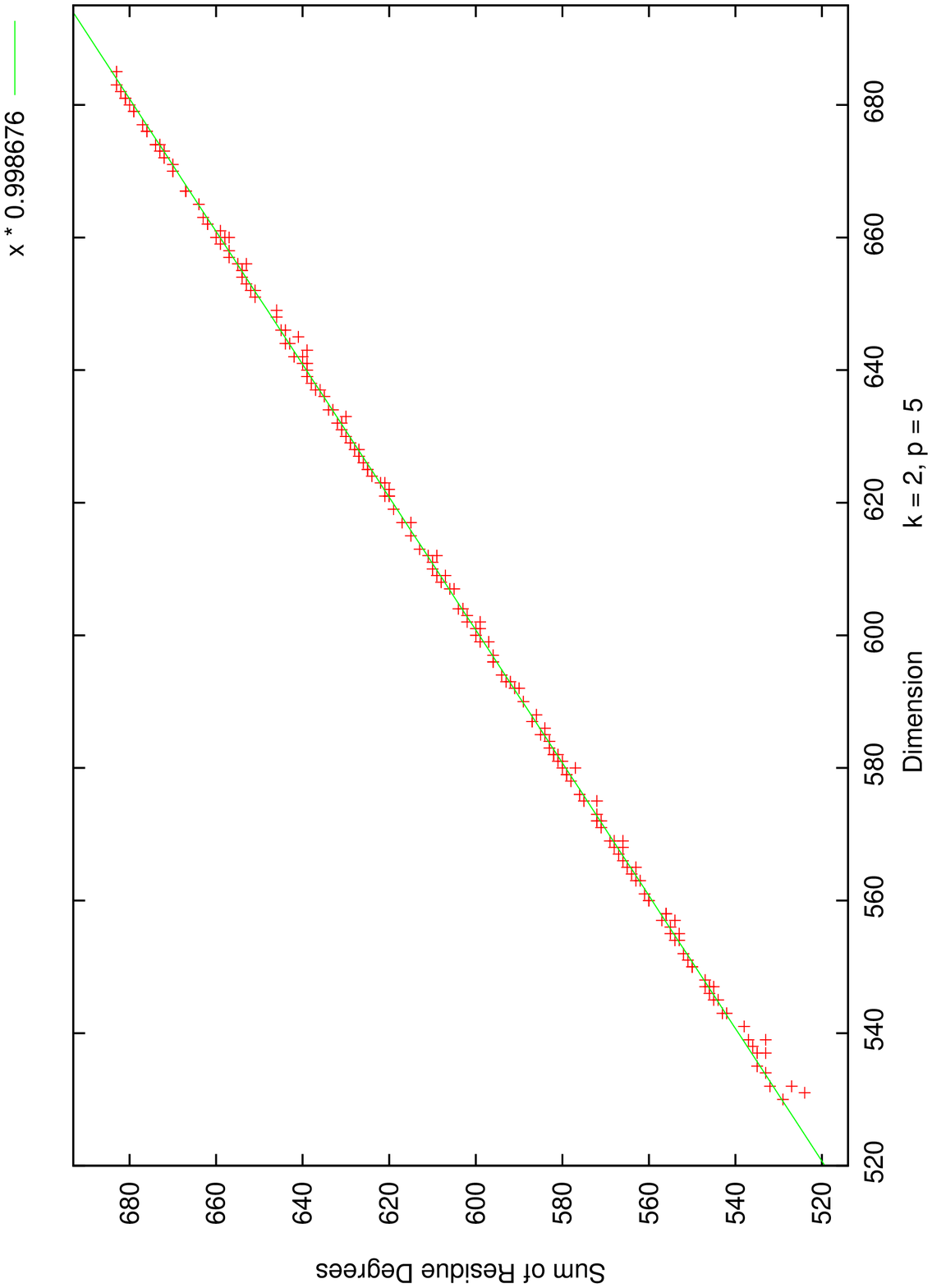} &
\mplot{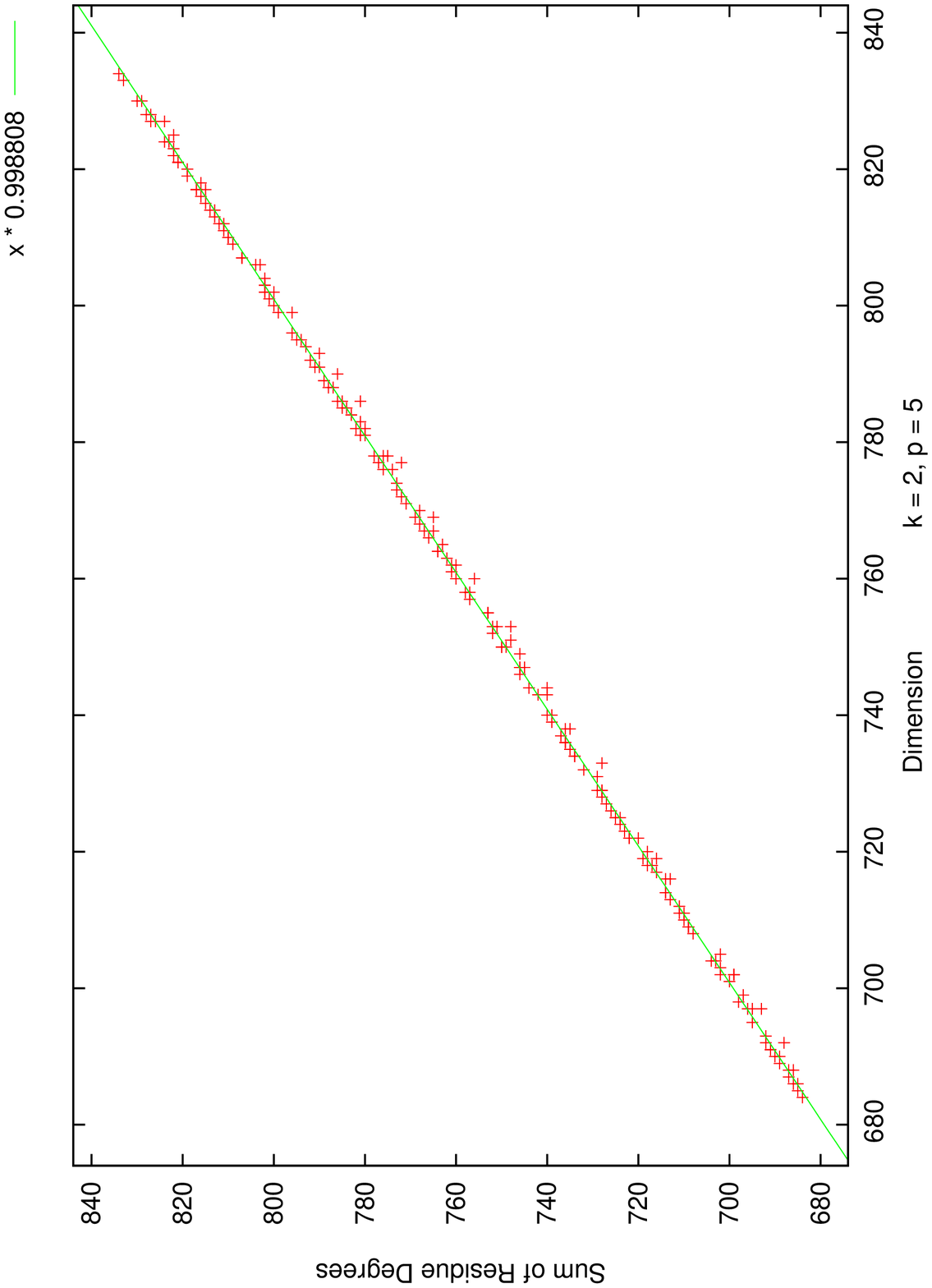} \\
\end{longtable}

One can observe that the slope of the best fitting line through the origin seems to be increasing
with growing dimension.
Although we only computed relatively litte data, we include two examples for weight~$4$.
They do not suggest any significant difference to the weight~$2$ case.

\noindent\begin{longtable}{cc}
\mplot{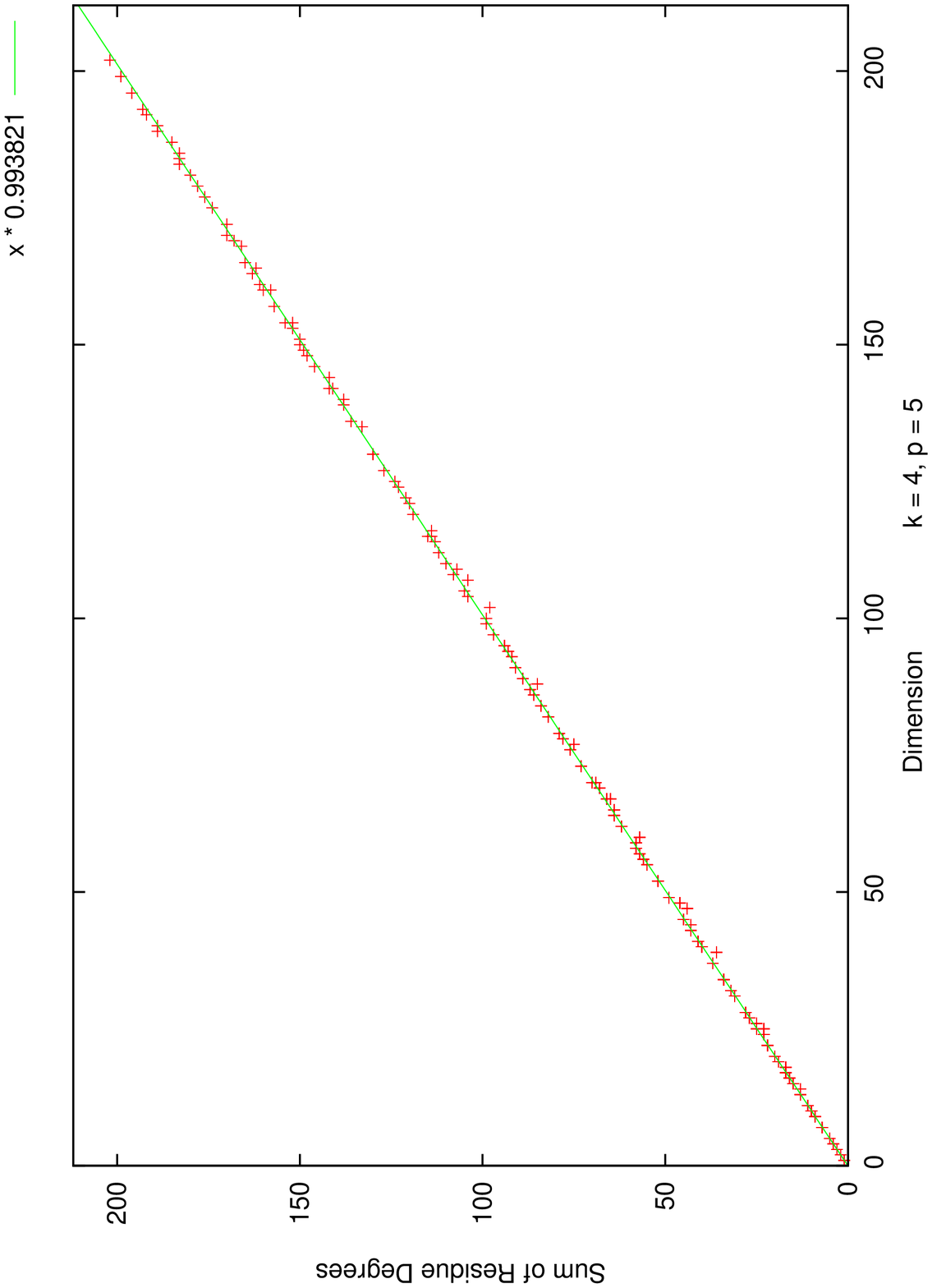} &
\mplot{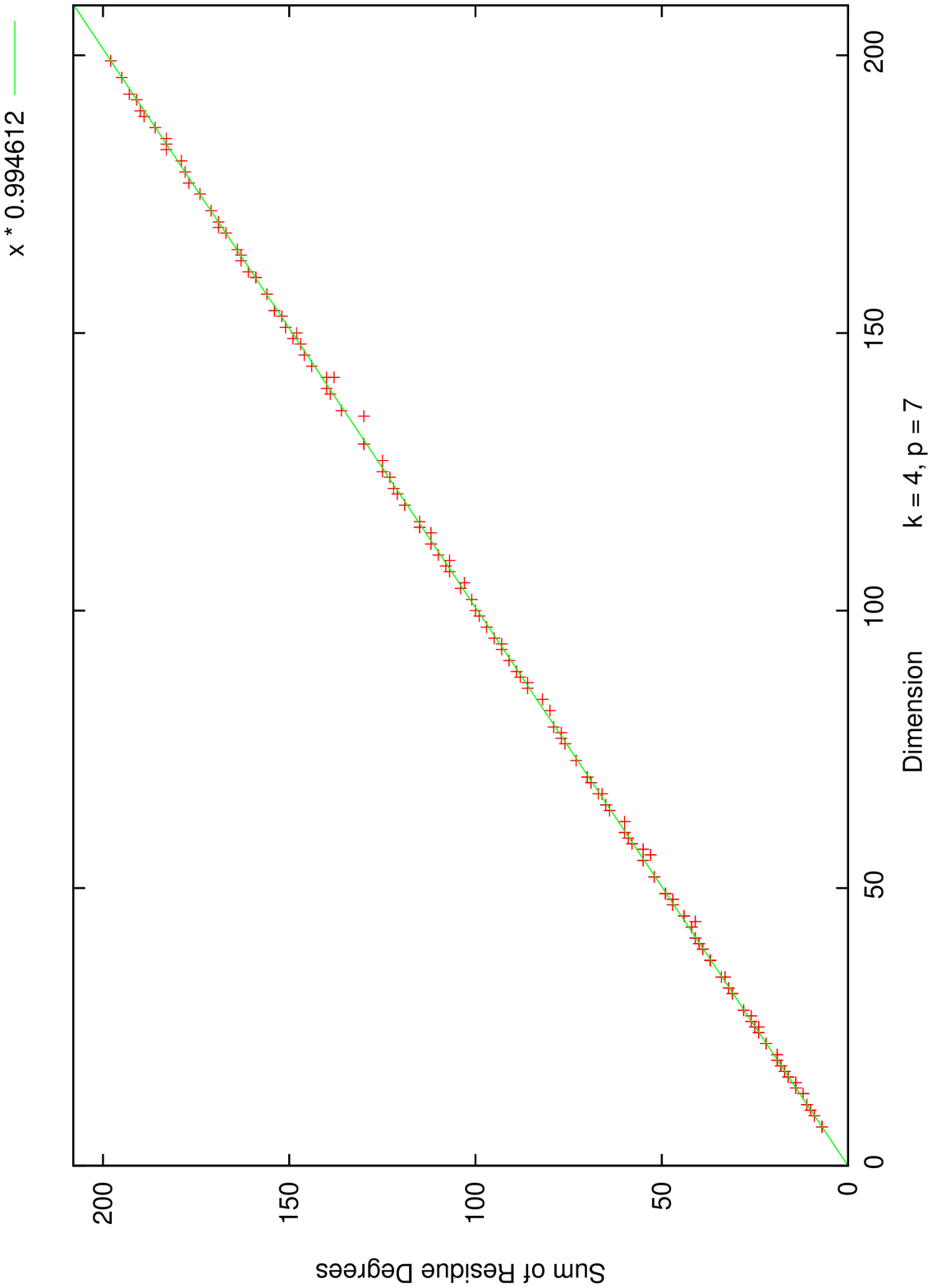}\\
\end{longtable}

We are led to ask the following question.

\begin{question}\label{question:a-odd}
Fix an odd prime~$p$ and an even $k \ge 2$.
Let $a(N) := a_{k,N}^{(p)}$ and $d(N) := \dim_{\Fbar_p} S_k(N;\Fbar_p)$.
Does the following statement hold?
\begin{quote}
For all $\epsilon>0$ there is $C_\epsilon > 0$
such that for all primes~$N$ the inequality
$$ a(N) > (1-\epsilon) d(N) - C_\epsilon$$
holds.
\end{quote}
\end{question}

We contrast the situation, which seems very similar for every odd prime,
with the one for $p=2$ and $k=2$. We do not consider any higher weights
due to the contributions from weight~$2$, which would 'disturb' the situation. The following
plots take prime numbers~$N$ into account that lie in six different intervals up to~$12000$.

\noindent\begin{longtable}{cc}
\mplotr{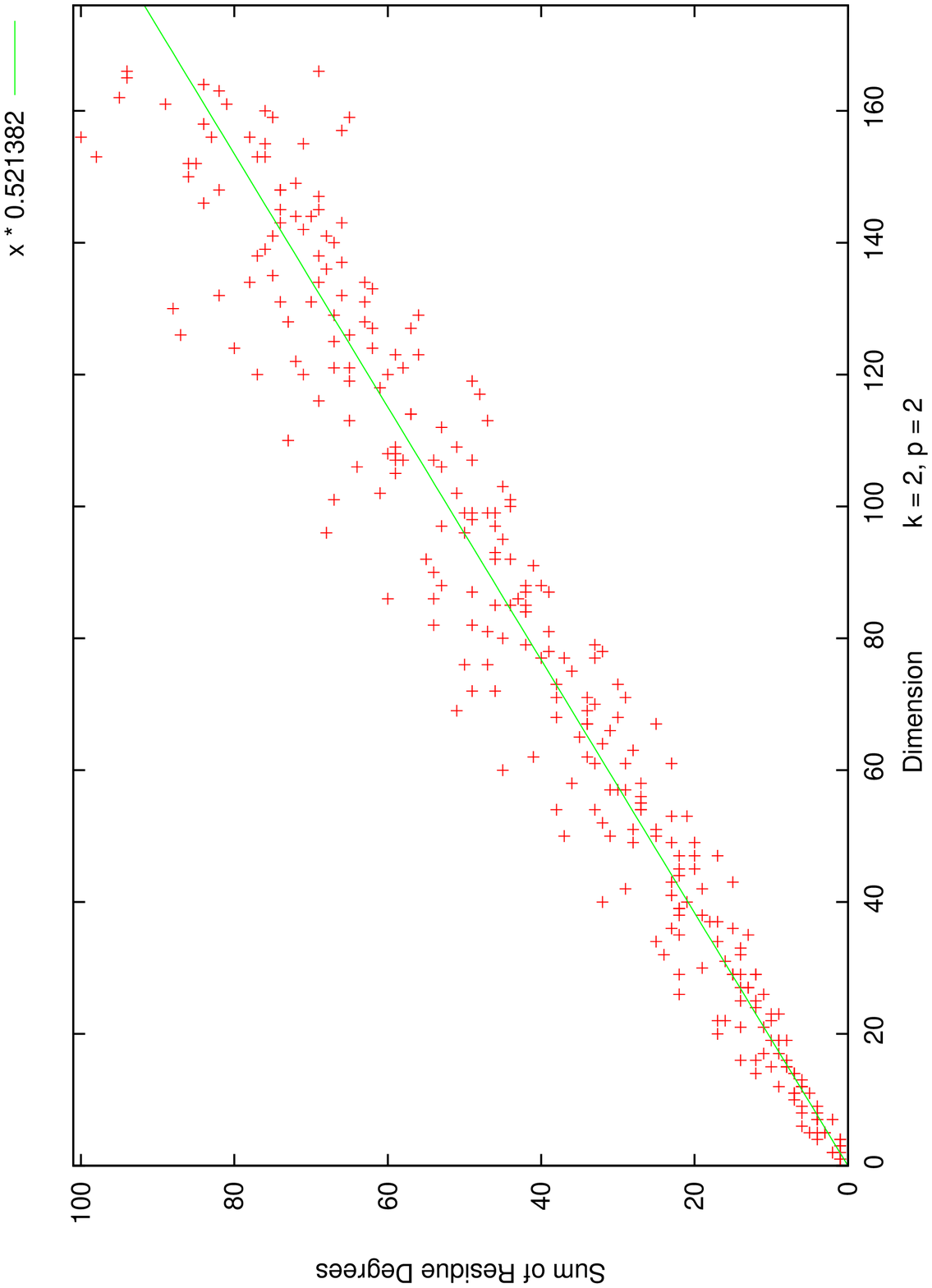}{a:2:2:klein} &
\mplotr{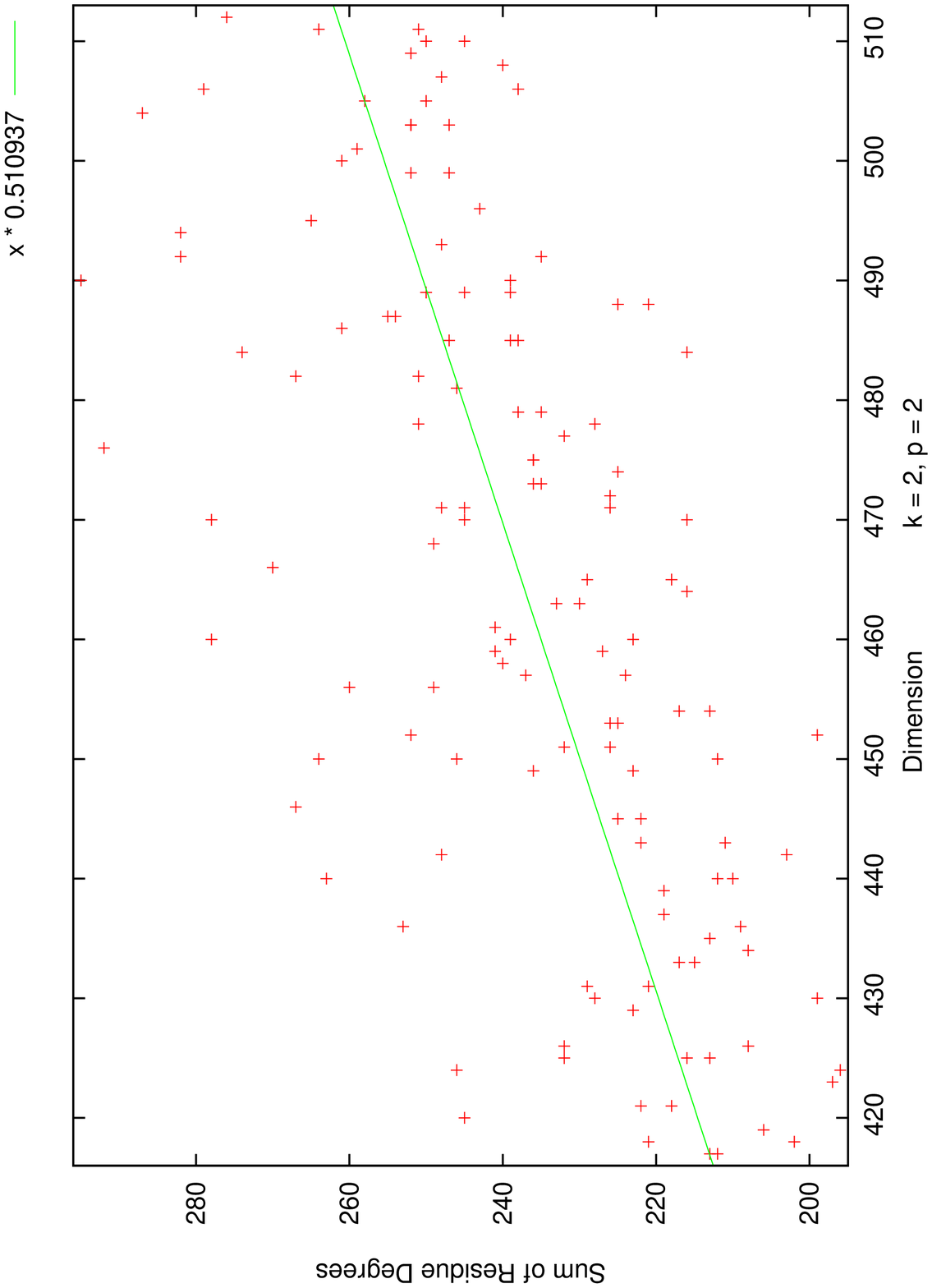}{a:2:2:r1} \\
\mplotr{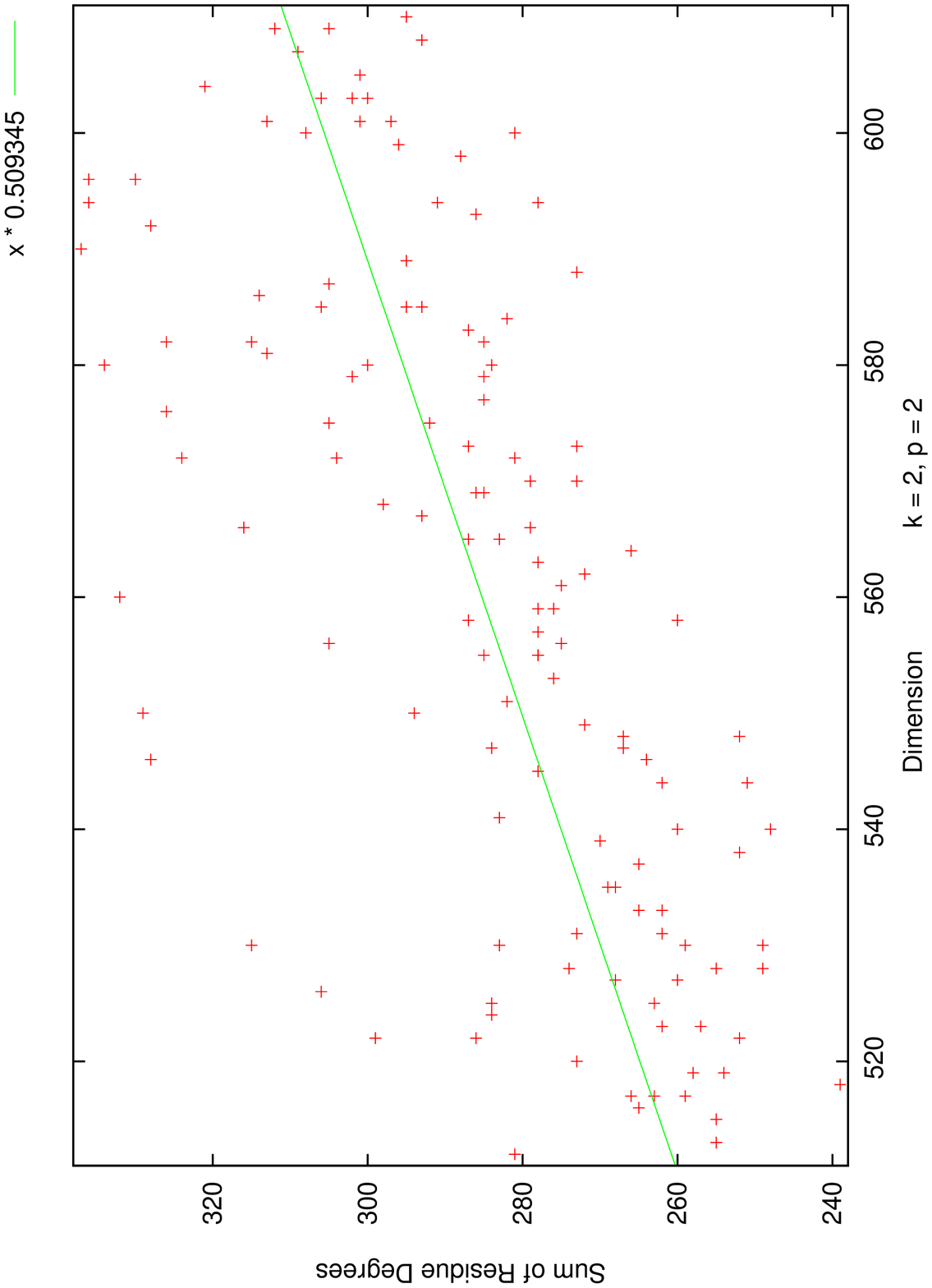}{a:2:2:r2} &
\mplotr{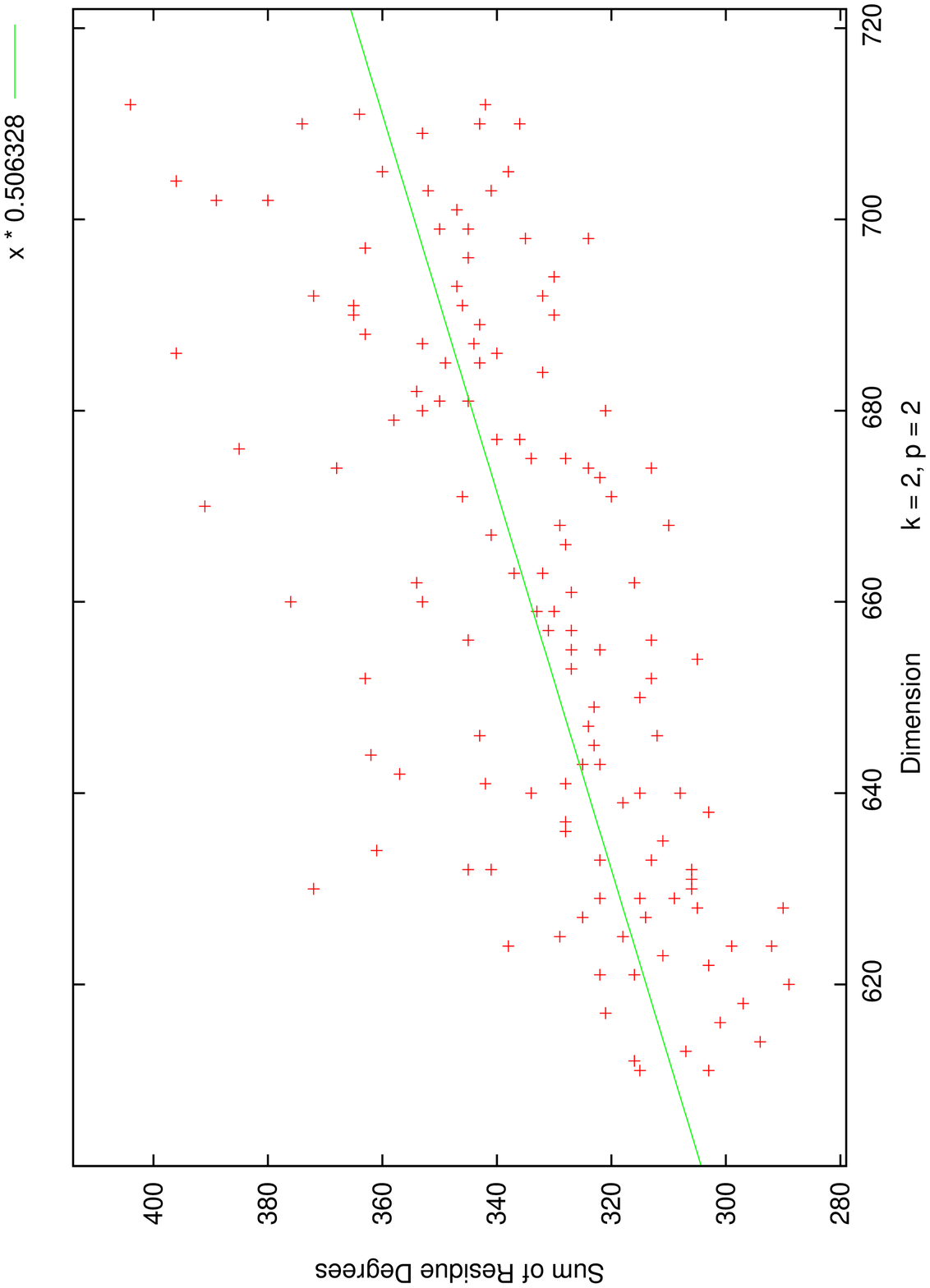}{a:2:2:r3} \\
\mplotr{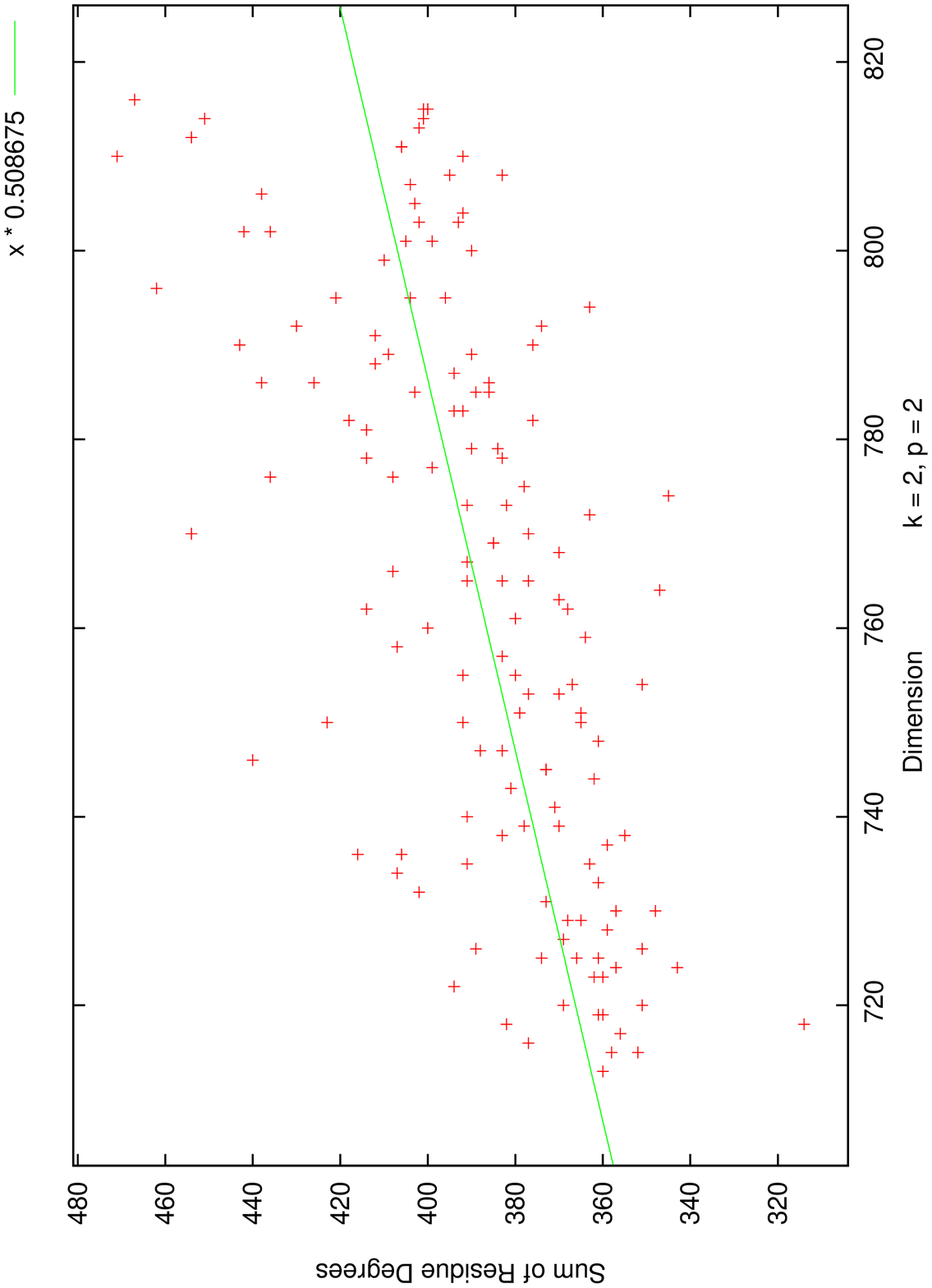}{a:2:2:r4} &
\mplotr{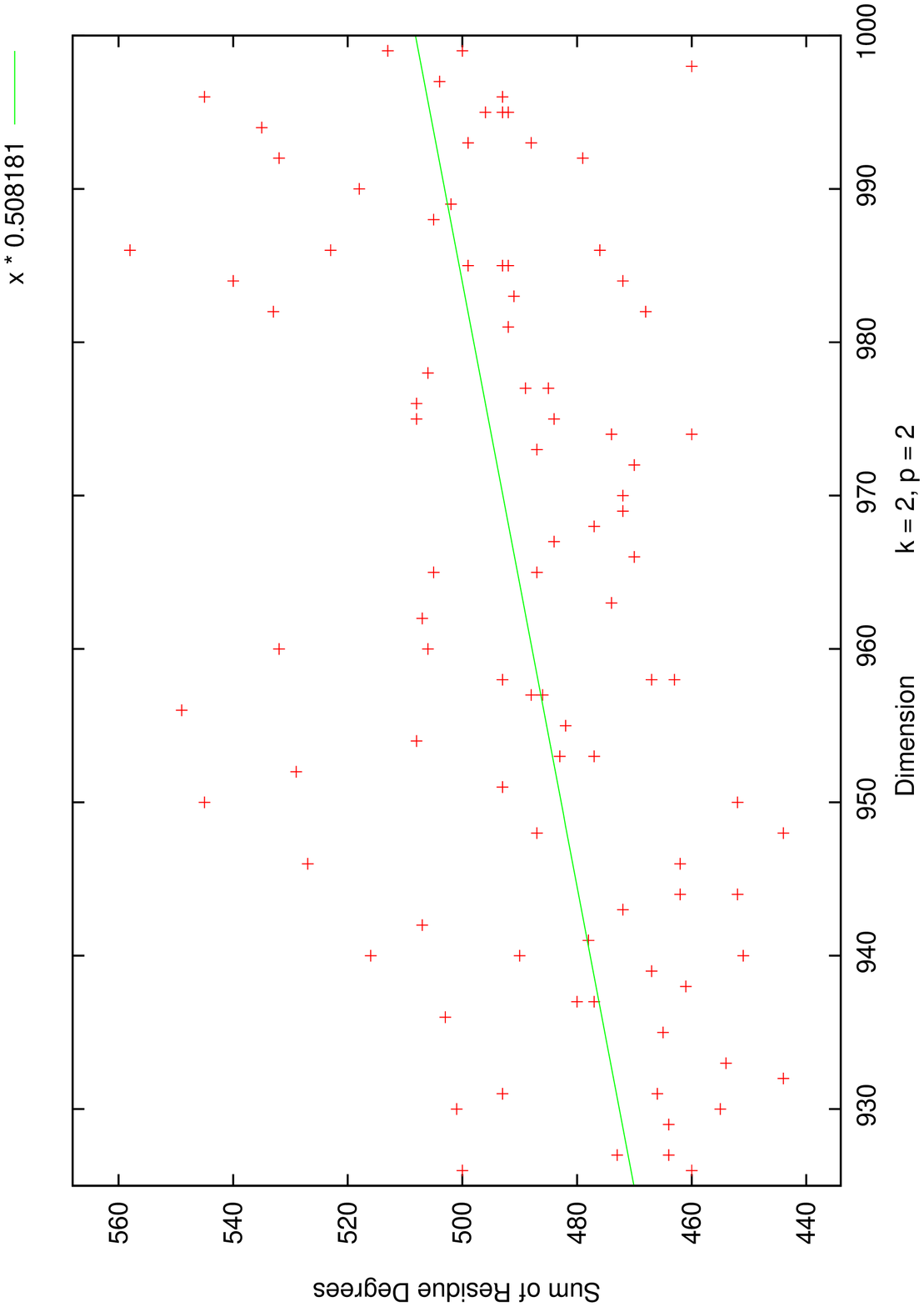}{a:2:2:gross} \\
\end{longtable}

In spite of the very irregular behaviour, it is remarkable that the slope of the best fitting
line through the origin is always just a little bigger than~$\frac{1}{2}$.
We do not attempt to explain this behaviour, but we believe that
the contribution from weight one is not alone responsible for the big difference between
the cases $p=2$ and odd~$p$.
We are led to ask the following question.

\begin{question}\label{question:a-two}
Fix an even weight~$k \ge 2$.
Let $a(N) := a_{k,N}^{(2)}$ and $d(N) := \dim_{\Fbar_2} S_k(N;\Fbar_2)$.

Are there $1 > \alpha \ge \beta > 0$ and constants $C,D > 0$ such that the inequality
$$ \alpha \cdot d(N) + C >a(N) > \beta \cdot d(N) - D$$
holds?
\end{question}

\subsection{Average Residue Degree}\label{sec:b}

We now study the {\em average residue degree},
which we define for given level~$N$, weight~$k$ and prime~$p$ as
$$ b_{N,k}^{(p)} = \frac{1}{\#\Spec(\Tbar_k(N))}\sum_{\fm \in \Spec(\Tbar_k(N))} [\FF_\fm:\FF_p]
= \frac{a_{N,k}^{(p)}}{\#\Spec(\Tbar_k(N))}.$$

We made computations for weight~$2$ and all primes~$p$ less than~$100$, where $N$
runs through the same ranges as previously. We again plot the dimension $d(N)$ on the $x$-axis
and the function~$b_{N,k}^{(p)}$ on the $y$-axis and the green line is again the best fitting
function $\alpha \cdot d(N)$, although we believe that this is not the right function to take (see below).
We again present a selection of prime numbers as before, however, including $p=2$ from the beginning.
The graphs that we leave out have very similar shapes.

\noindent\begin{longtable}{cc}
\mplot{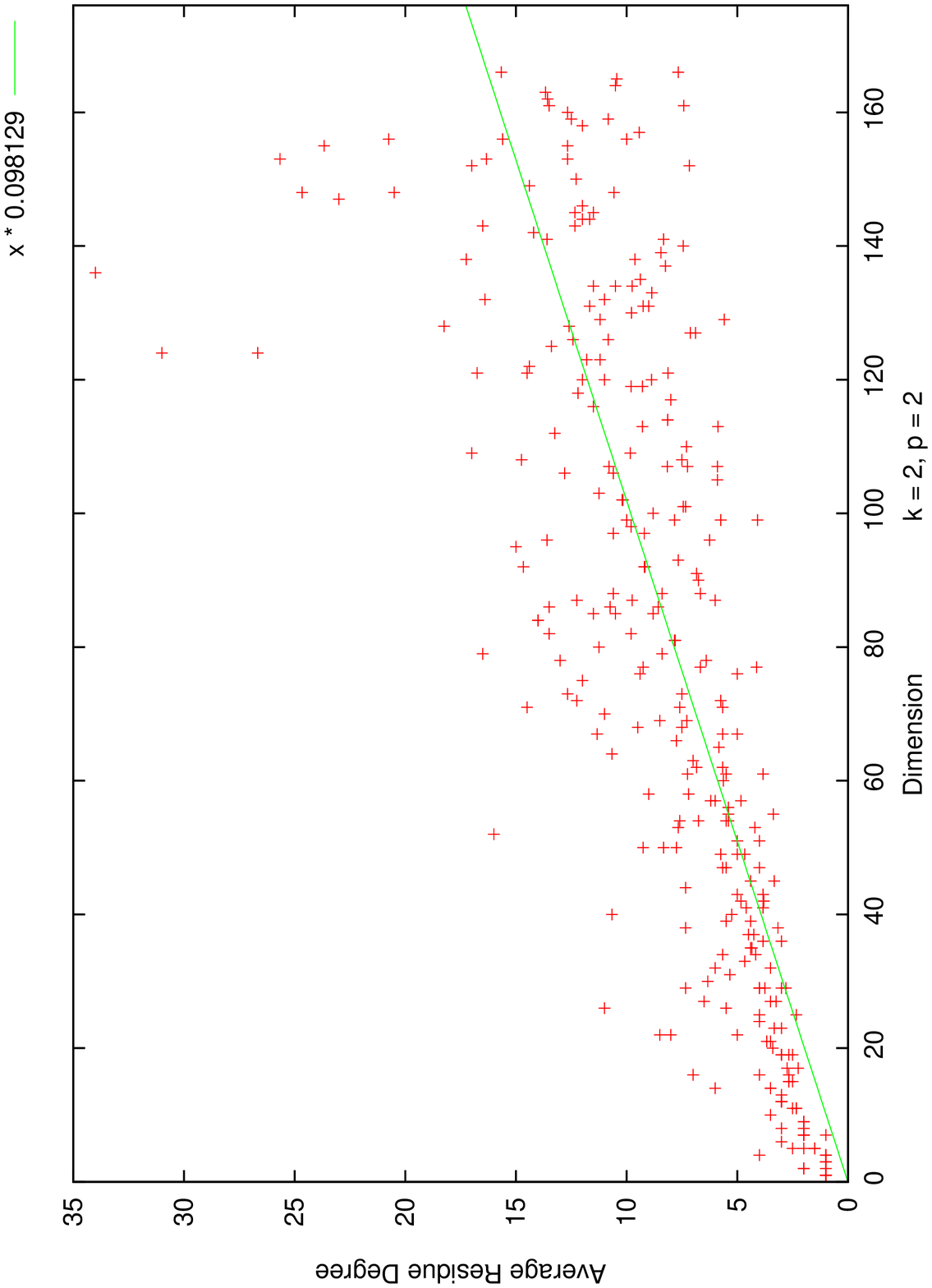} &
\mplot{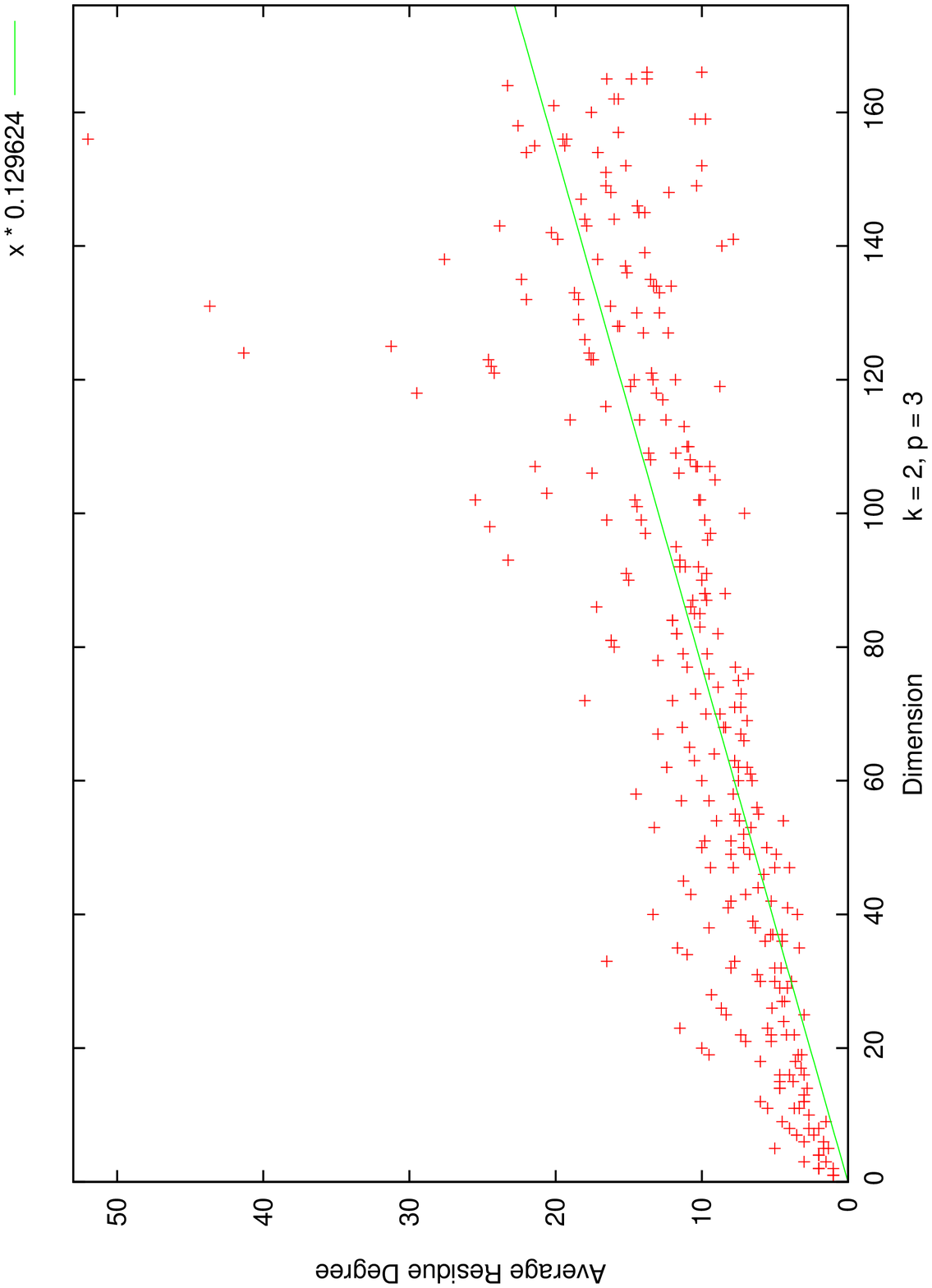} \\
\mplot{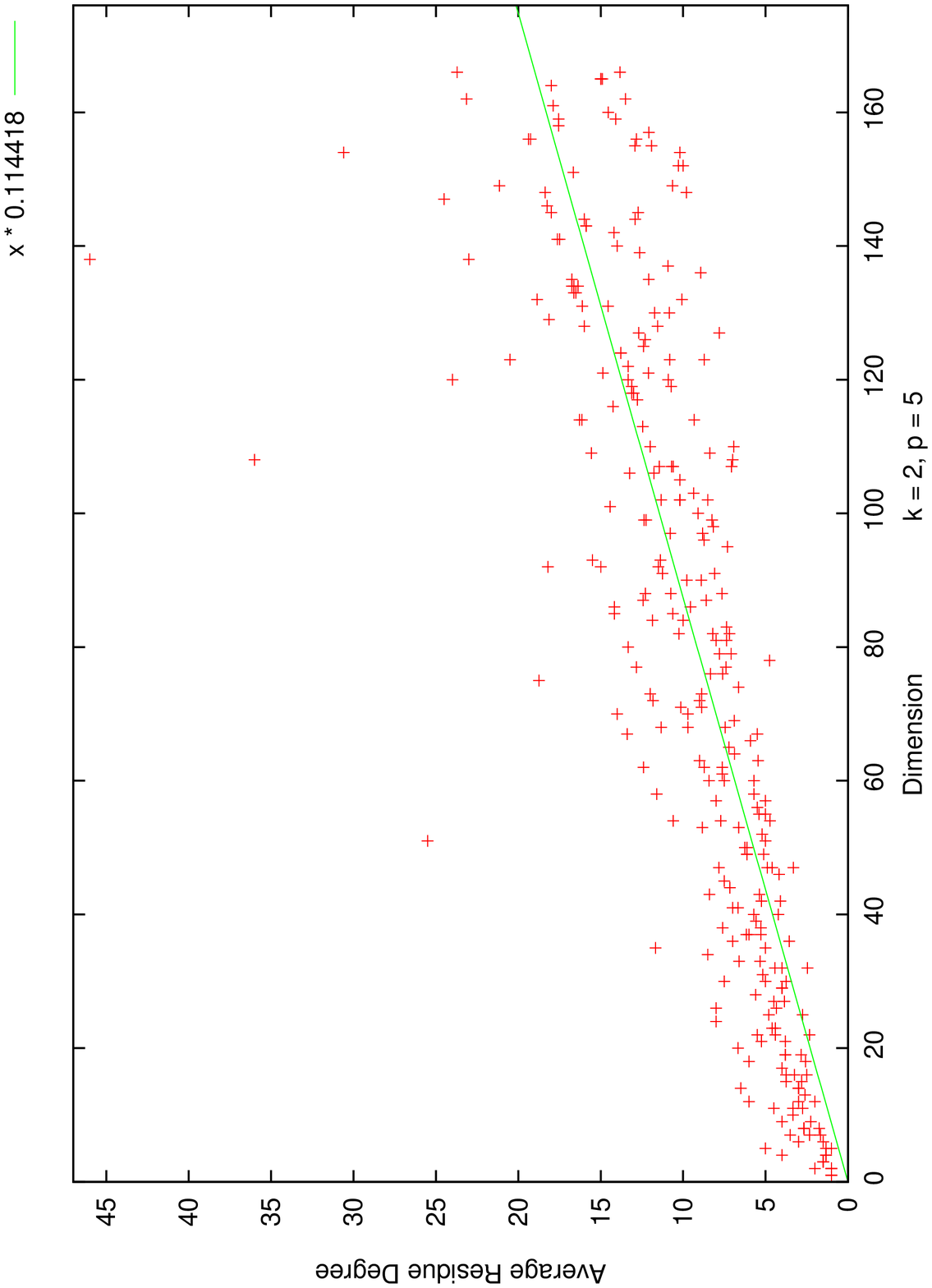} &
\mplot{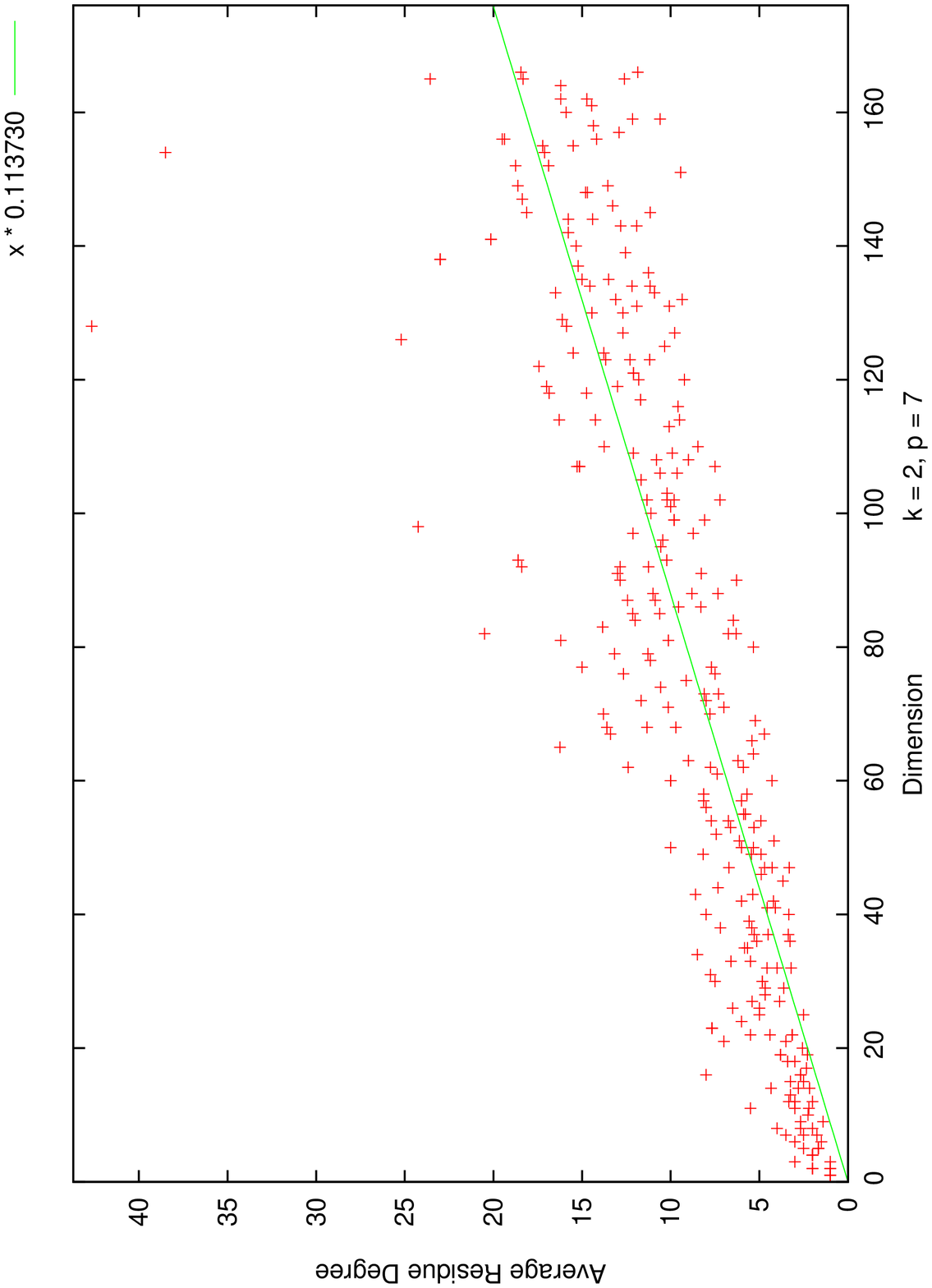} \\
\mplot{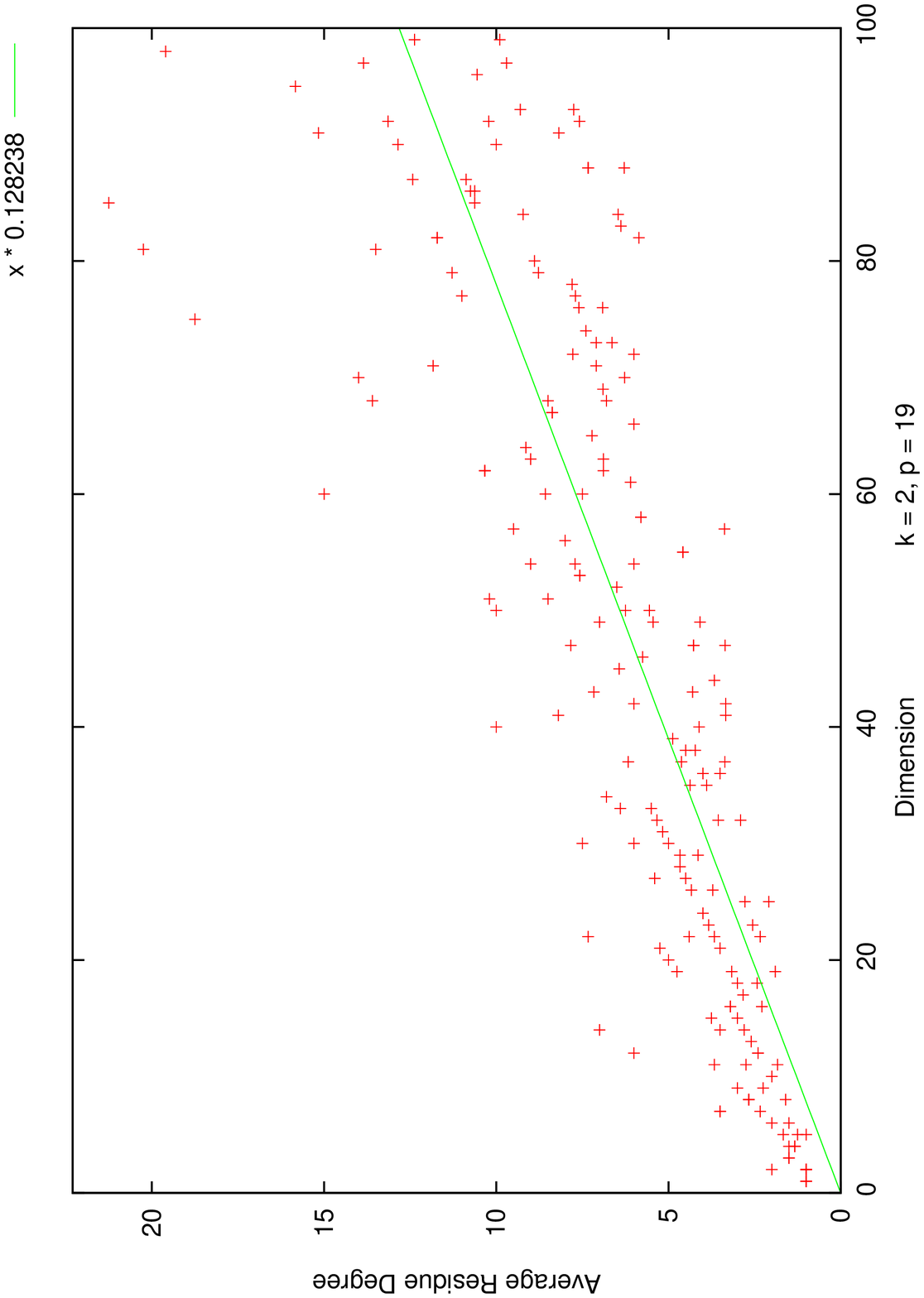}&
\mplot{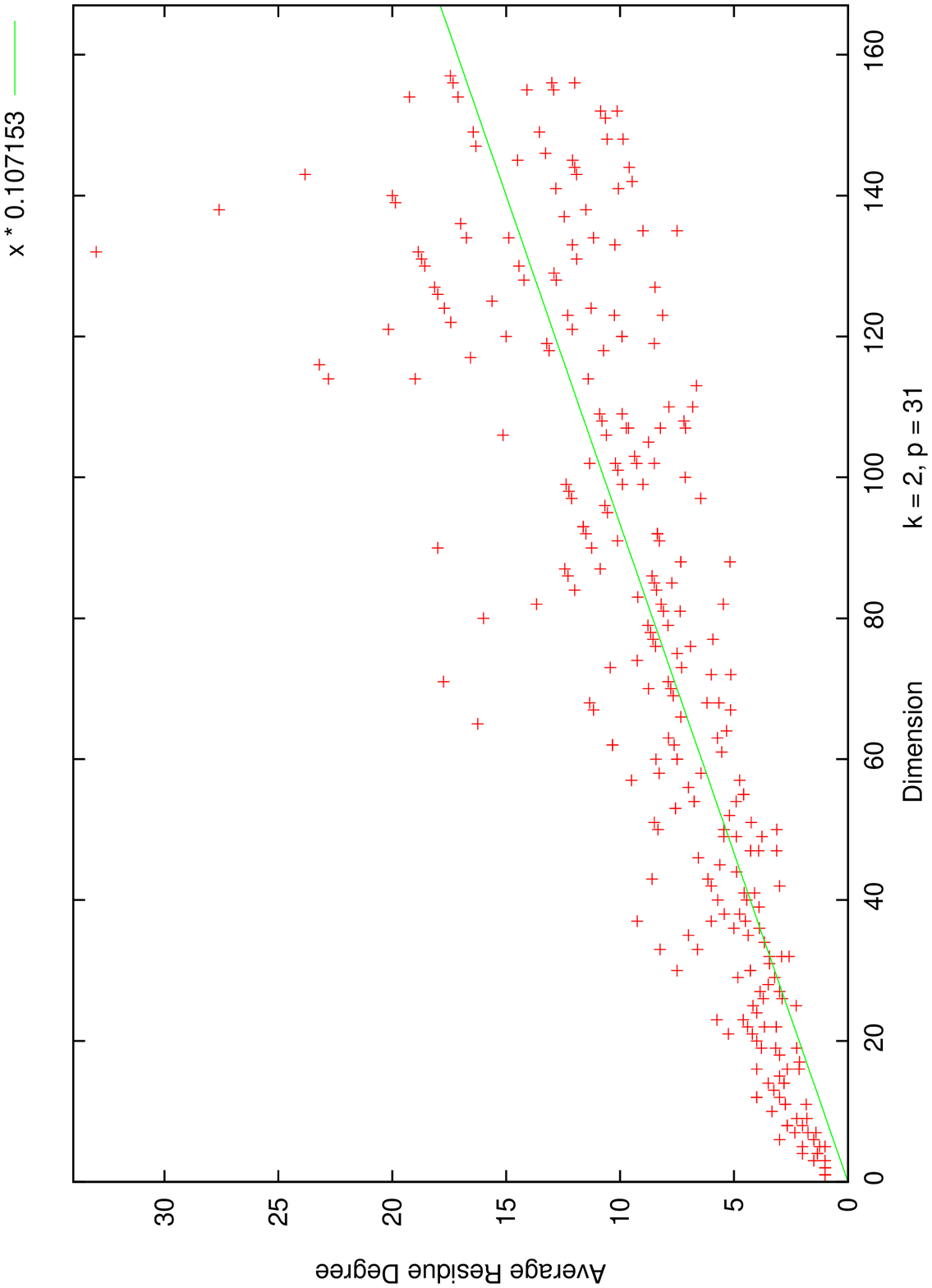} \\
\mplot{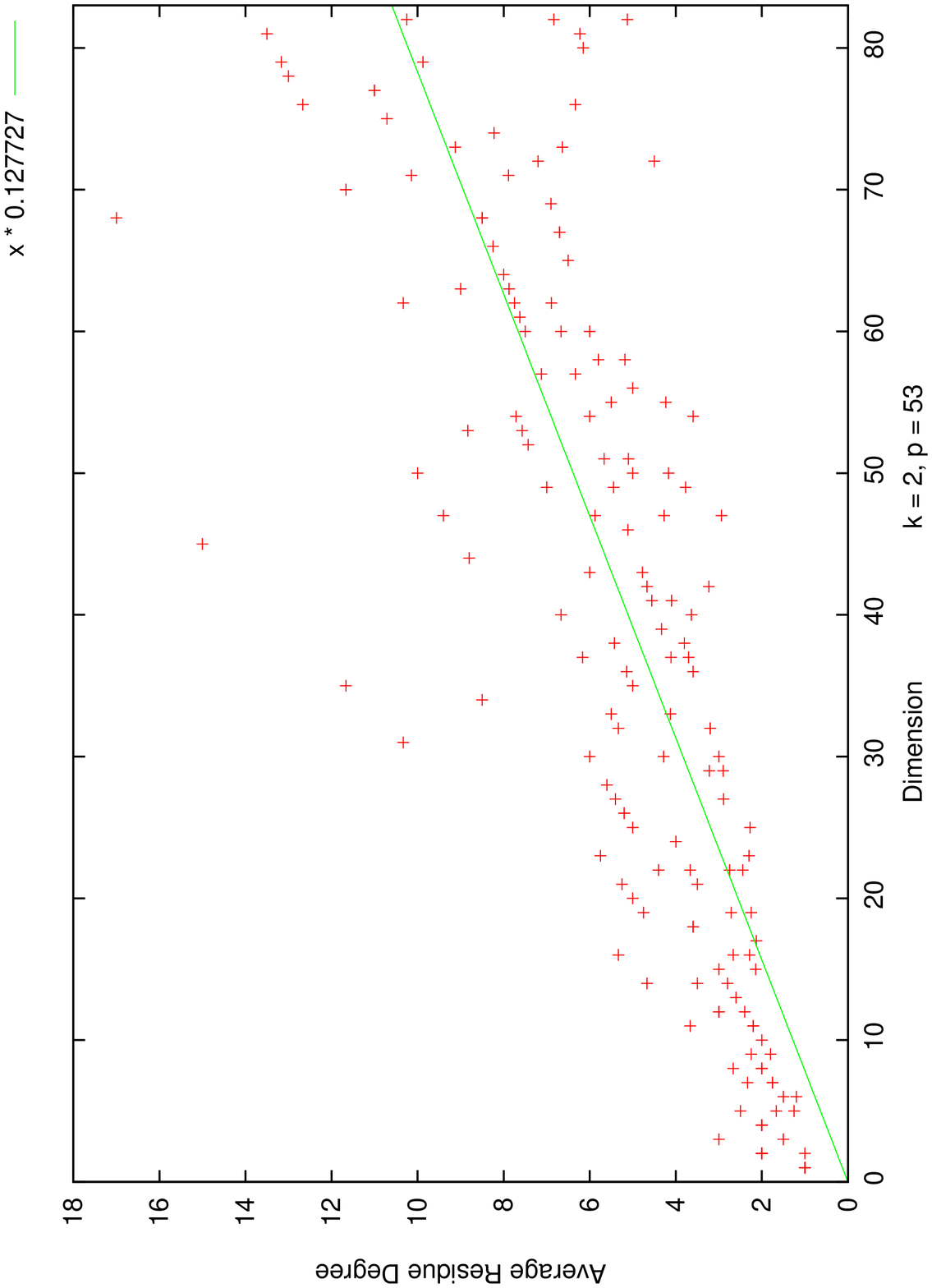}&
\mplot{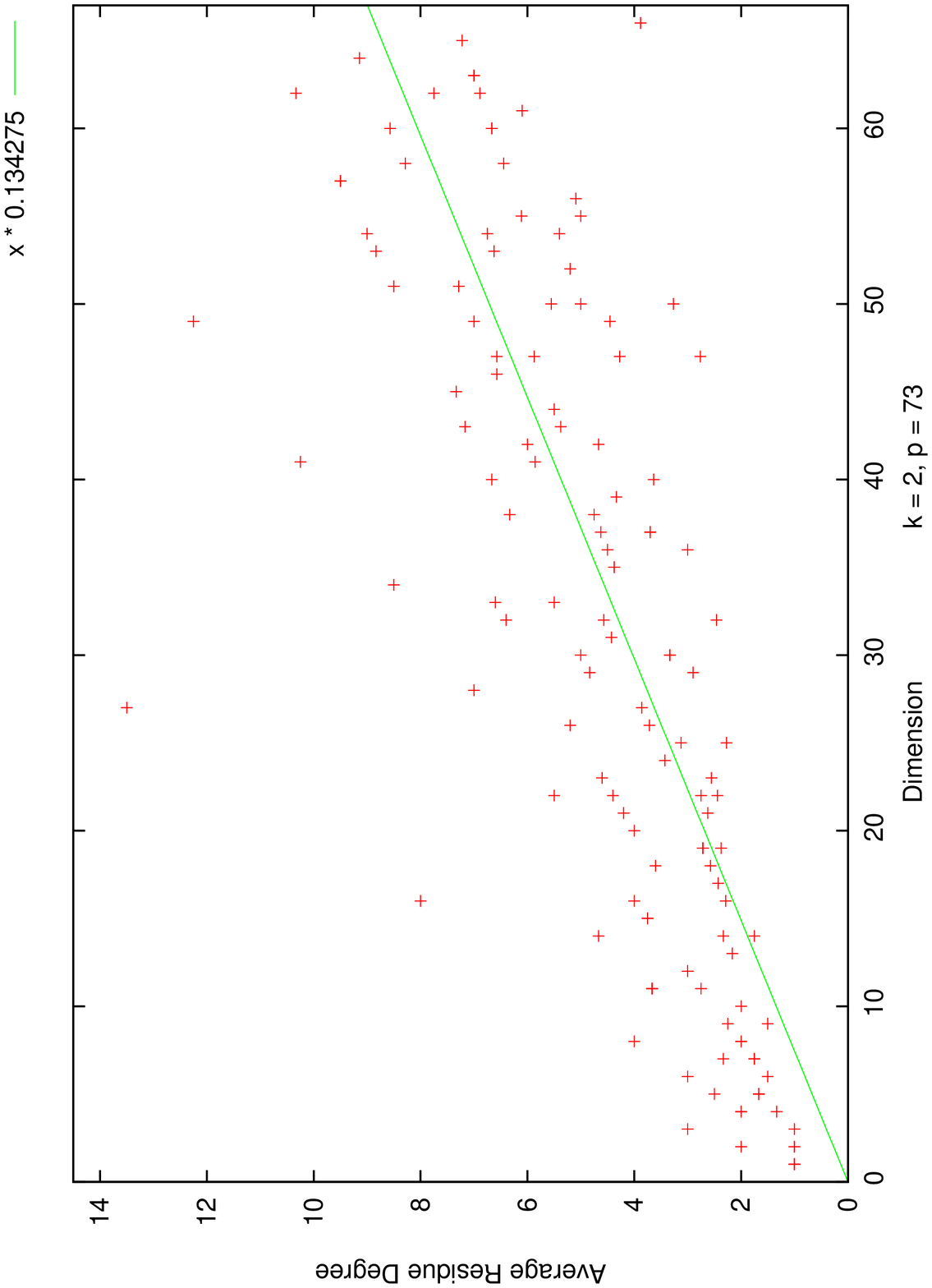} \\
\mplot{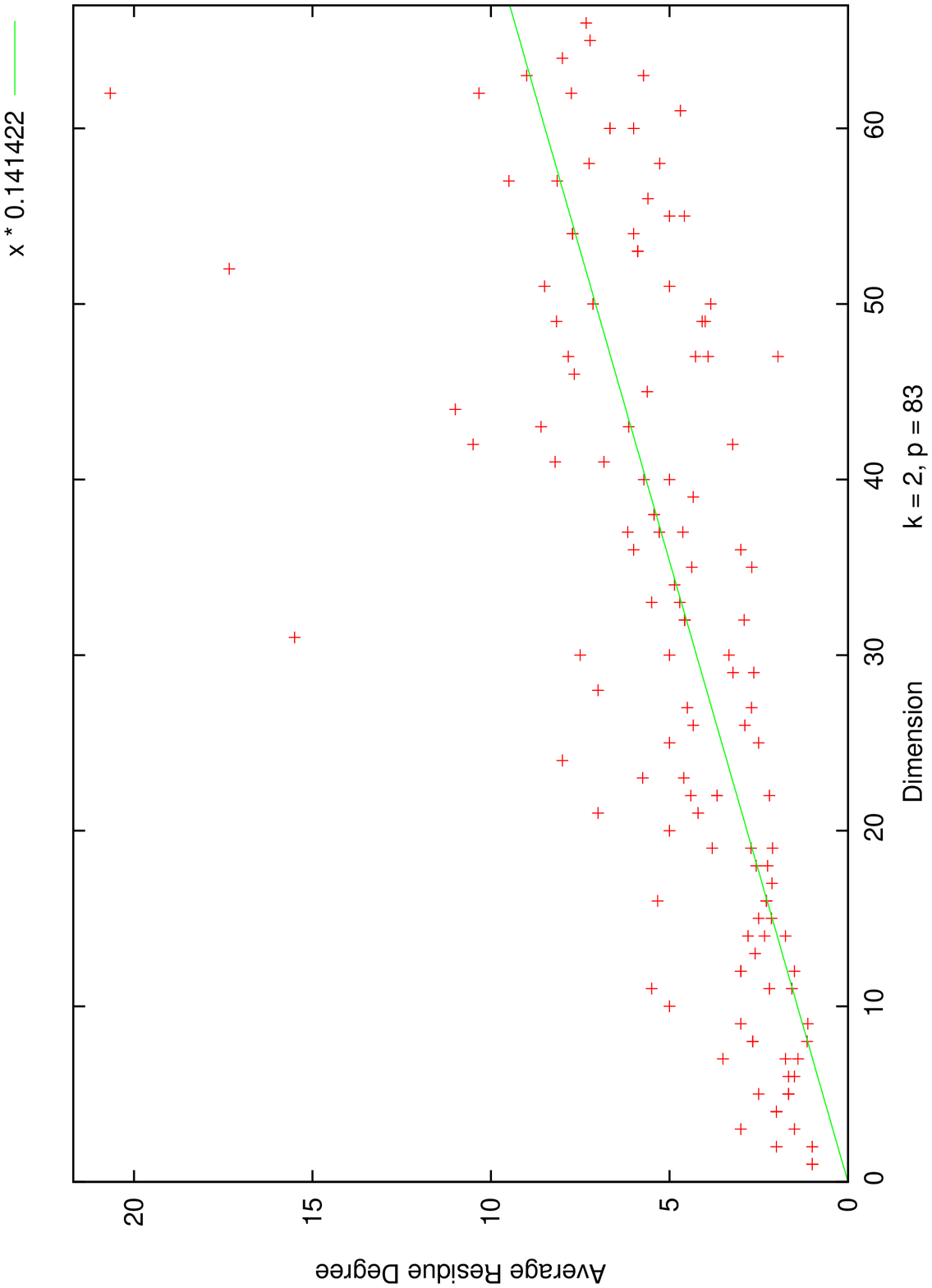} &
\mplot{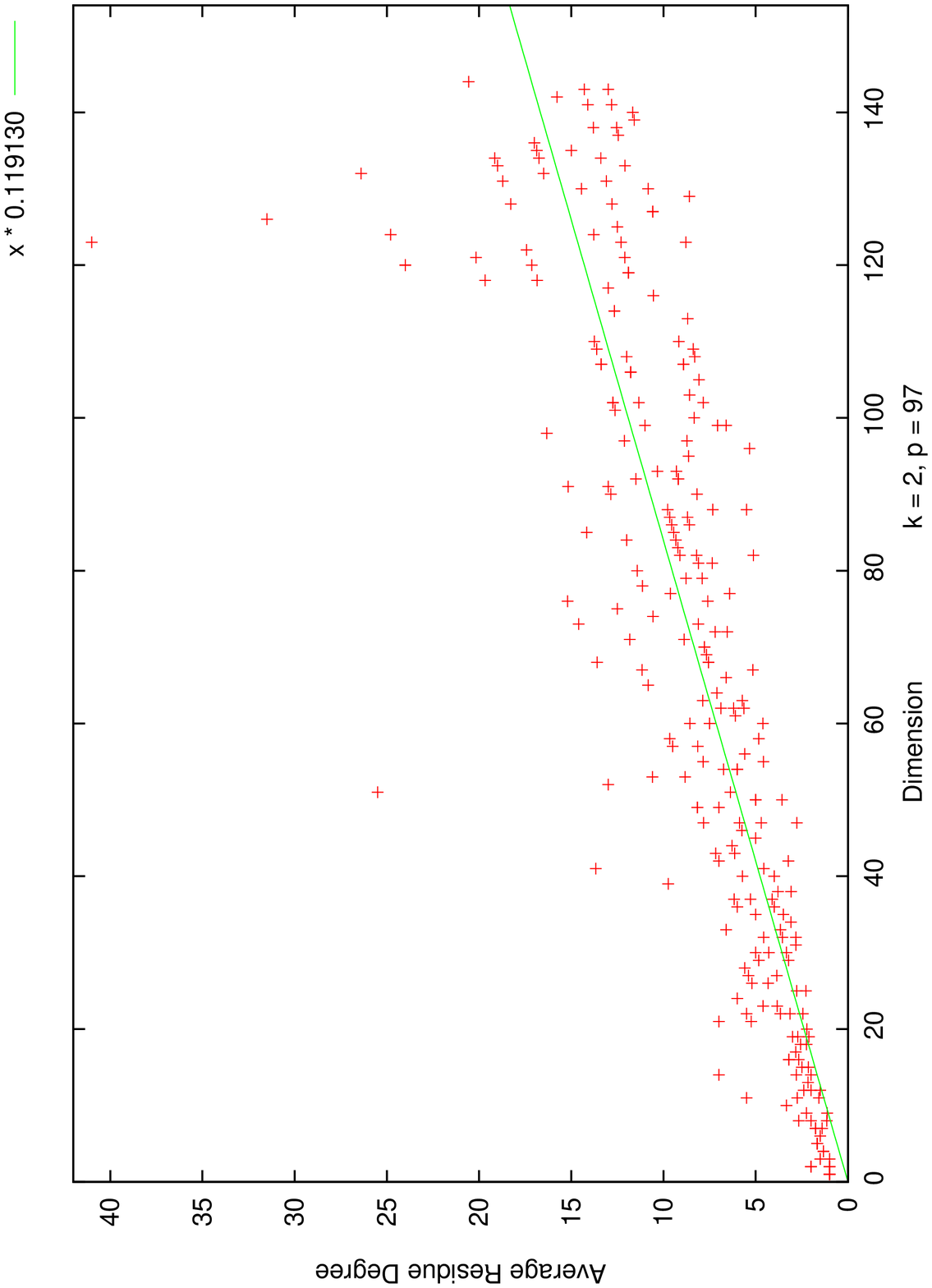} \\
\end{longtable}

Here are again two examples for weight~$4$.
\noindent\begin{longtable}{cc}
\mplot{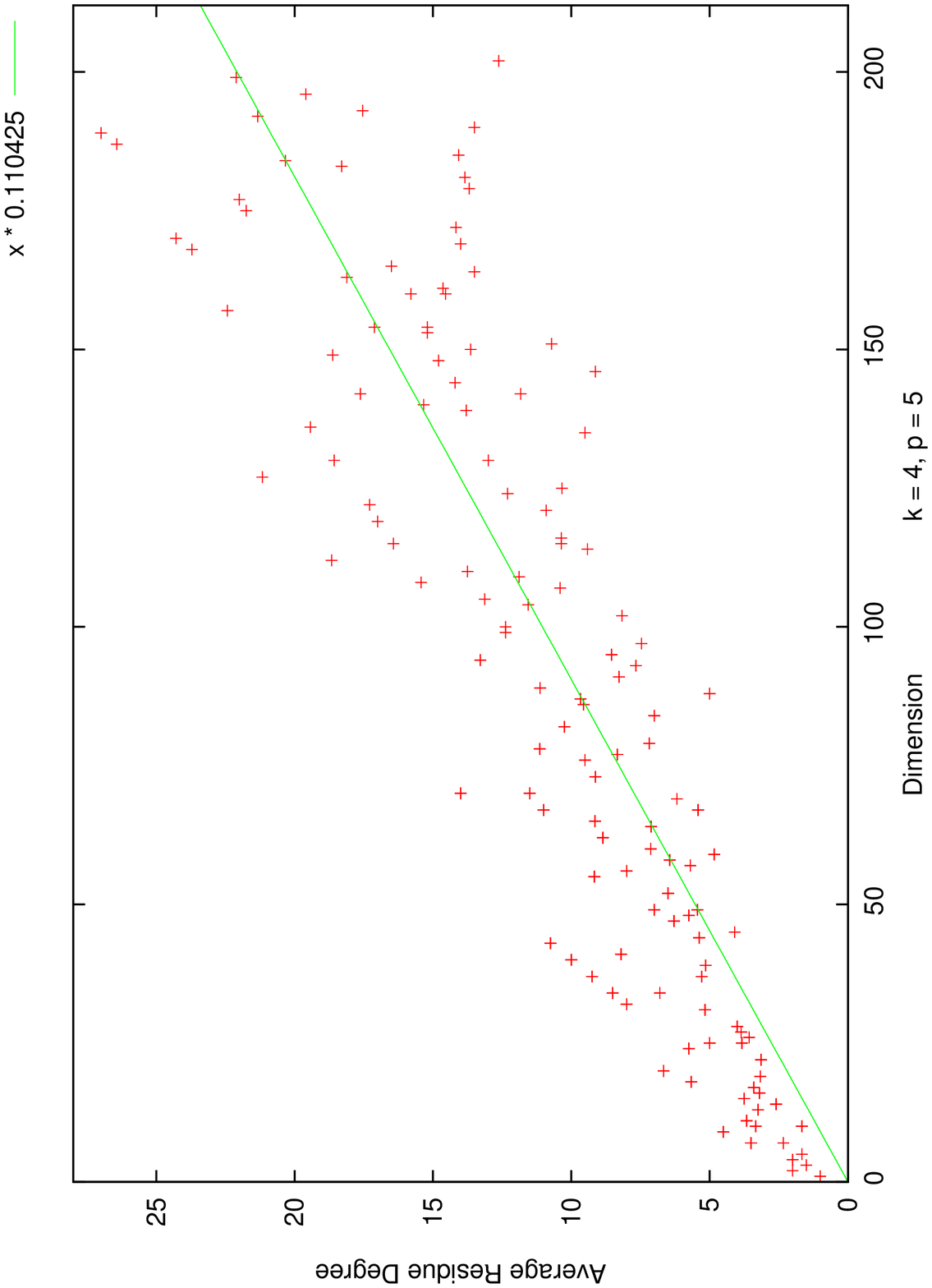} &
\mplot{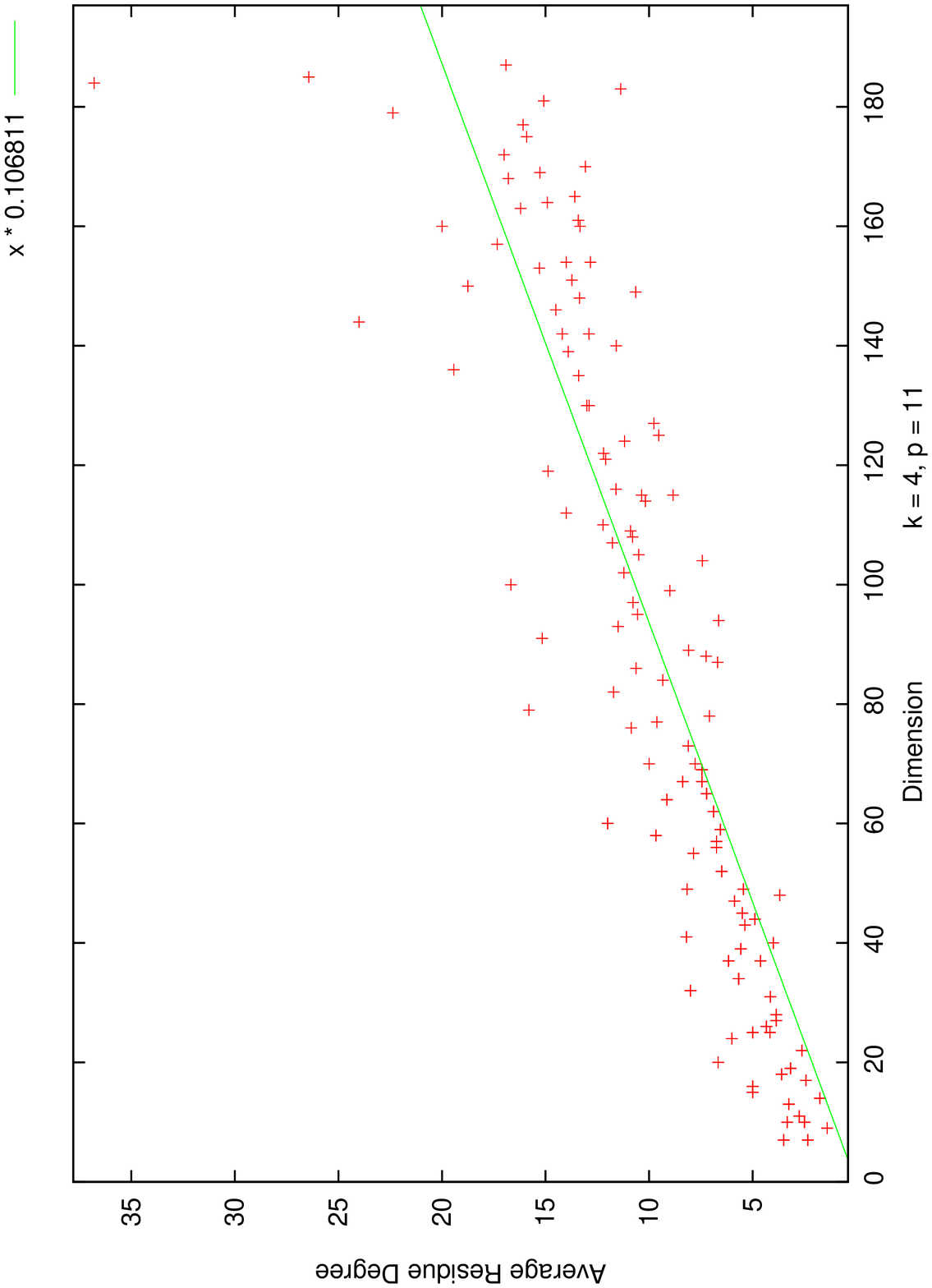} \\
\end{longtable}

Very roughly speaking the data suggest that the average residue degree grows
with the dimension, as is certainly to be expected.
We also conducted a closer analysis for the primes~$2$, $3$ and~$5$.
For $p=2$ we used all primes in different intervals up to~$12000$ and obtained
these plots:

\noindent\begin{longtable}{cc}
\mplot{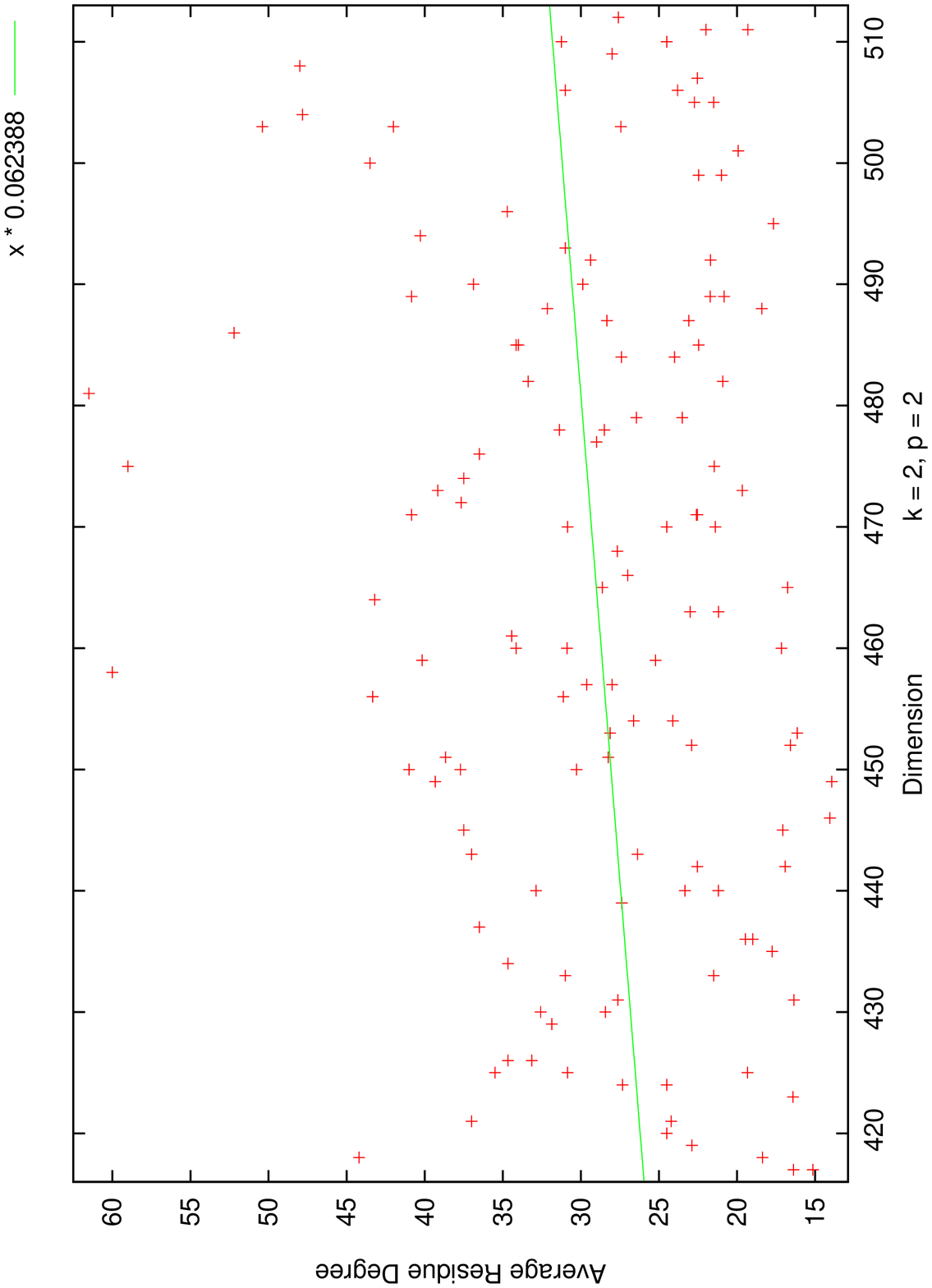} &
\mplot{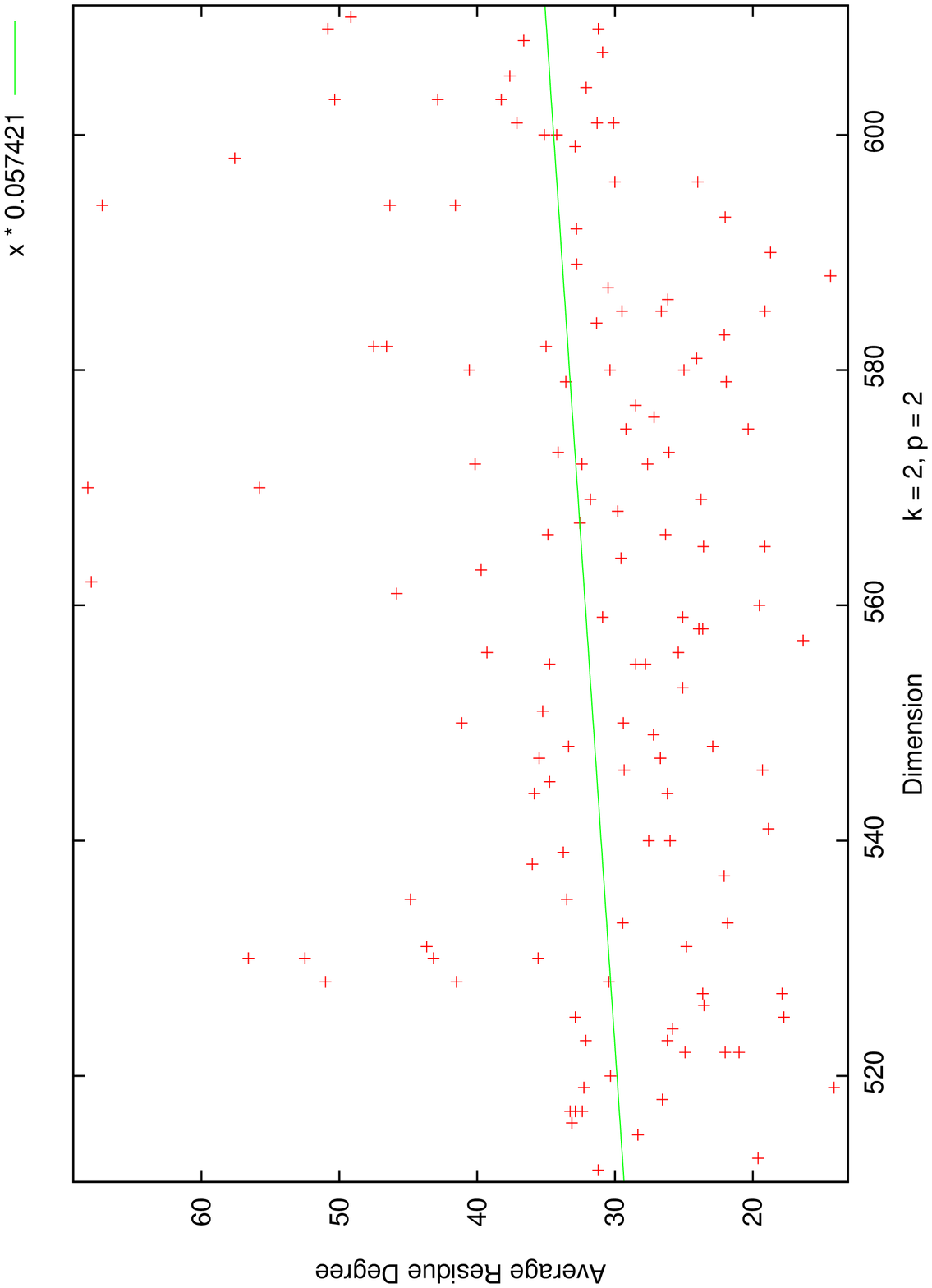} \\
\mplot{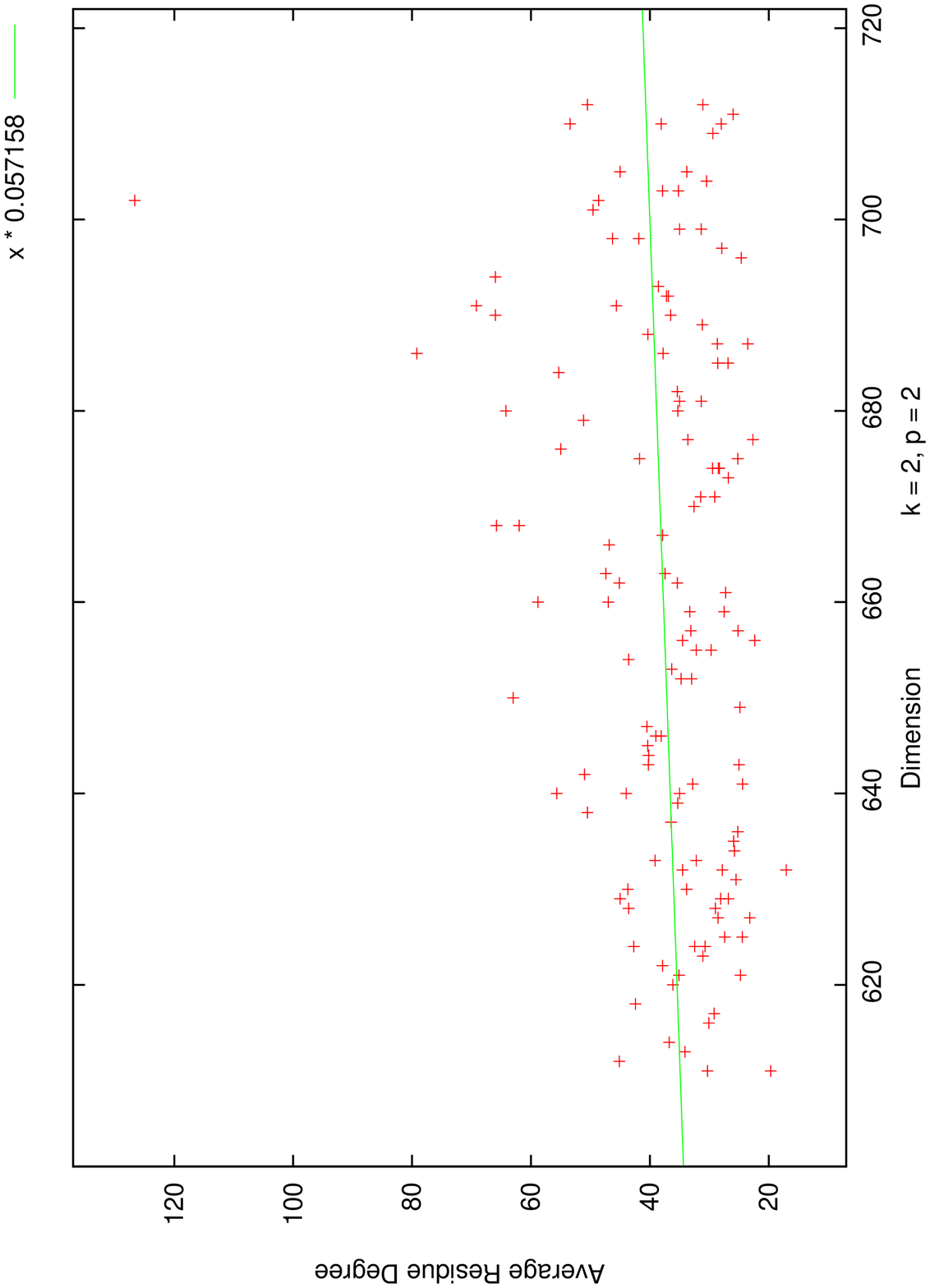} &
\mplot{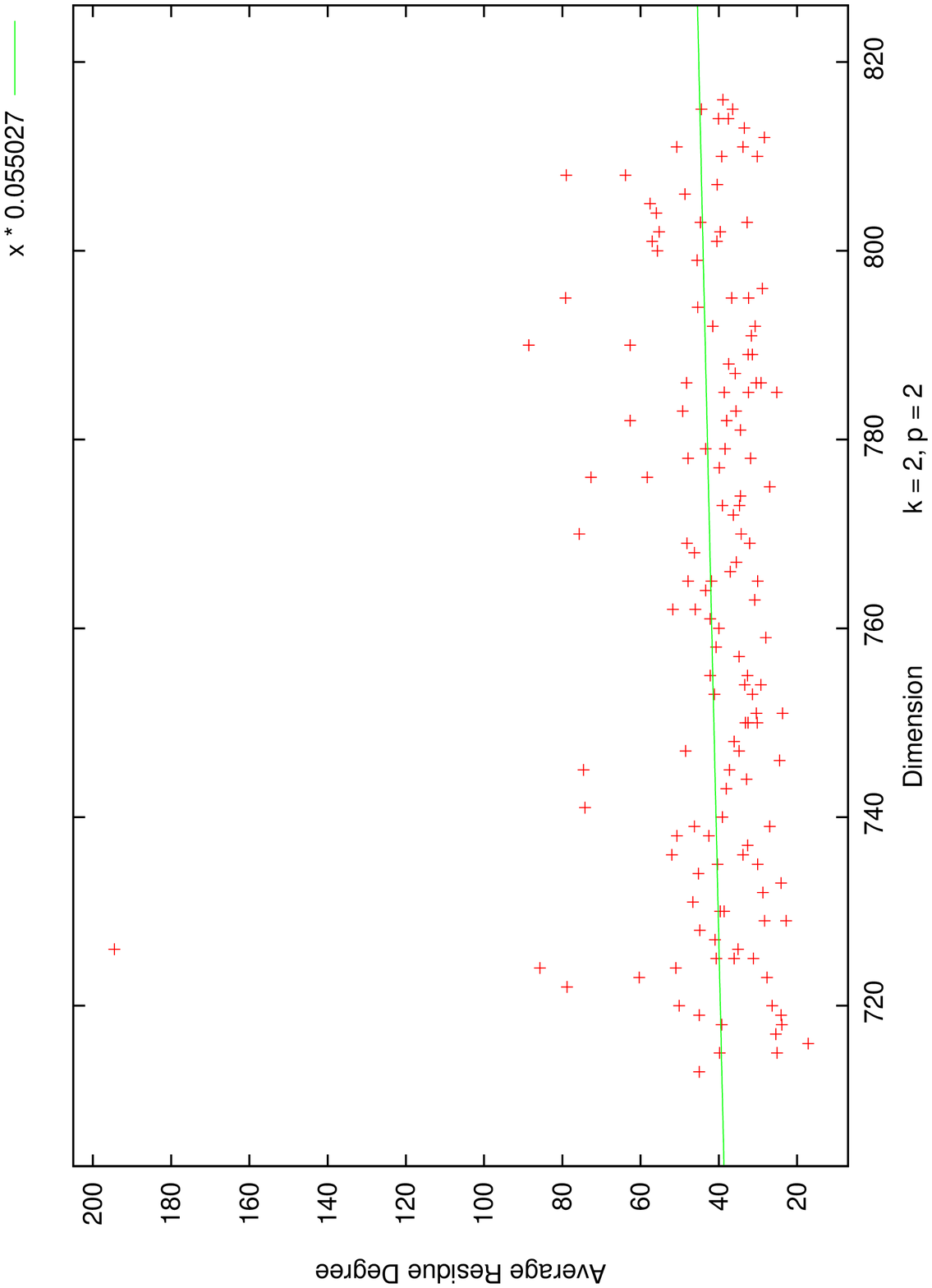} \\
\mplot{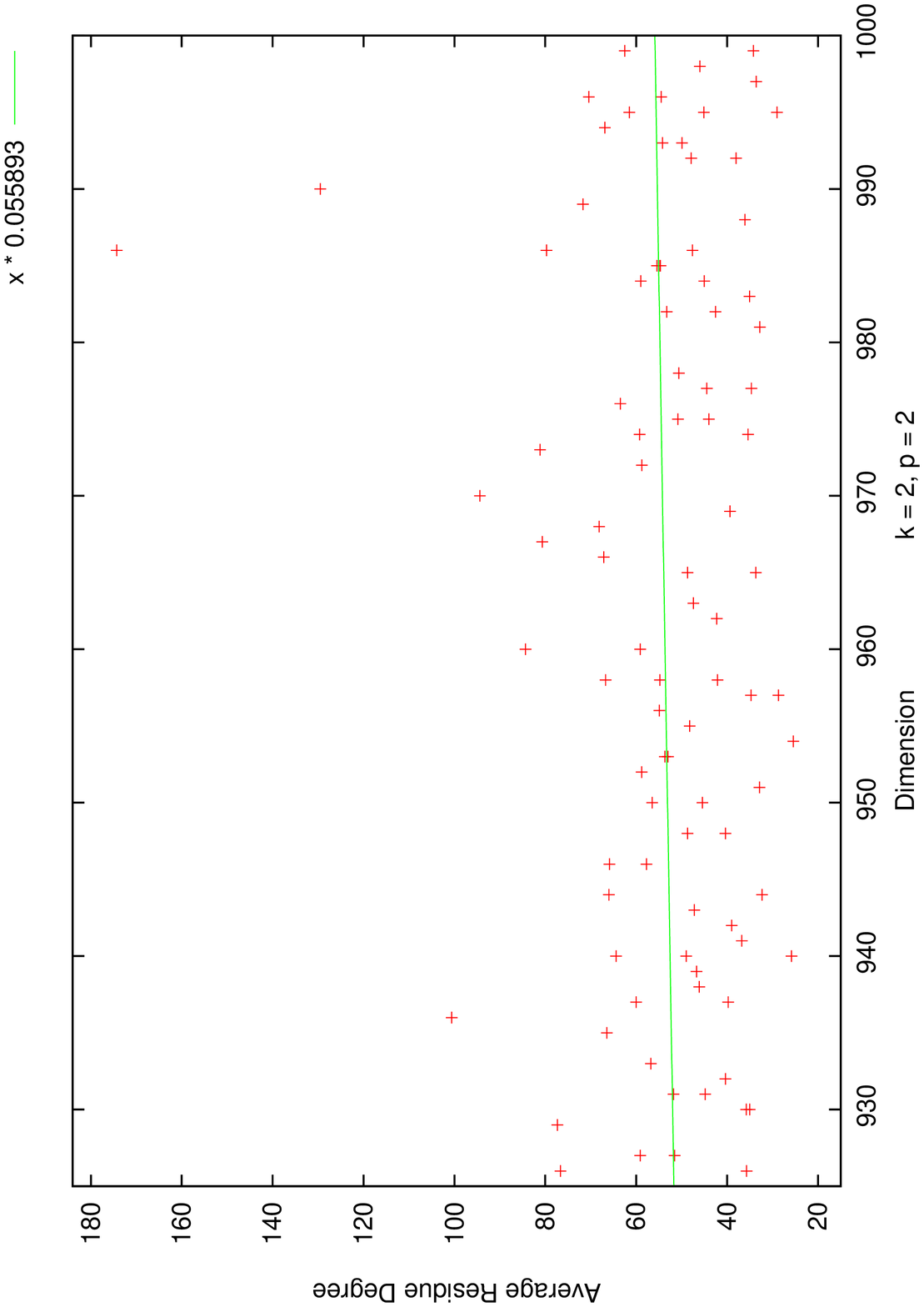} \\
\end{longtable}

Here are the plots for $p=3,5$ and the primes between $3000$ and $10009$ subdivided into four intervals.

\noindent\begin{longtable}{cc}
\mplot{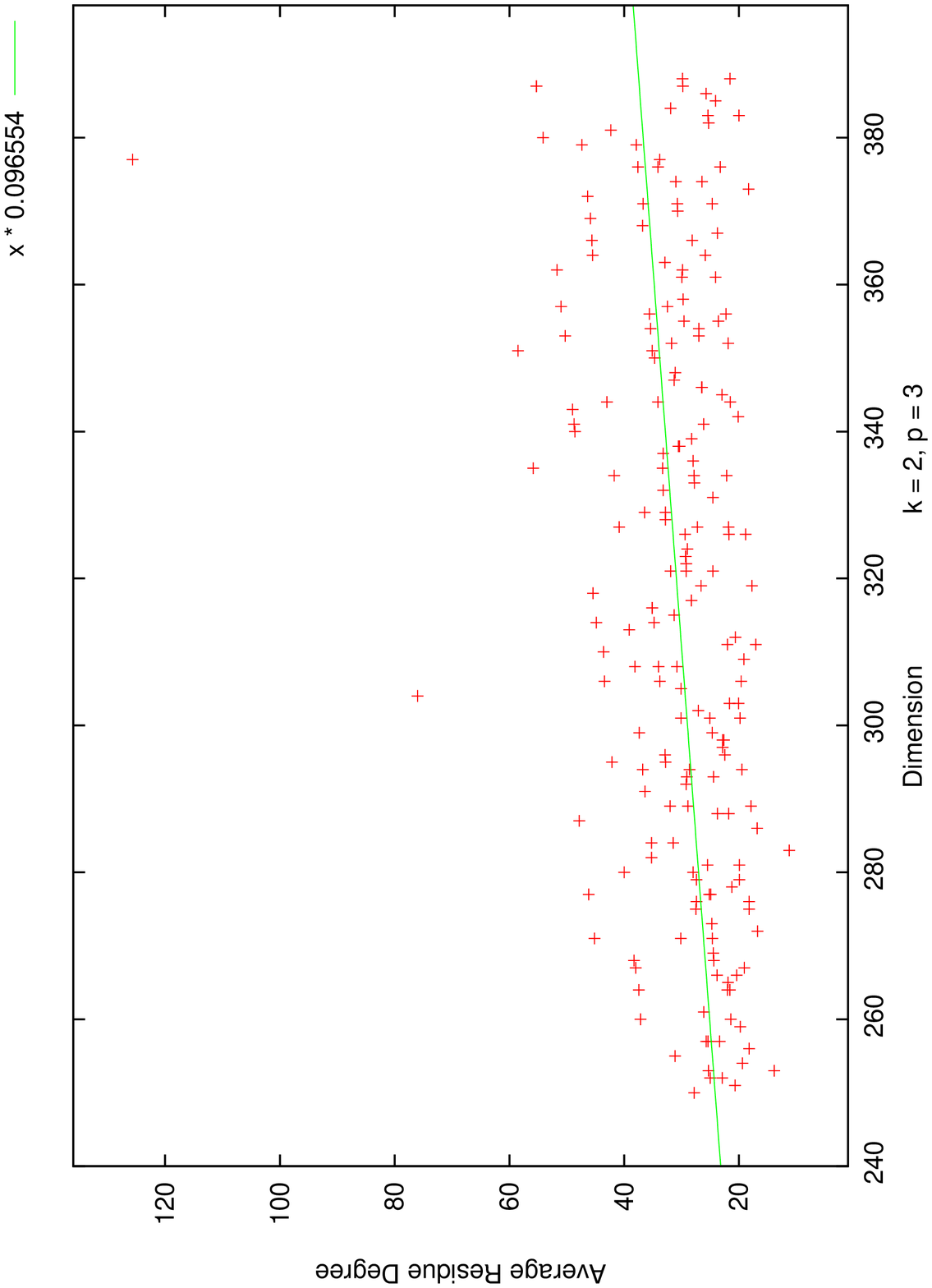} &
\mplot{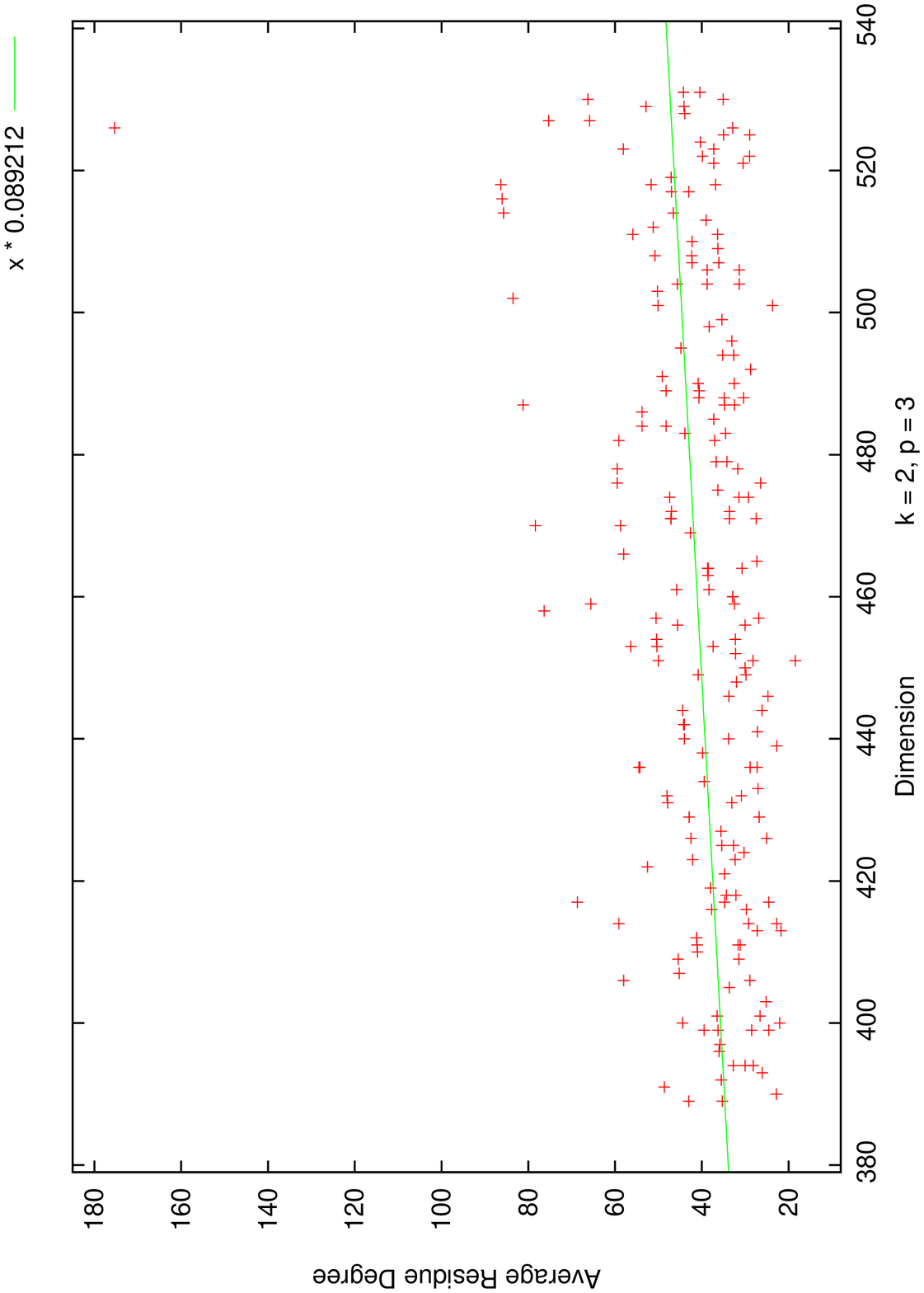} \\
\mplot{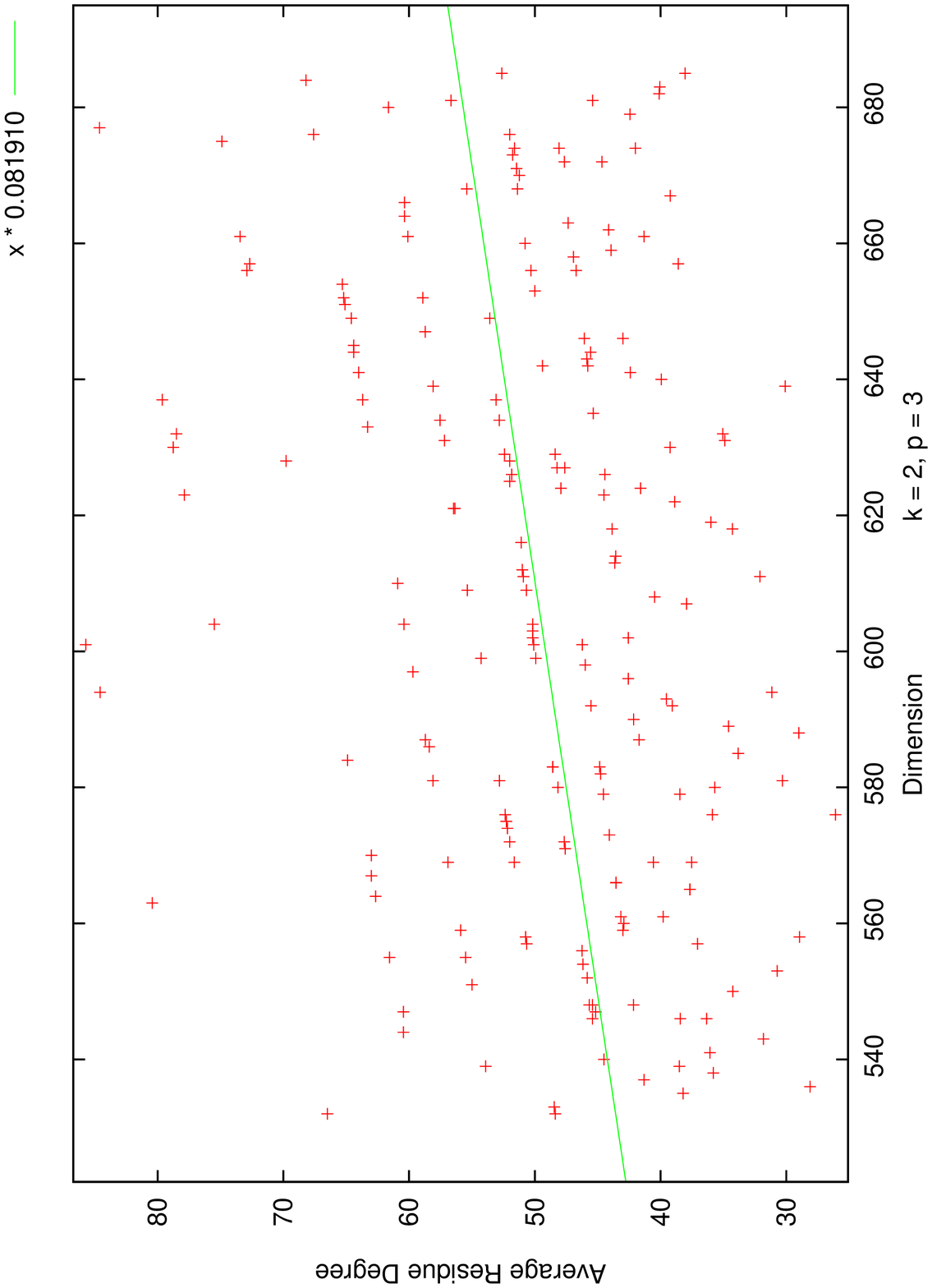} &
\mplot{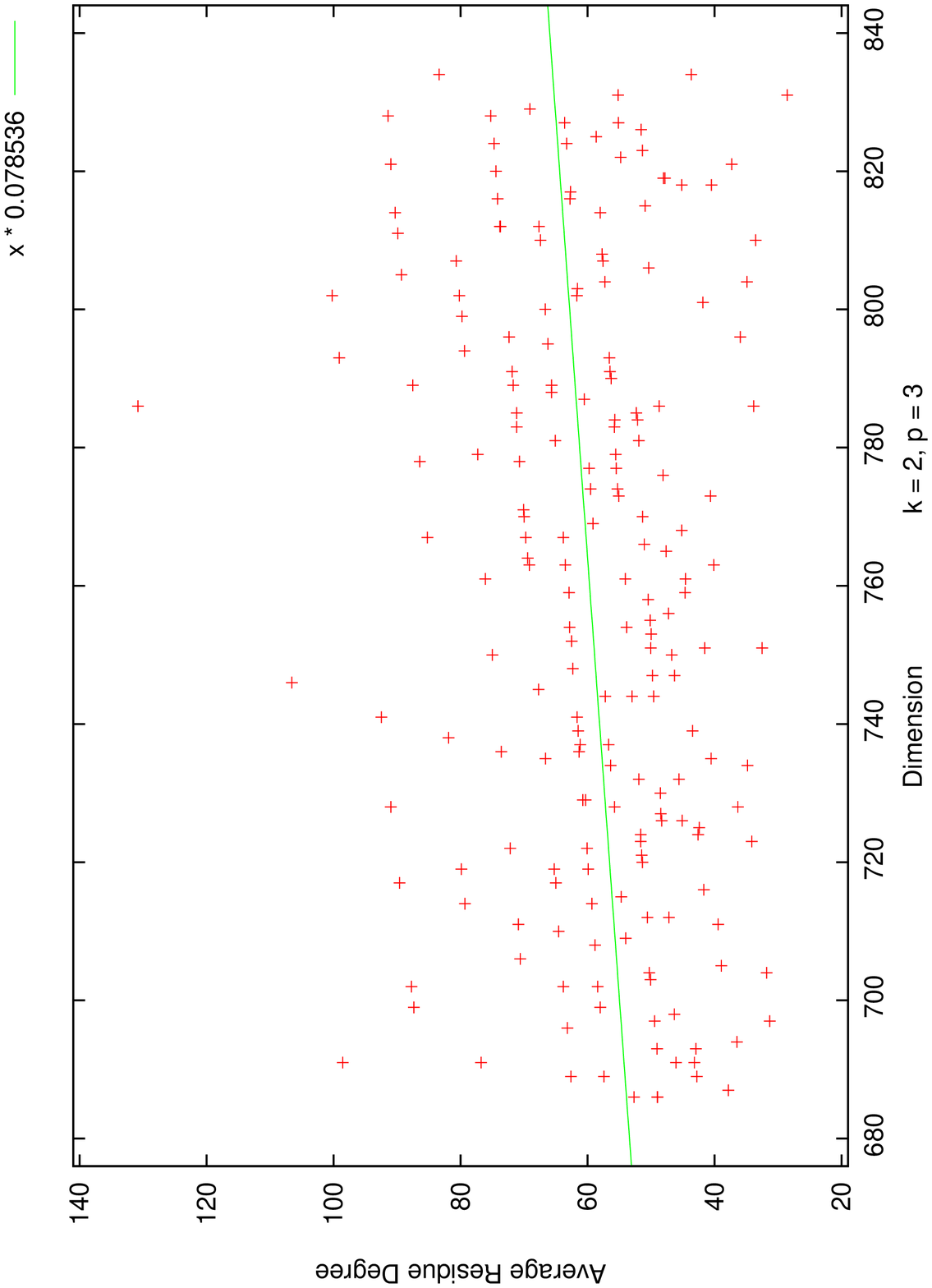} \\
\mplot{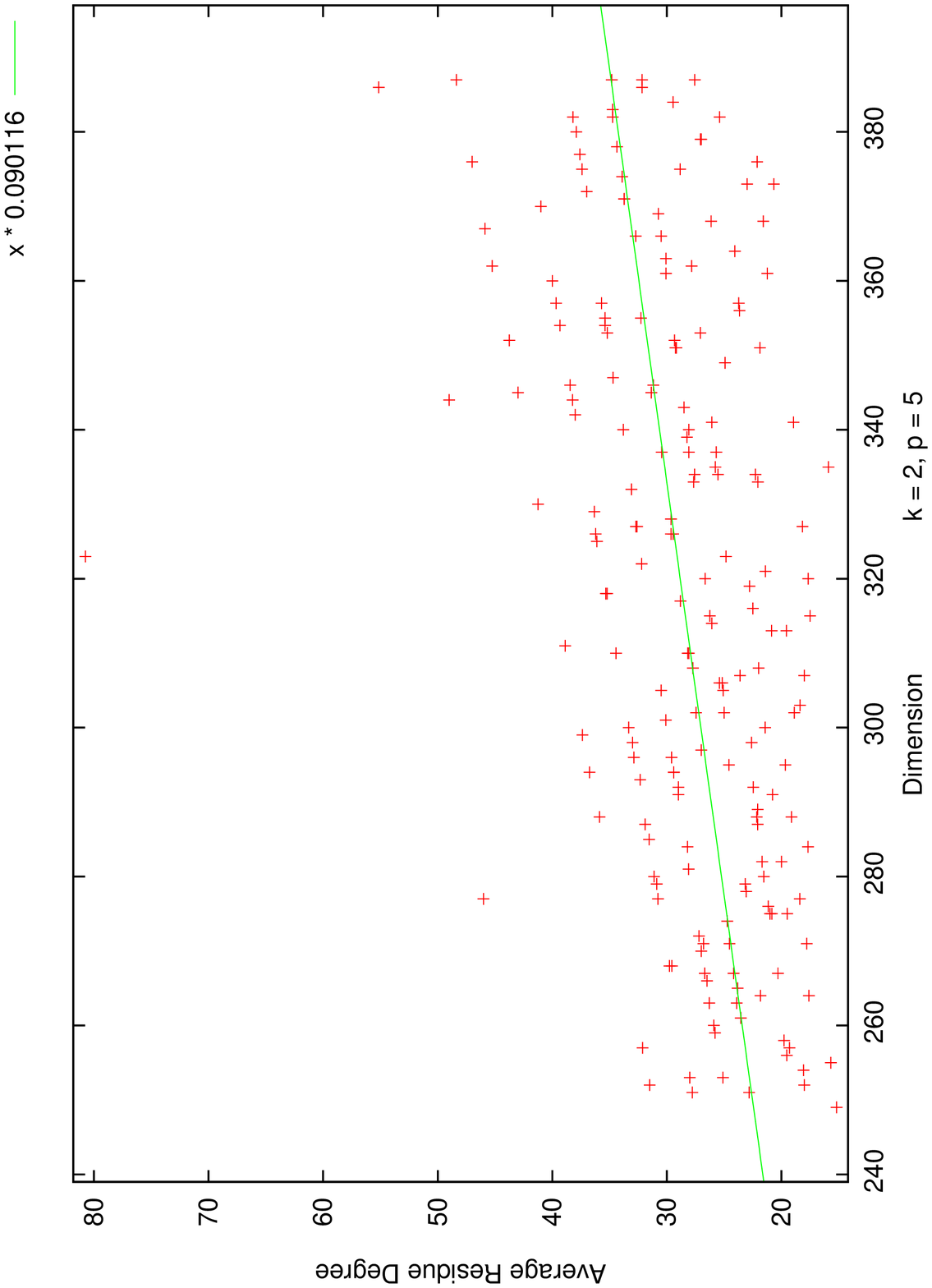} &
\mplot{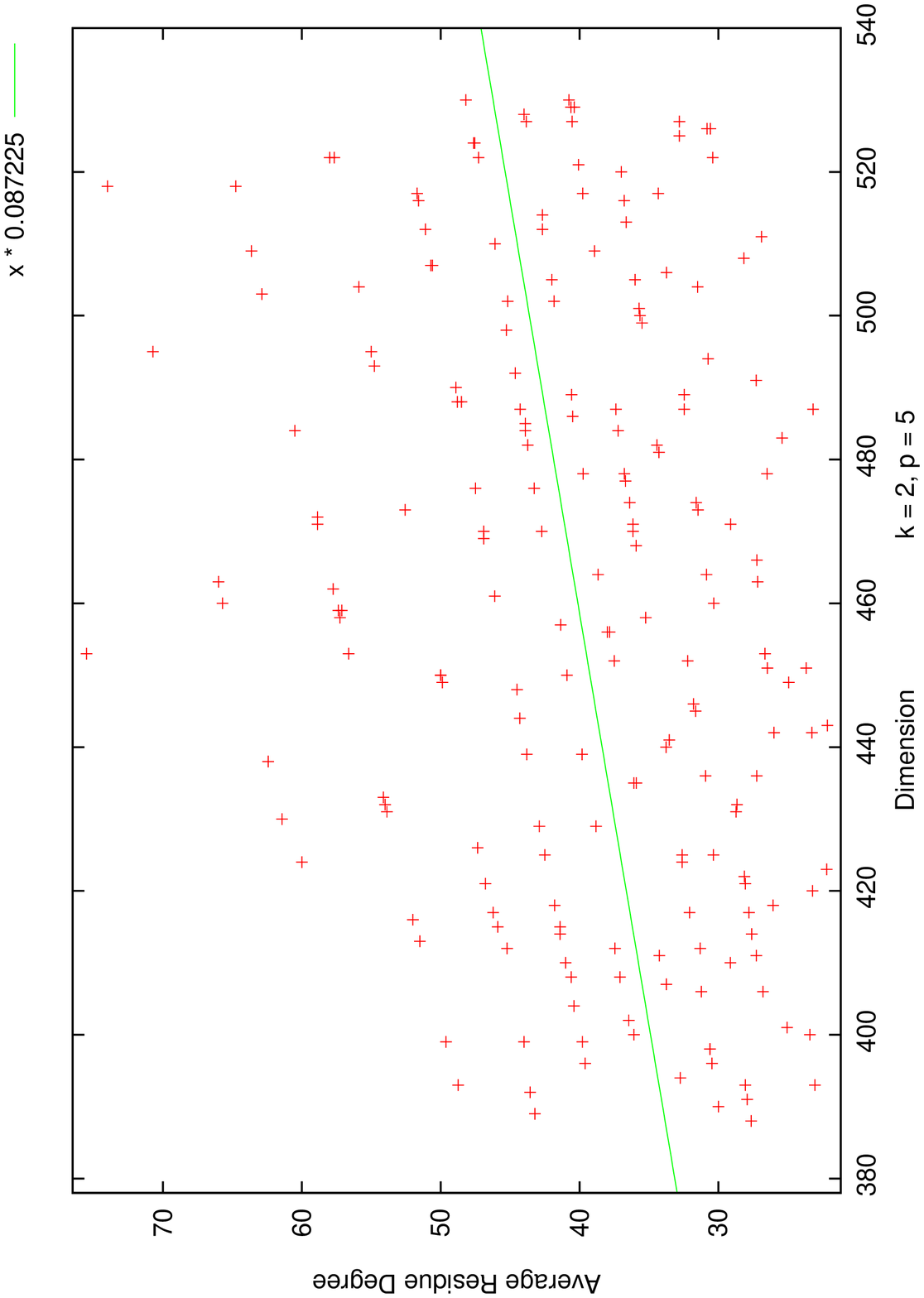} \\
\mplot{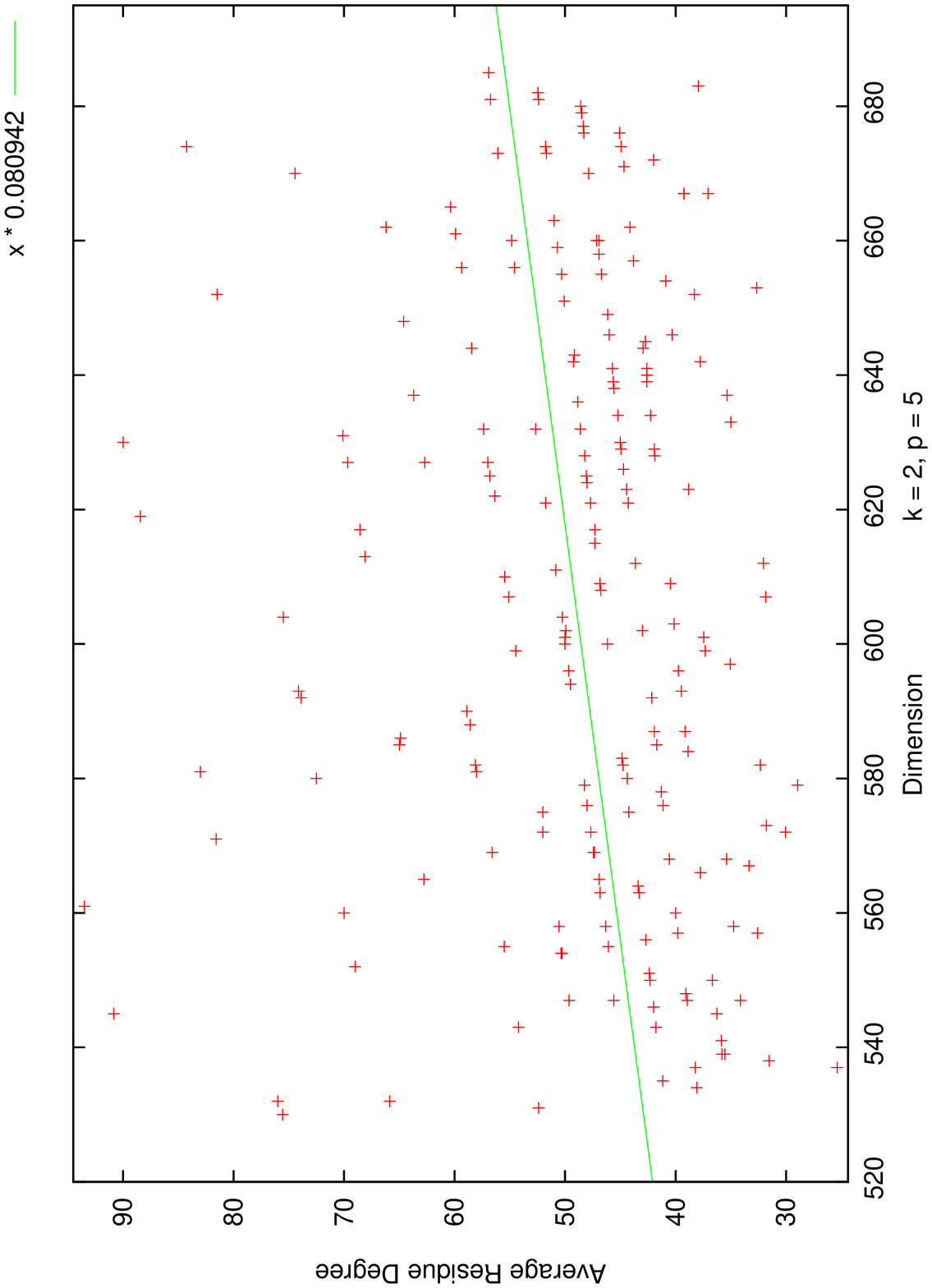} &
\mplot{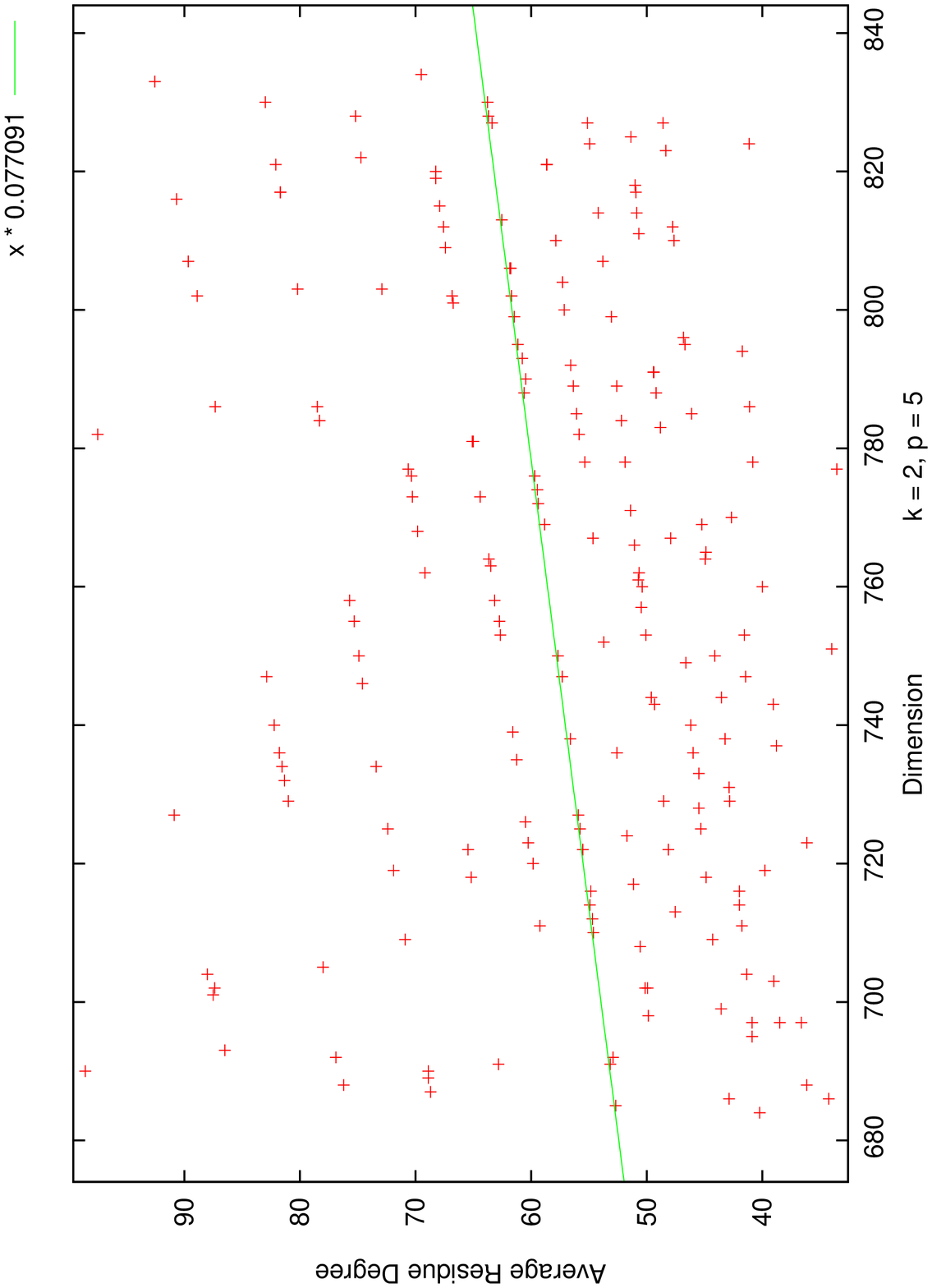} \\
\end{longtable}

We observe that the slope of the best fitting line goes slightly down with the dimension. This
strongly suggests that taking a straight line does not seem to be quite correct.
We also made logarithmic plots, which we do not reproduce here; they seemed to suggest to us
that a behaviour of the form $b_{N,k}^{(p)} \sim \const \cdot d(N)^\alpha$ is not quite correct either
(the best choice of~$\alpha$ seems to be close to~$1$ in accordance with the previous discussion).

These computations suggest the following question.

\begin{question}\label{question:b}
Fix a prime~$p$ and an even weight~$k \ge 2$.
Let $b(N) := b_{k,N}^{(p)}$ and $d(N) := \dim_{\Fbar_p} S_k(N;\Fbar_p)$.
Do there exist constants $C_1,C_2$ and $0<\alpha \le \beta < 1$ such that the inequality
$$ C_1 + \alpha \frac{d(N)}{\log(d(N))} \le b(N) \le C_2 + \beta \cdot d(N)$$
holds?
\end{question}

We remark that if $a_{N,k}^{(p)}$ behaves like $d(N)$, as suggested by Question~\ref{question:a-odd},
then Question~\ref{question:b} is equivalent to asking that $\#\Spec(\Tbar_k(N))$ does not
grow faster than a constant times $\log(d(N))$.

The phenomenon that for odd primes~$p$ most of the dots in the diagrams seem to lie on
or very close to certain distinguished lines through the origin is natural
in view of Question~\ref{question:a-odd}: the slope of the line on or close to which a dot lies
is just $\frac{1}{\#\Spec(\Tbar_k(N))}$.

\subsection{Maximum Residue Degree}\label{sec:c}

Now we turn our attention to the {\em maximum residue degree},
which we define for given level~$N$, weight~$k$ and prime~$p$ as
$$ c_{N,k}^{(p)} = \max \{[\FF_\fm:\FF_p] \;\mid\; \fm \in \Spec(\Tbar_k(N))\}.$$

We made computations for weight~$2$ and all primes~$p$ less than~$100$, where $N$
runs through the same ranges as previously. We again plot the dimension $d(N)$ on the $x$-axis
and the function~$c_{N,k}^{(p)}$ on the $y$-axis and the green line is again the best fitting
function $\alpha \cdot d(N)$. This time we believe that this might be the right function to take.
Here is again a selection of the plots that we obtained.

\noindent\begin{longtable}{cc}
\mplot{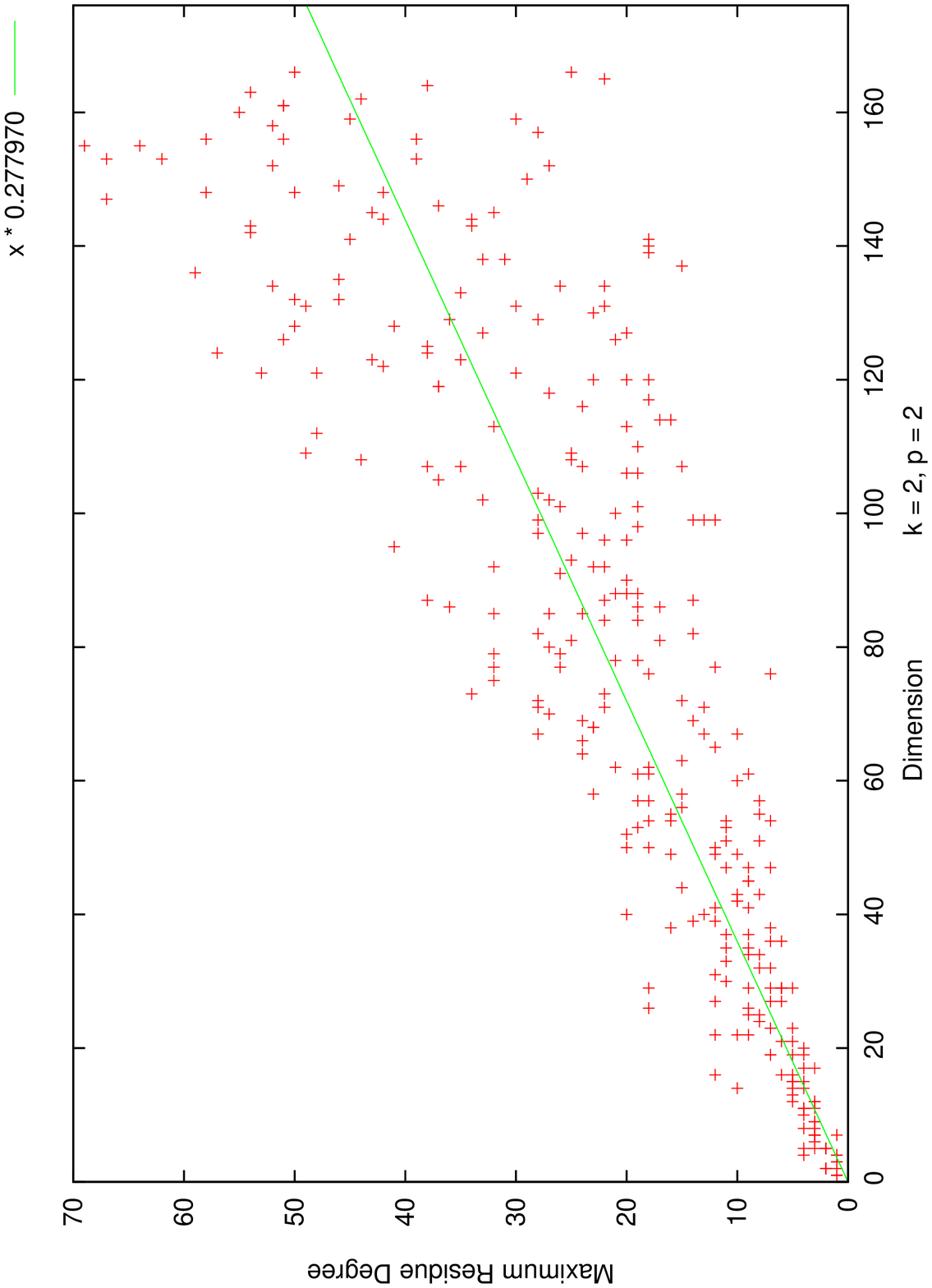} &
\mplot{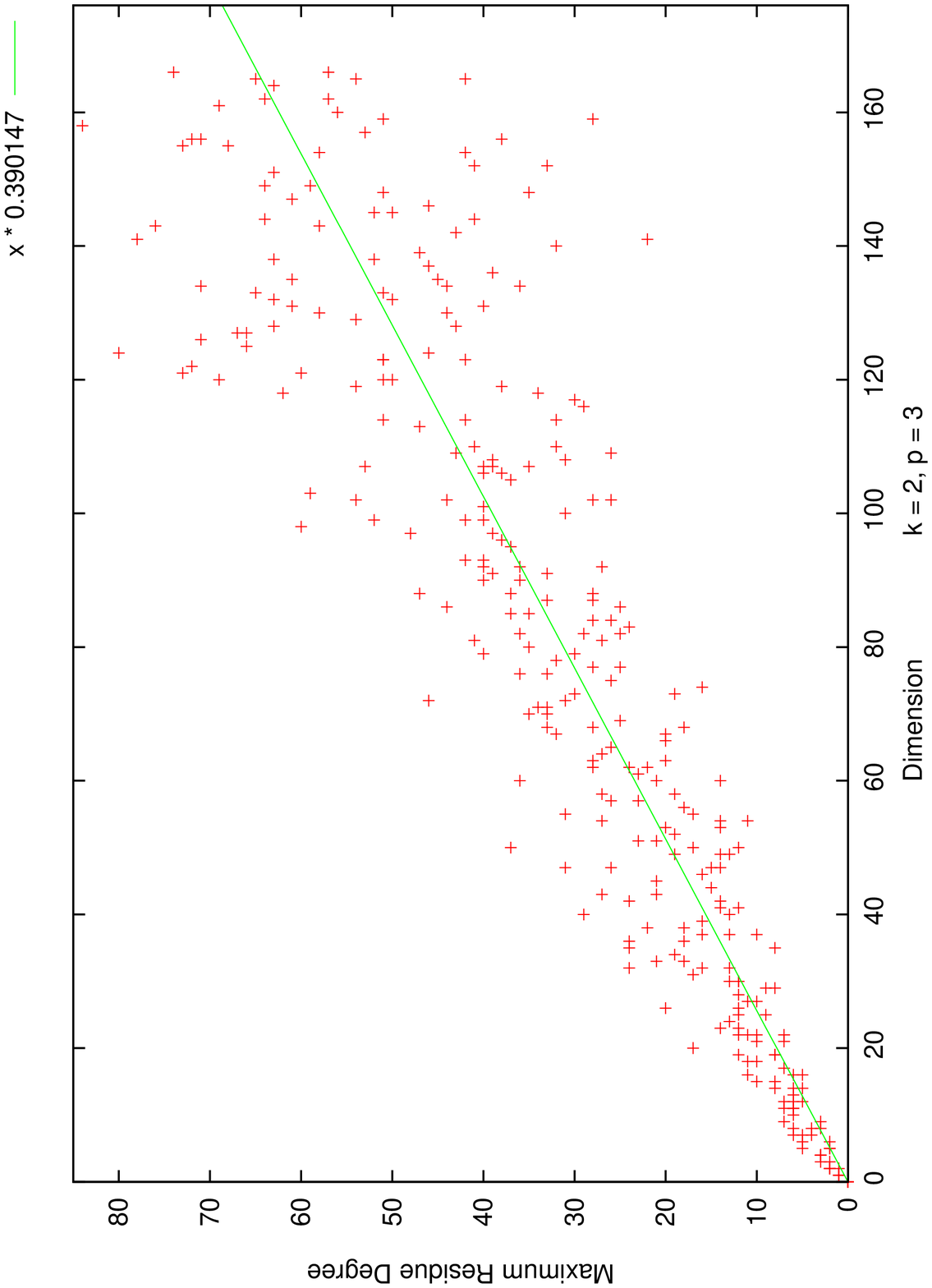} \\
\mplot{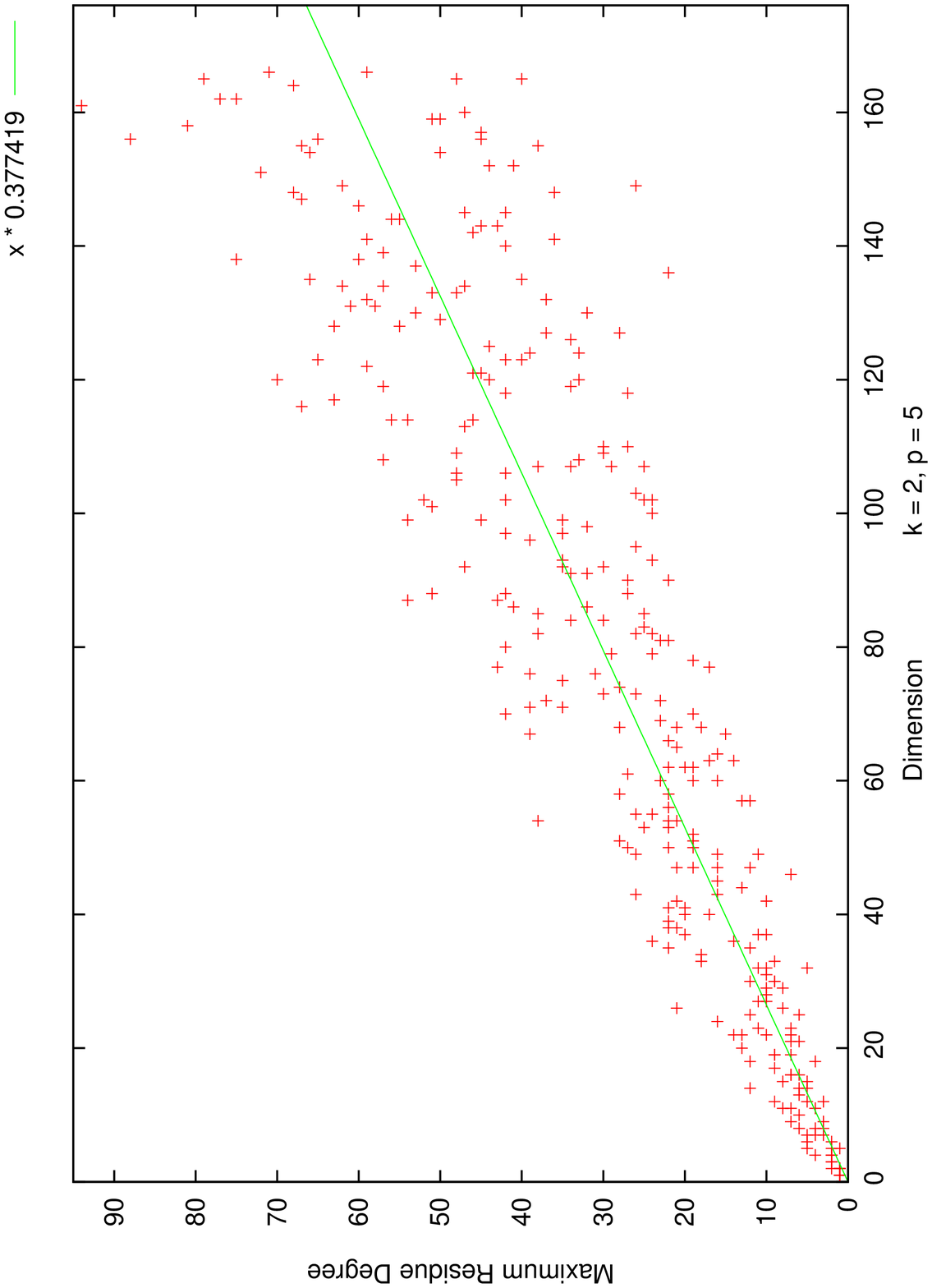} &
\mplot{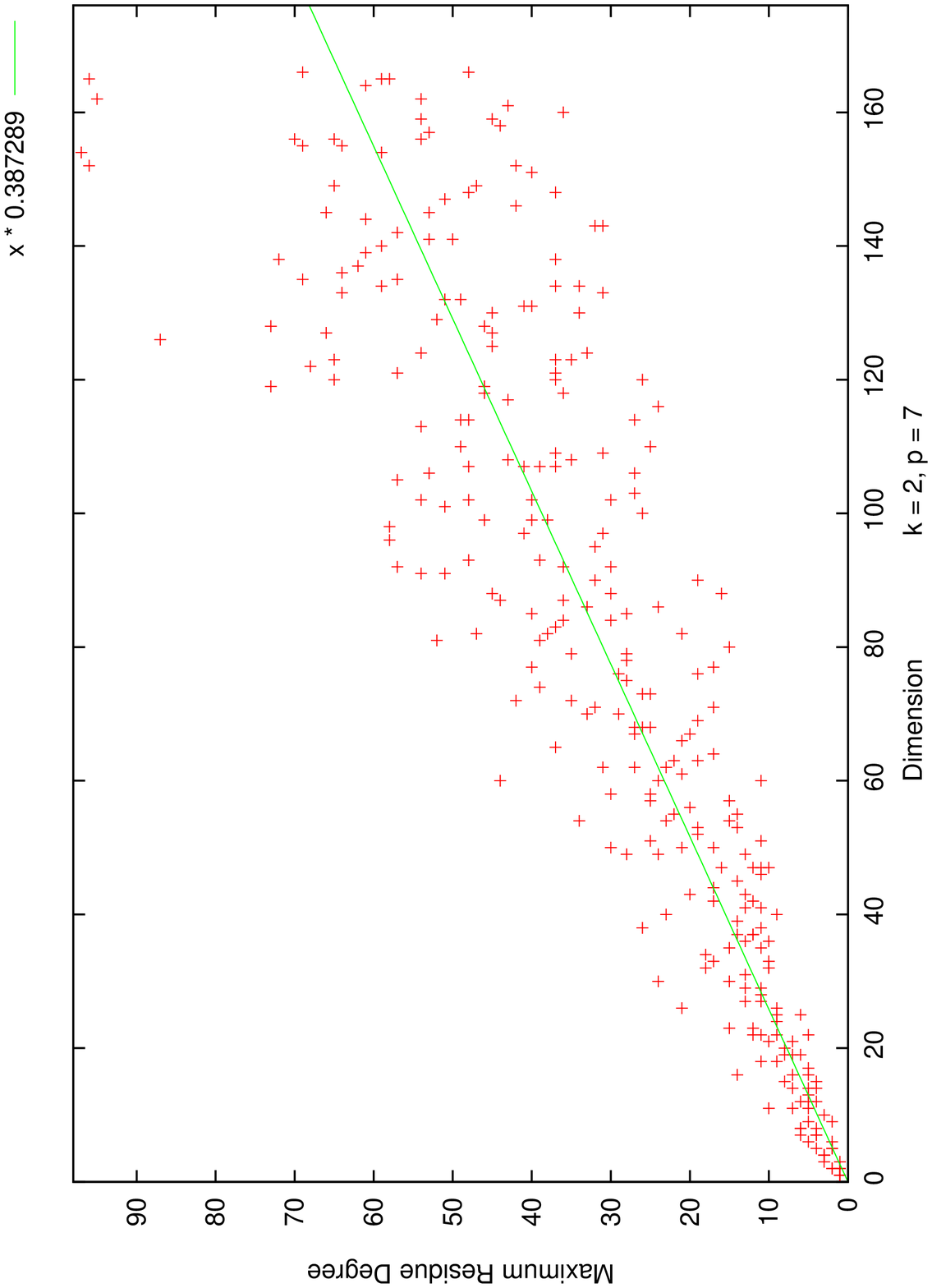} \\
\mplot{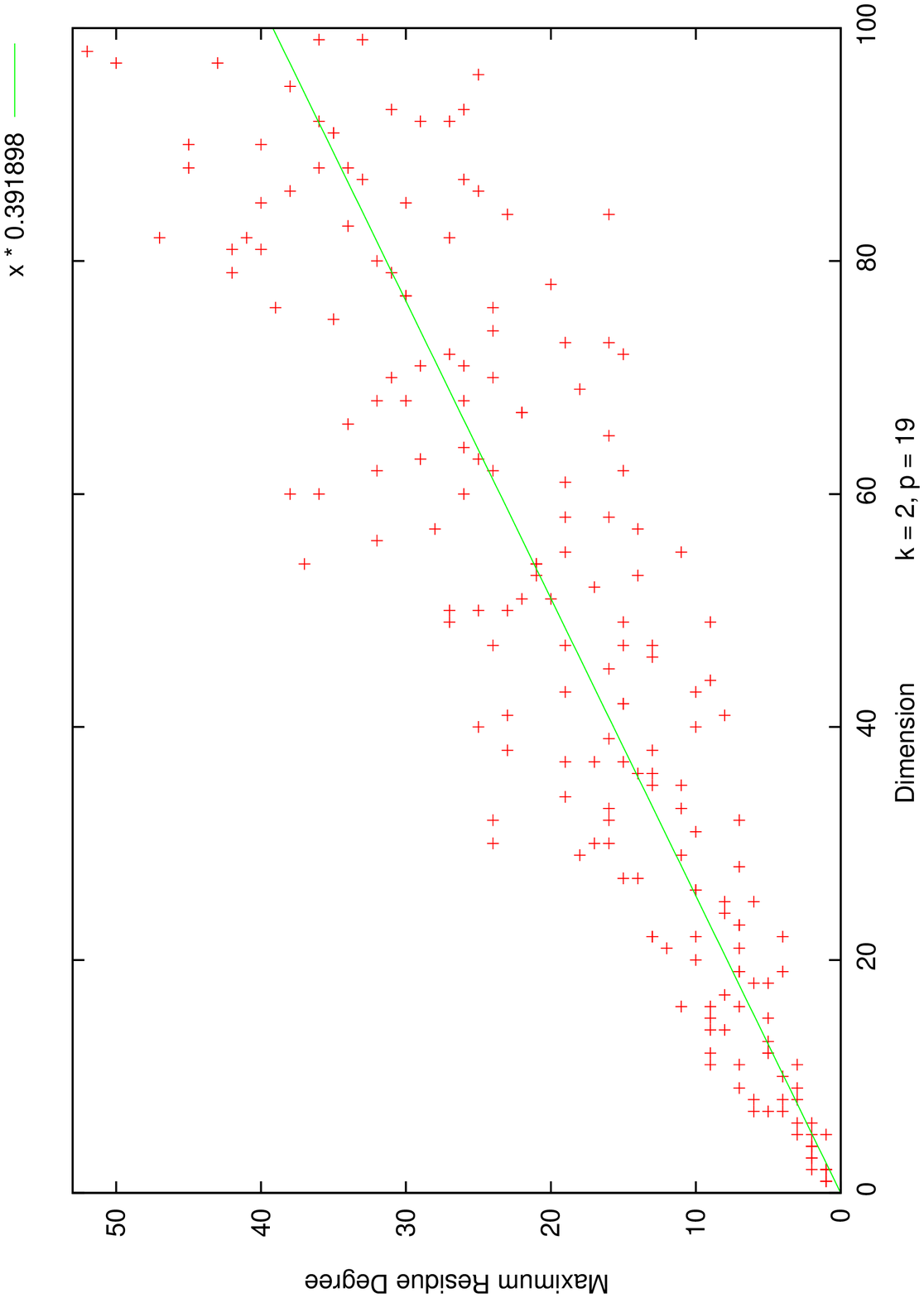}&
\mplot{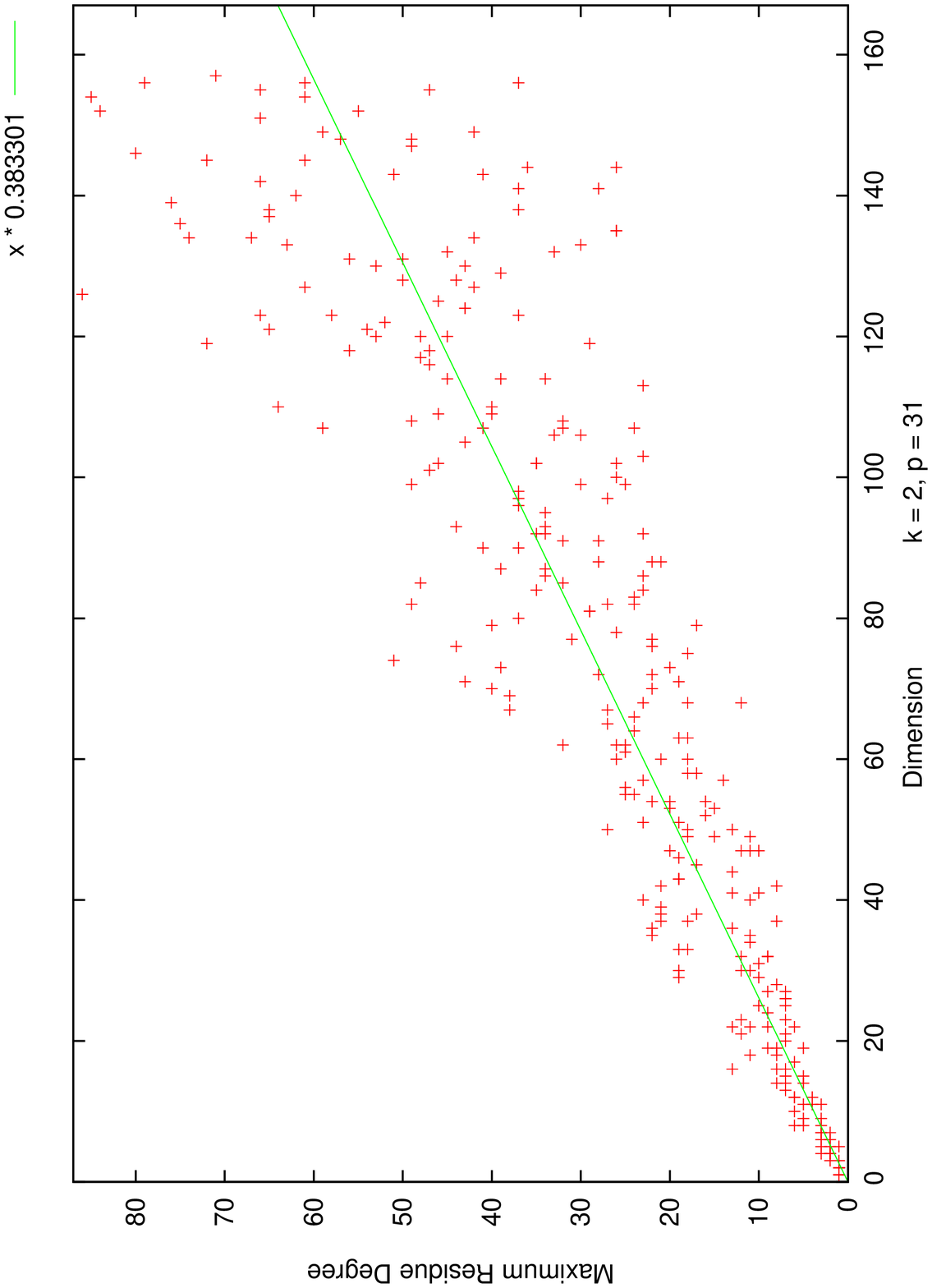} \\
\mplot{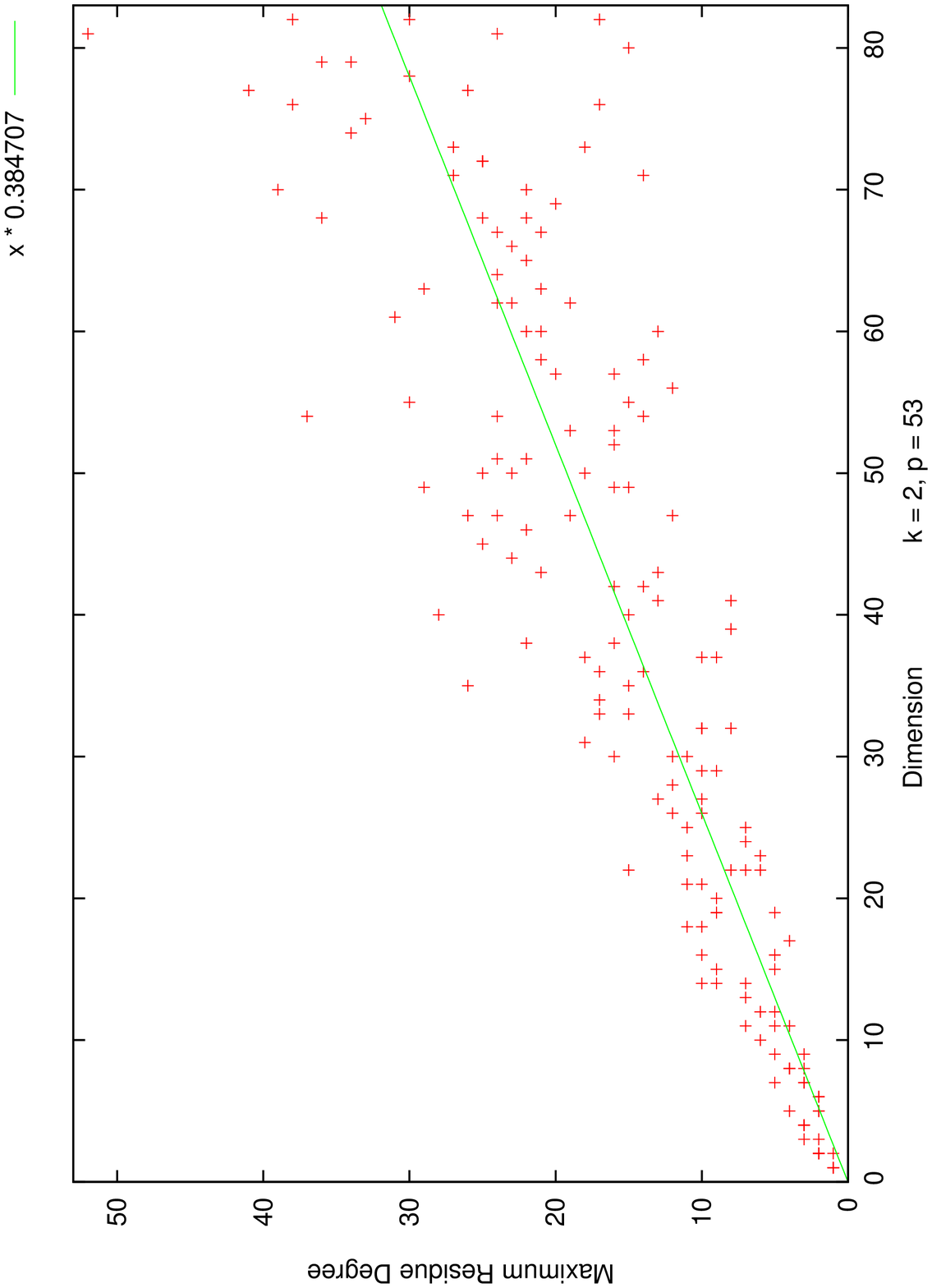}&
\mplot{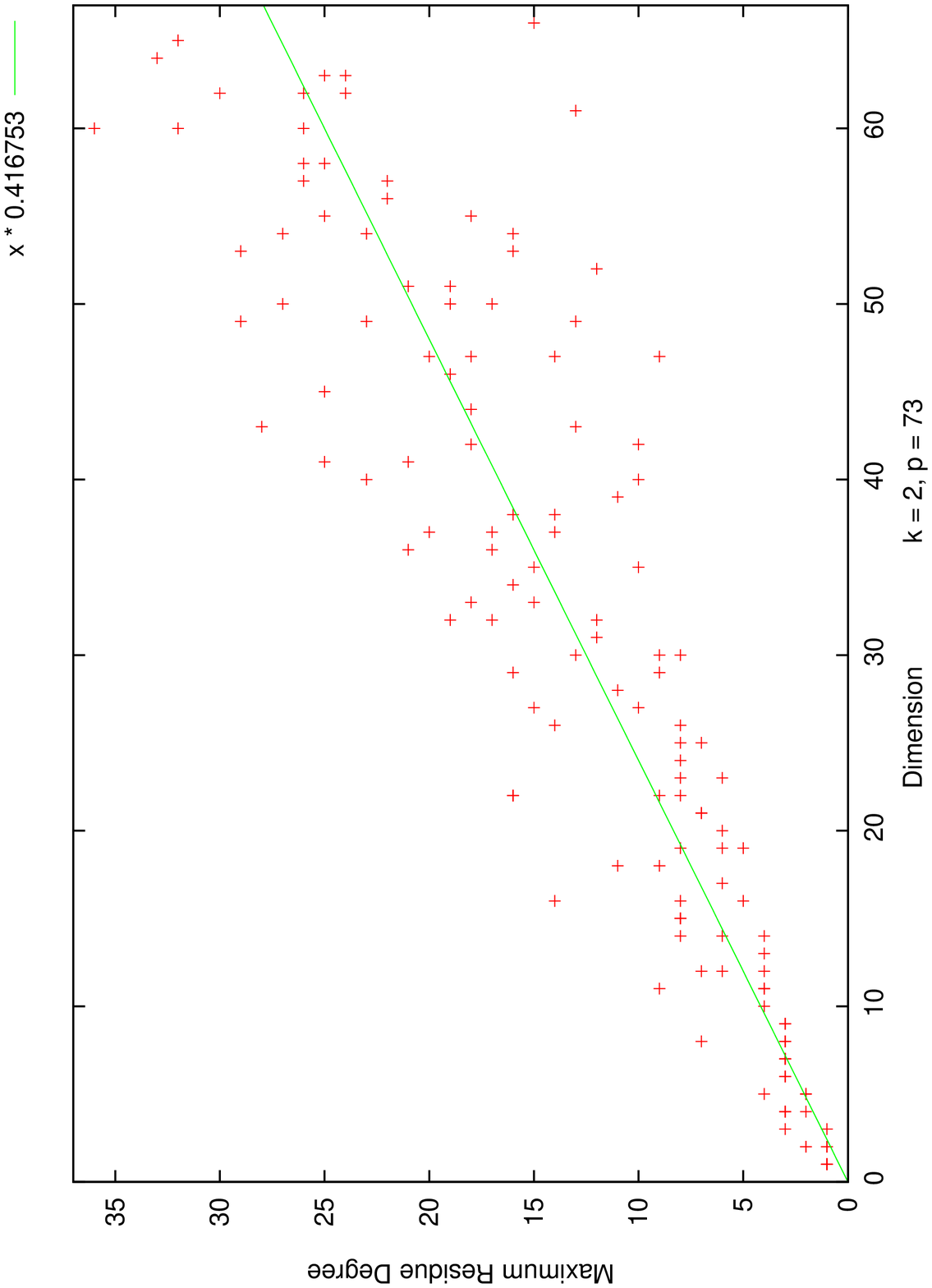} \\
\mplot{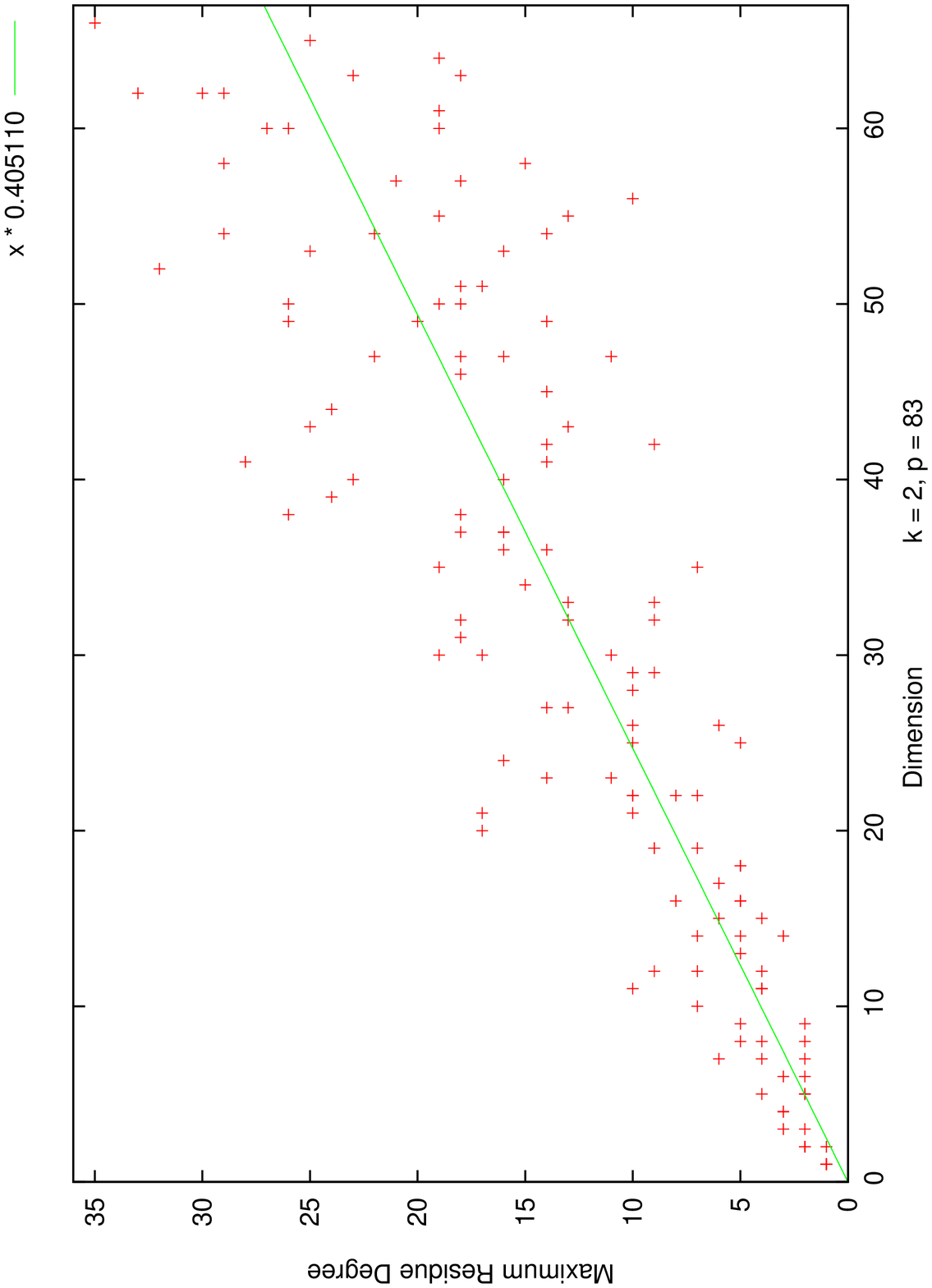} &
\mplot{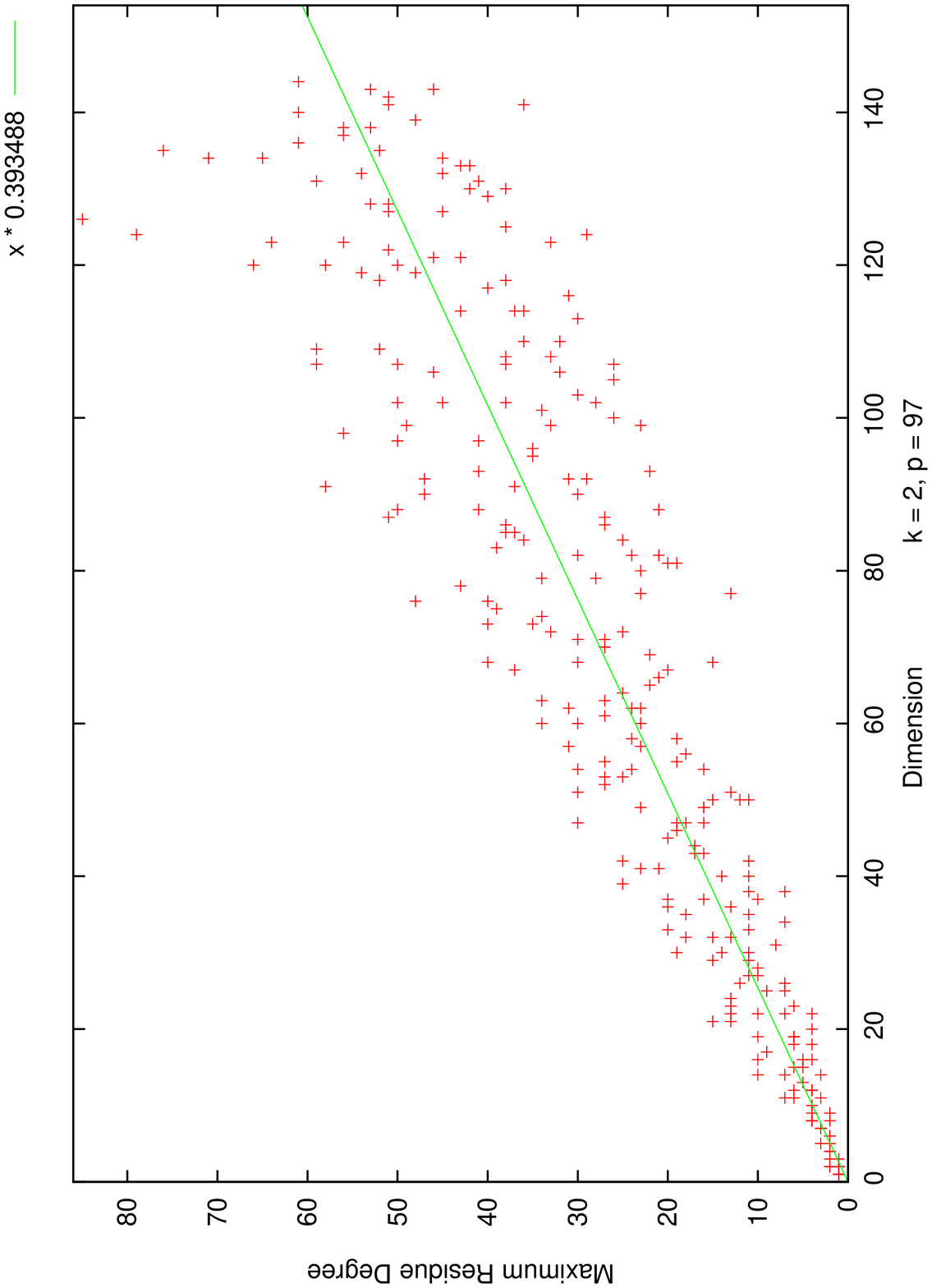} \\
\end{longtable}

Here are again two examples for weight~$4$.
\noindent\begin{longtable}{cc}
\mplot{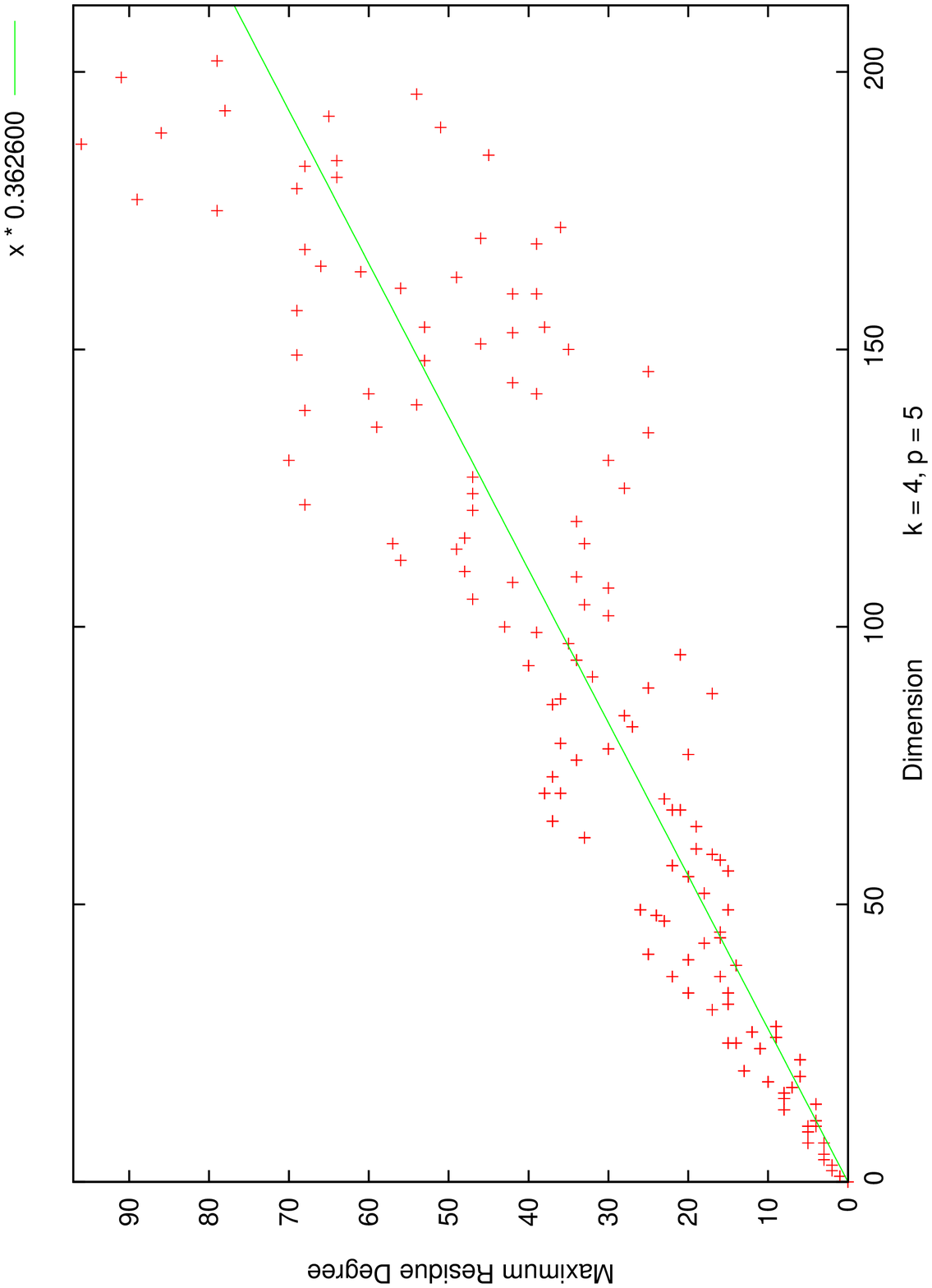} &
\mplot{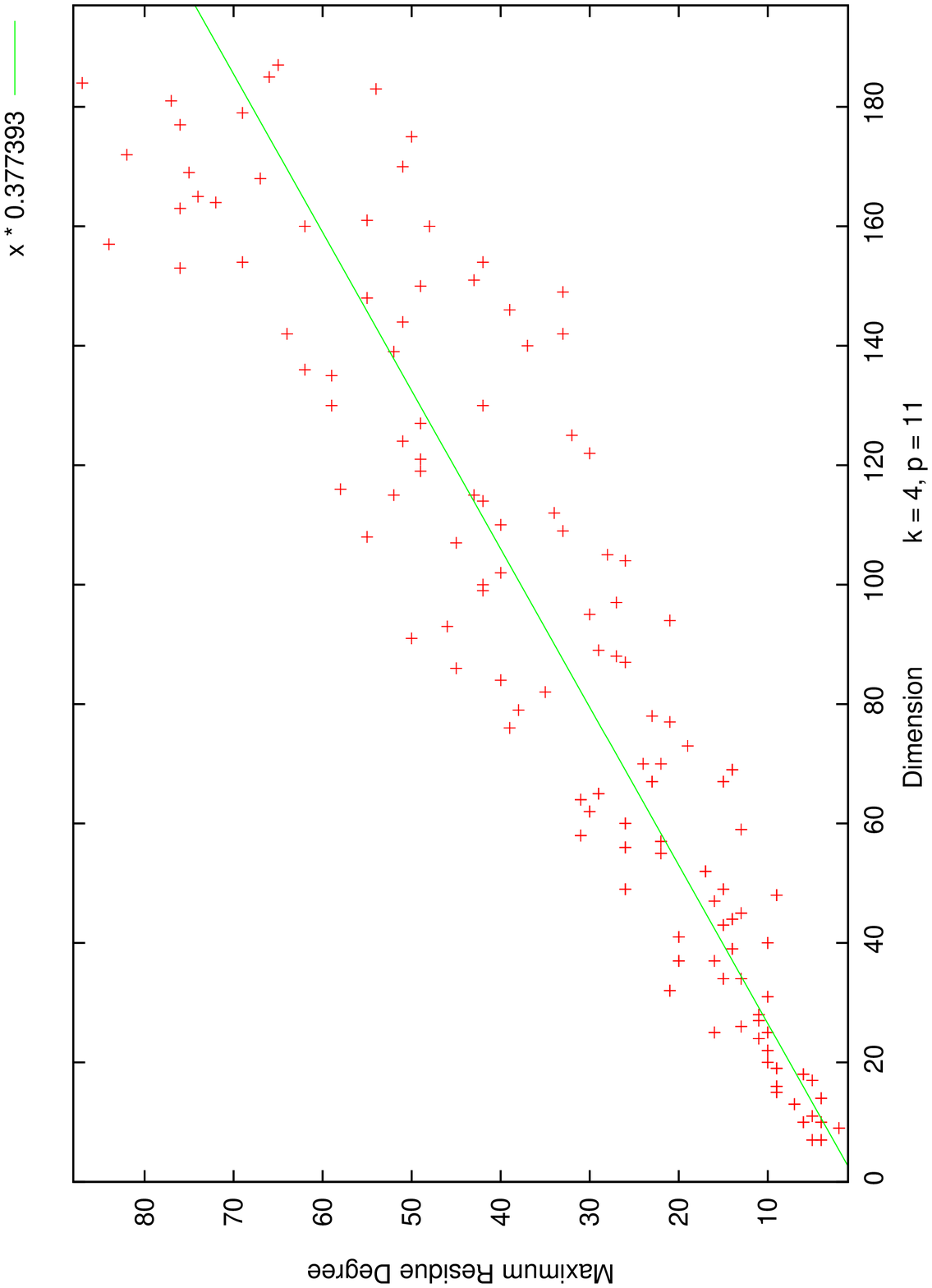} \\
\end{longtable}

The data certainly suggest that the maximum residue degree grows
with the dimension.
It is remarkable to see that the slopes of the best fitting lines all seem to be very close
to each other -- with the single exception of the case $p=2$, which might be caused
by the same yet unknown phenomenon as earlier.
Also in this case we conducted a closer analysis for the primes~$2$, $3$ and~$5$.
For $p=2$ we used all primes in different intervals up to~$12000$ and obtained
these plots:

\noindent\begin{longtable}{cc}
\mplot{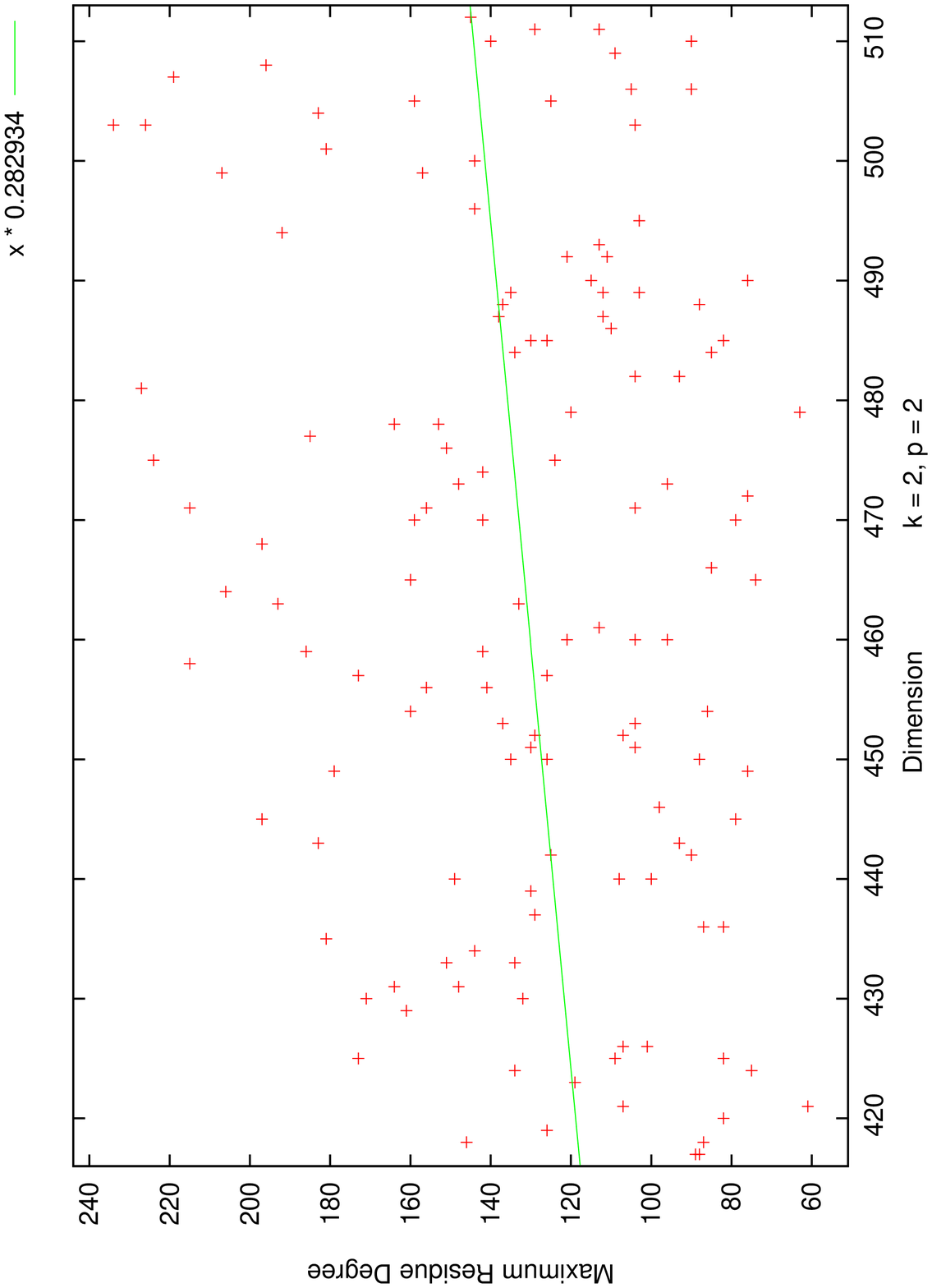} &
\mplot{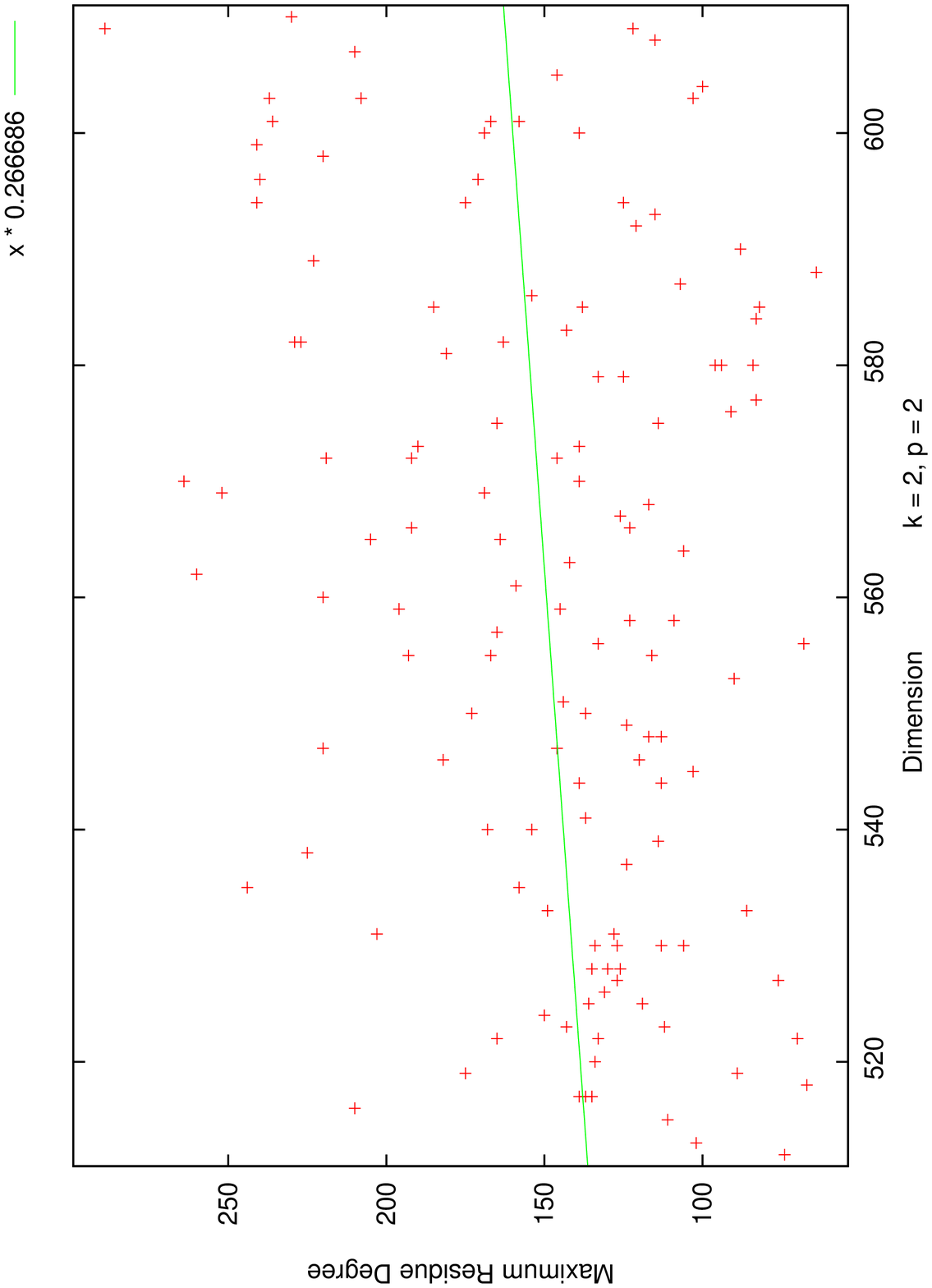} \\
\mplot{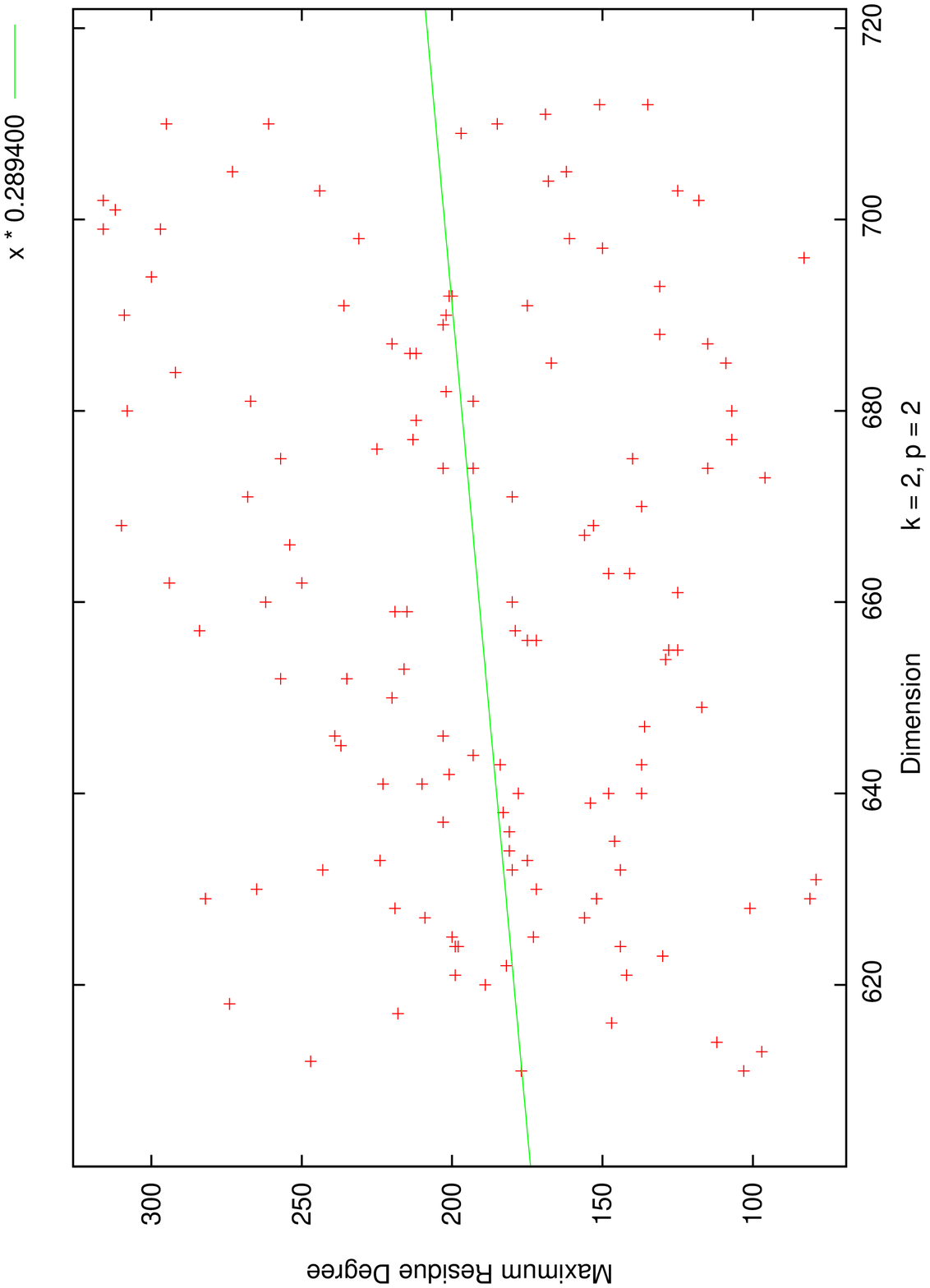} &
\mplot{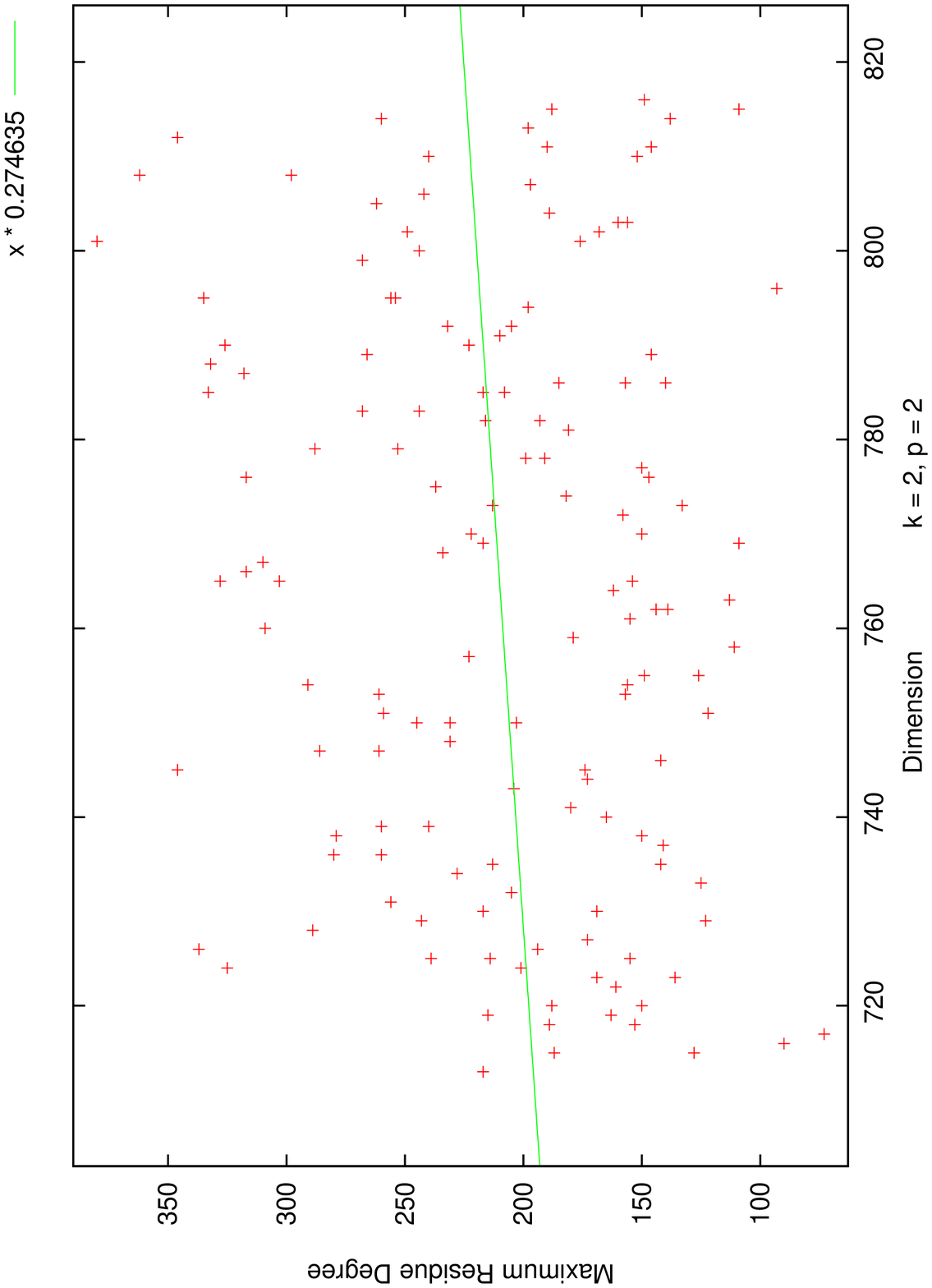} \\
\mplot{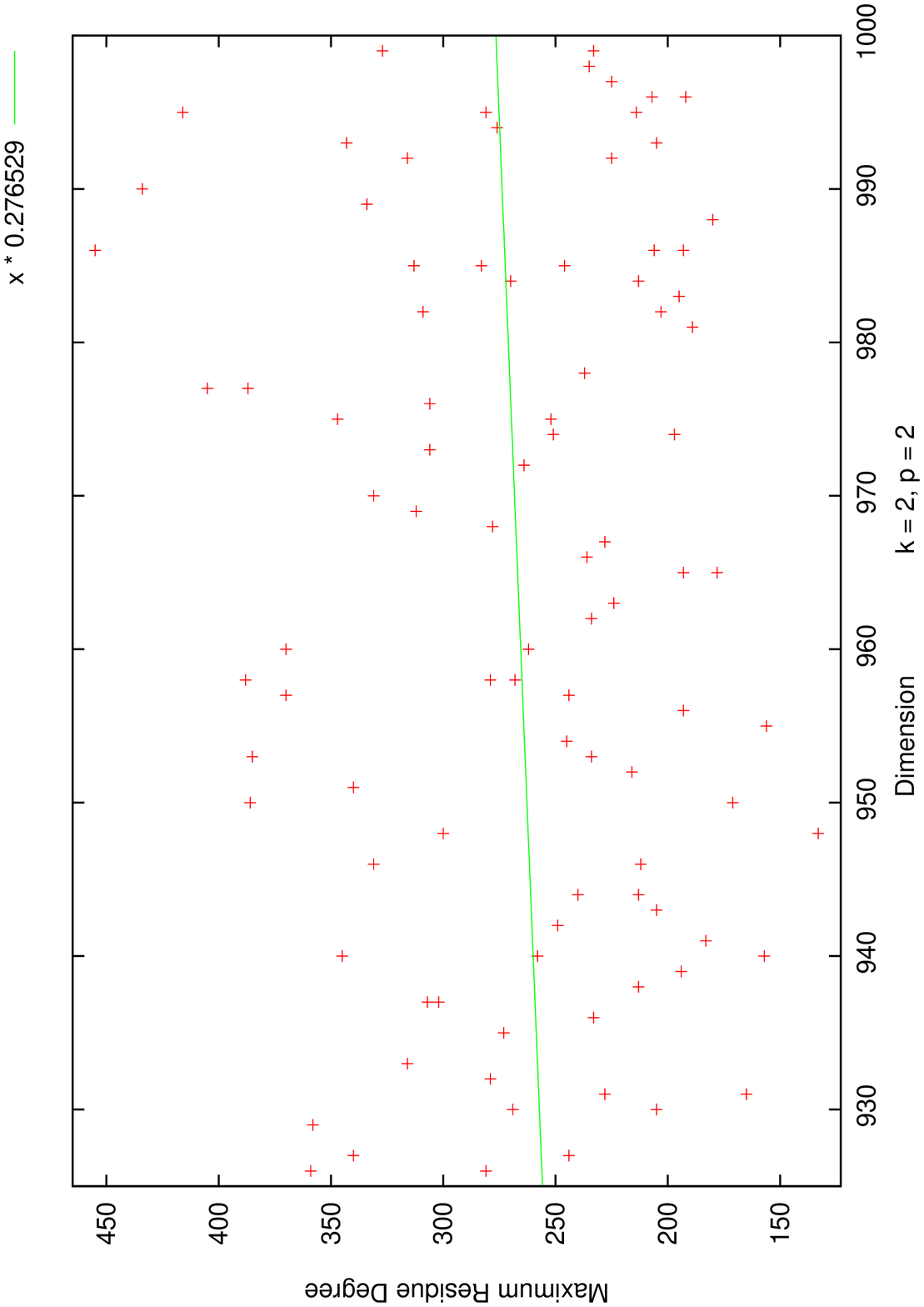} \\
\end{longtable}

Here are the plots for $p=3,5$ and the primes between $3000$ and $10009$ subdivided into four intervals.

\noindent\begin{longtable}{cc}
\mplot{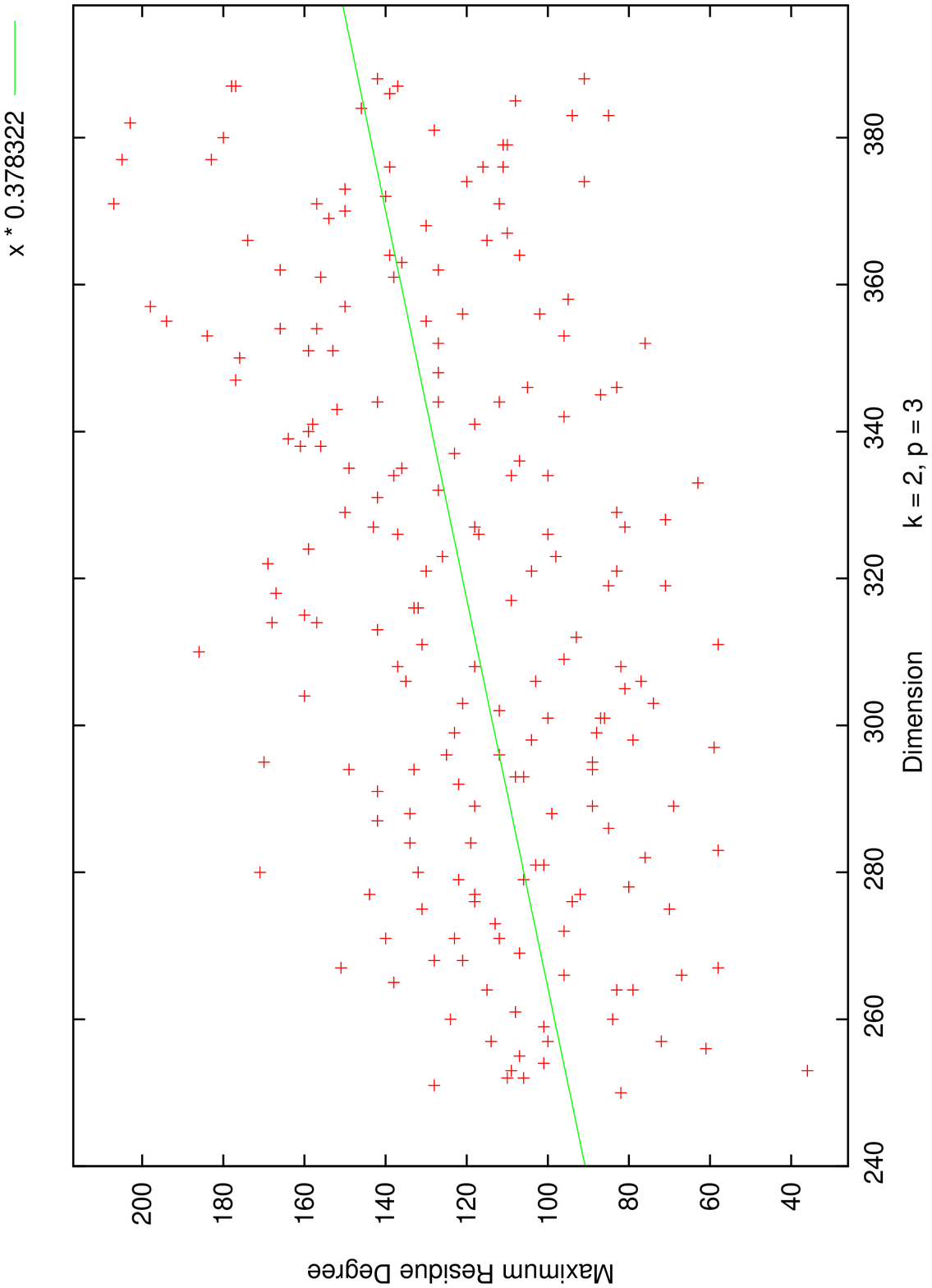} &
\mplot{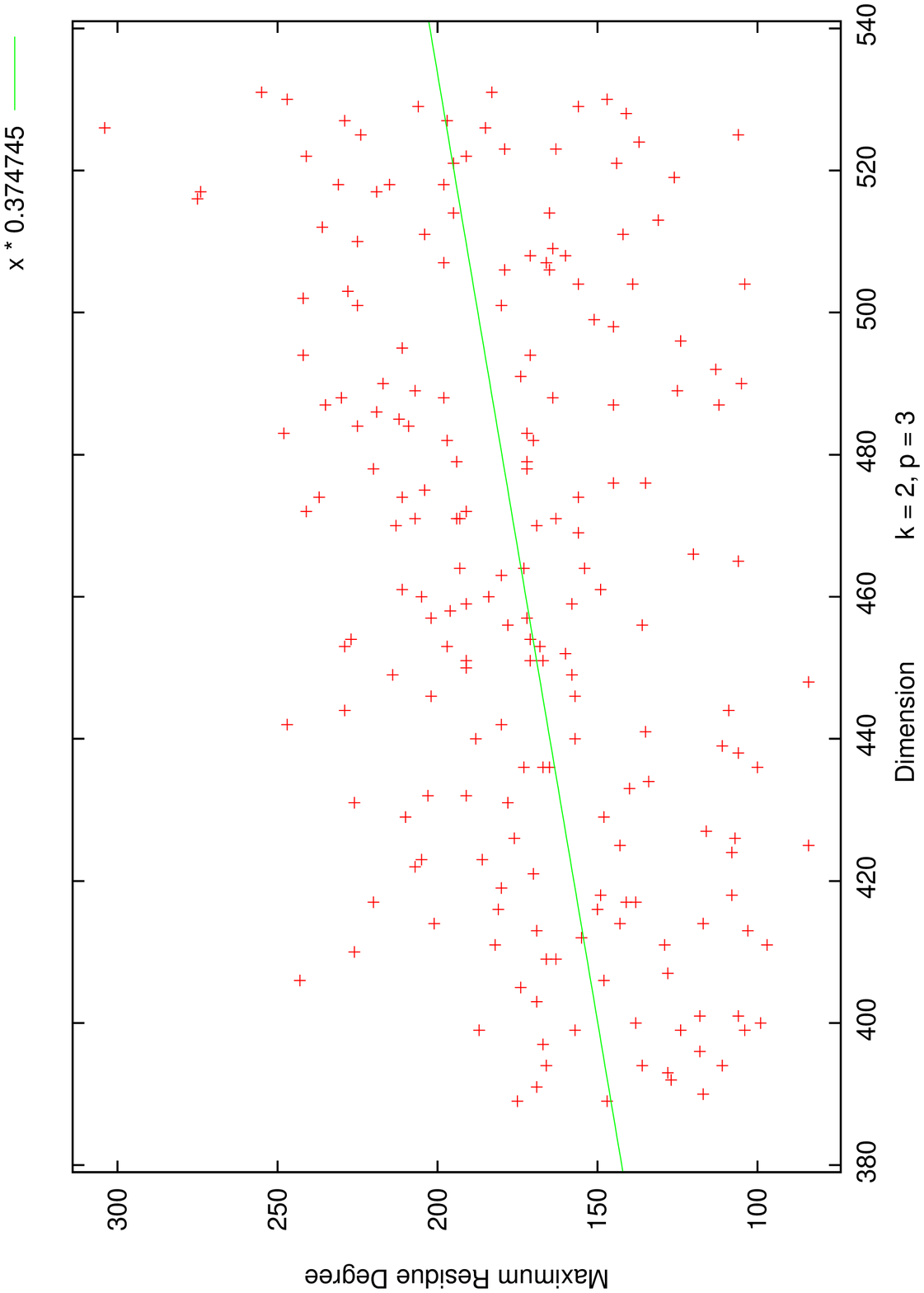} \\
\mplot{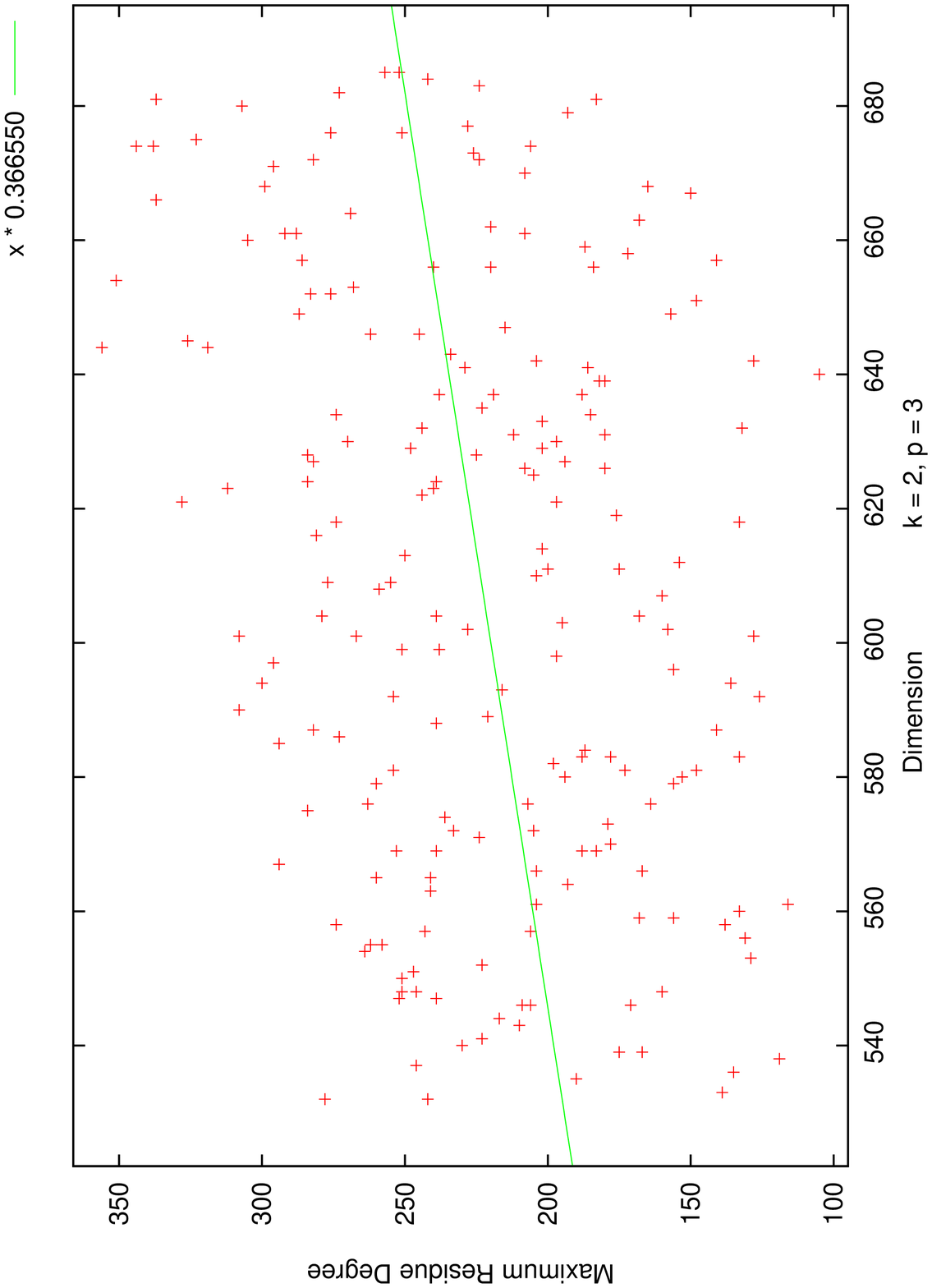} &
\mplot{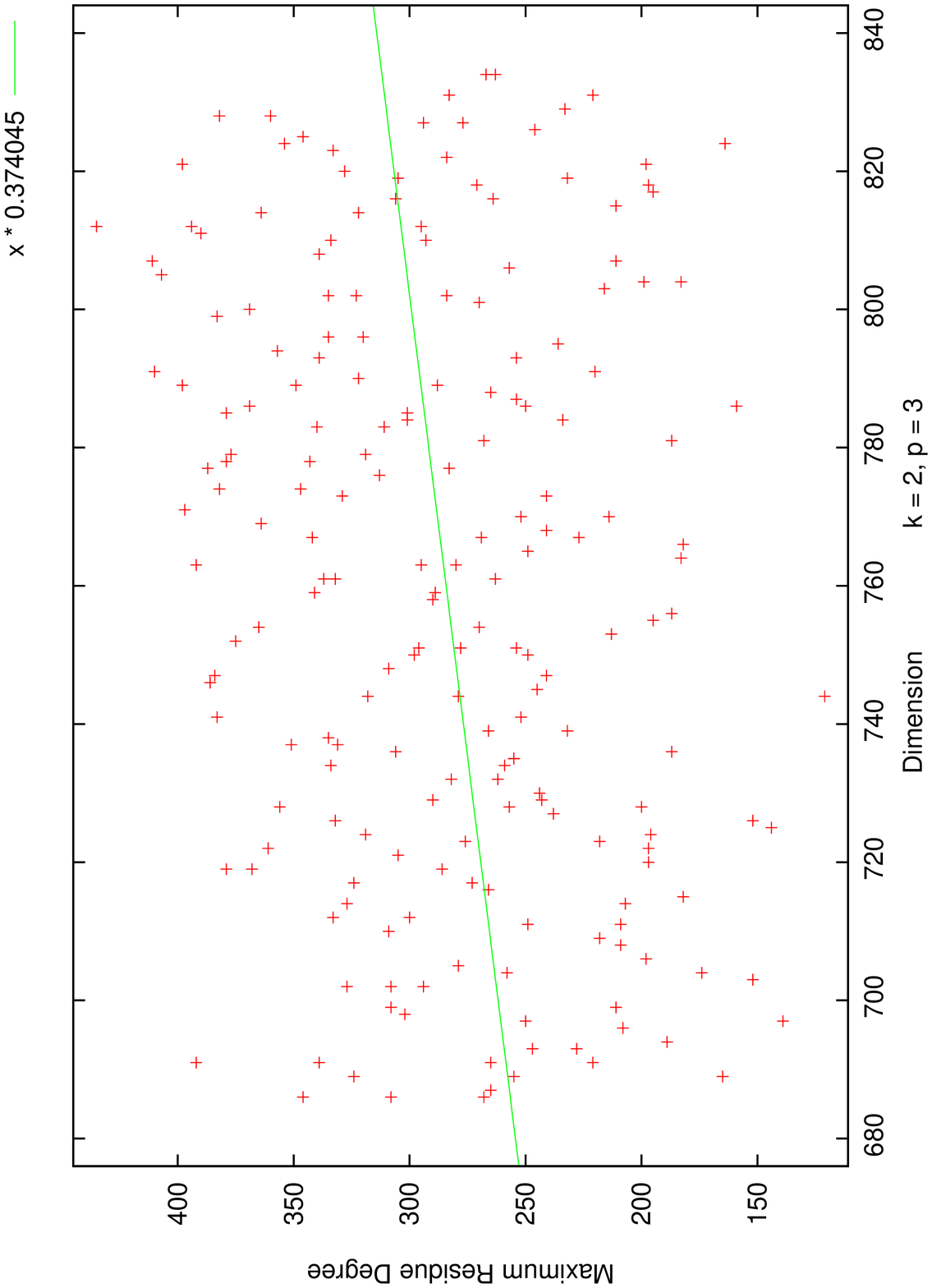} \\
\mplot{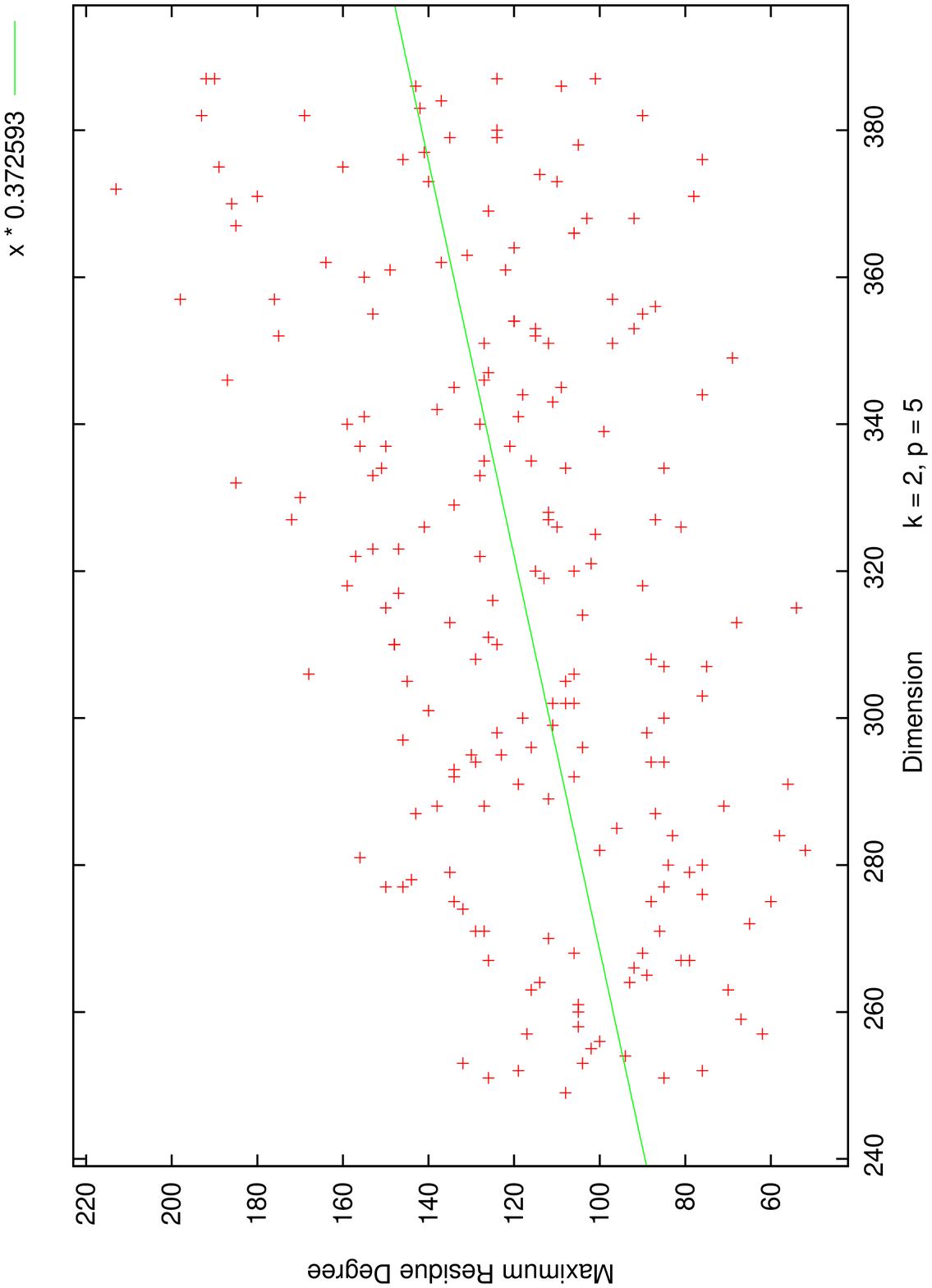} &
\mplot{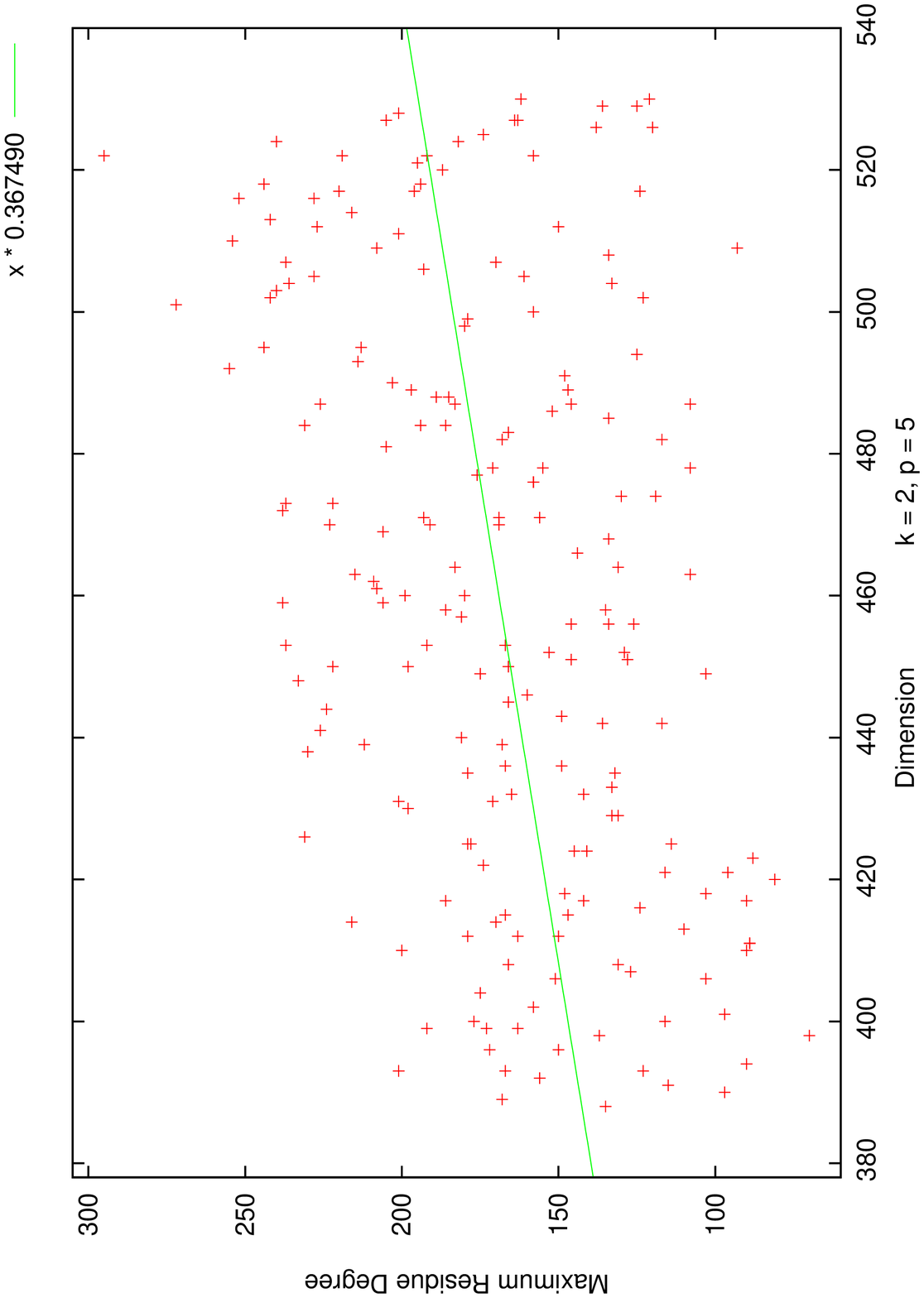} \\
\mplot{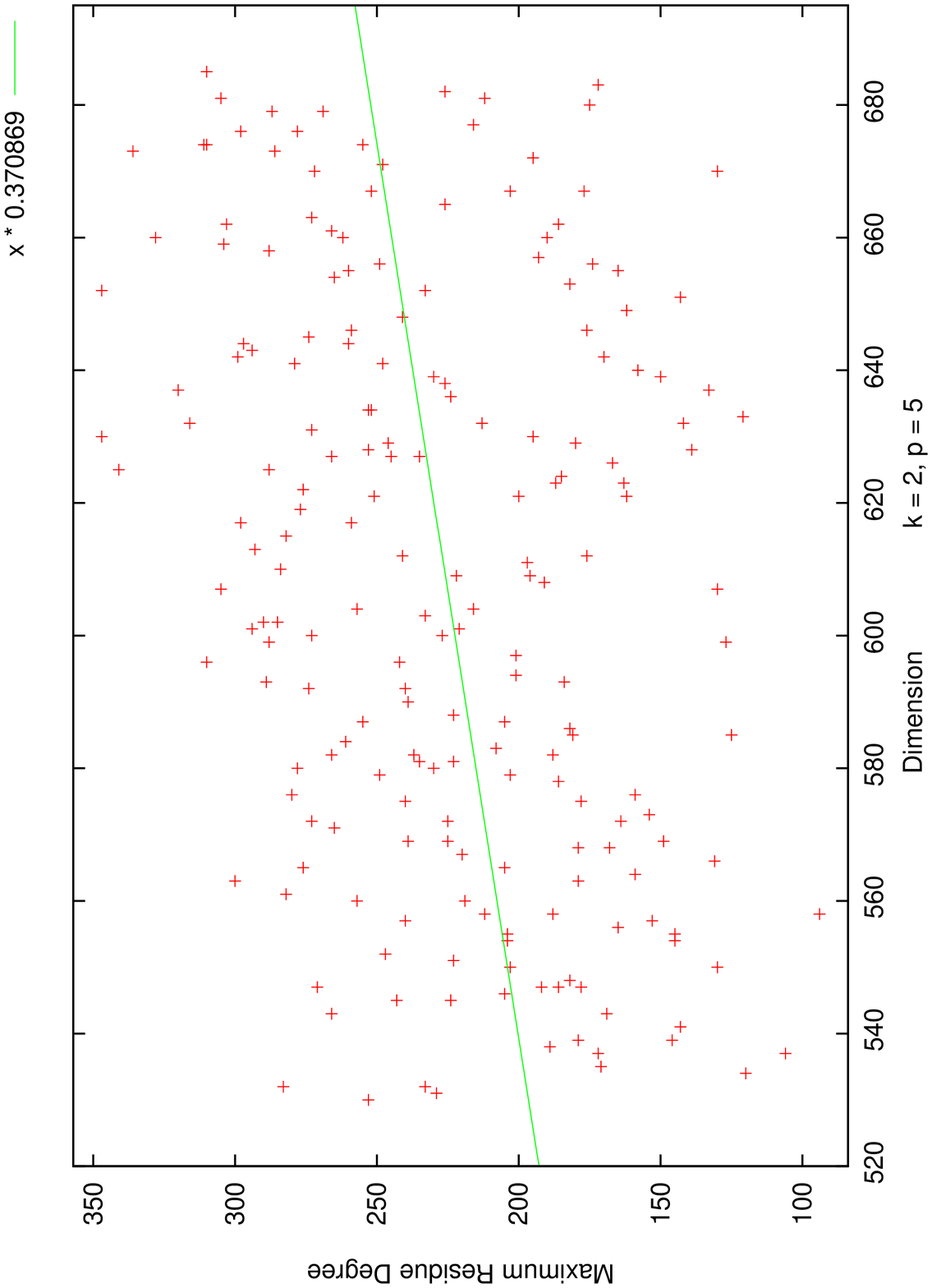} &
\mplot{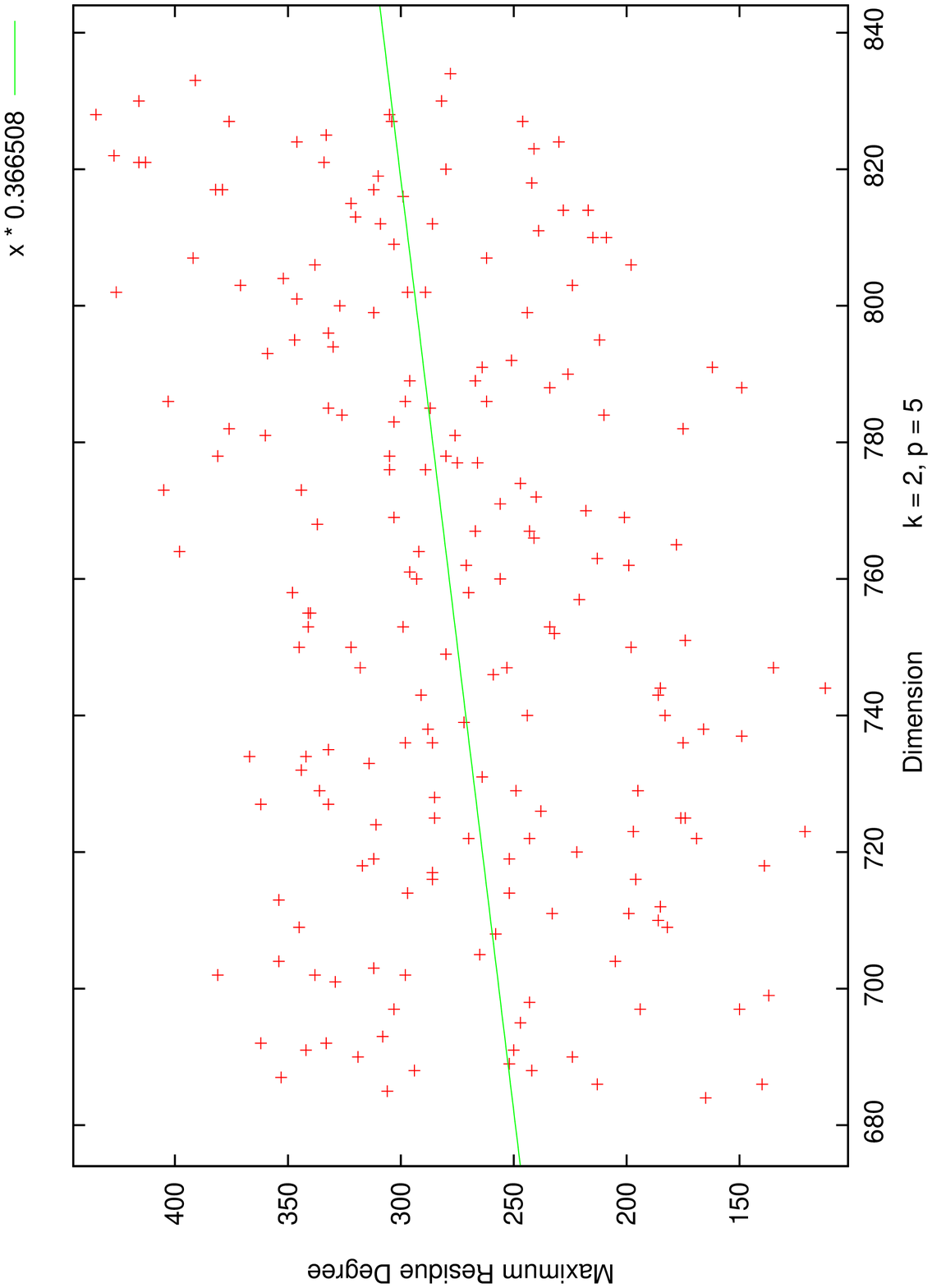} \\
\end{longtable}

We observe that, although the best fitting lines were computed using different intervals,
their slopes are very close to each other.
These computations suggest the following question.

\begin{question}\label{question:c}
Fix a prime~$p$ and an even weight~$k \ge 2$.
Let $c(N) := c_{k,N}^{(p)}$ and $d(N) := \dim_{\Fbar_p} S_k(N;\Fbar_p)$.
Do there exist constants $C_1,C_2$ and $0<\alpha \le \beta < 1$ such that the inequality
$$ C_1 + \alpha \cdot d(N) \le c(N) \le C_2 + \beta \cdot d(N)$$
holds?
\end{question}

\end{document}